%% file: main.tex
\documentclass{article}

\usepackage{microtype}
\usepackage{graphicx}
\usepackage{subfigure}
\usepackage{booktabs} 

\usepackage{hyperref}

\usepackage{icml2021_arxiv}

\usepackage{xr}
\makeatletter
\newcommand*{\addFileDependency}[1]{
  \typeout{(#1)}
  \@addtofilelist{#1}
  \IfFileExists{#1}{}{\typeout{No file #1.}}
}
\makeatother
 
\newcommand*{\myexternaldocument}[1]{%
    \externaldocument{#1}%
    \addFileDependency{#1.tex}%
    \addFileDependency{#1.aux}%
}

\myexternaldocument{SI_ICML}

\input{custom_settings}

\icmltitlerunning{Efficient Discretizations of Optimal Transport}

\begin{document}
\twocolumn[
\icmltitle{Efficient Discretizations of Optimal Transport}



\icmlsetsymbol{equal}{*}

\begin{icmlauthorlist}
\icmlauthor{Junqi Wang}{ru}
\icmlauthor{Pei Wang}{ru}
\icmlauthor{Patrick Shafto}{ru}
\end{icmlauthorlist}

\icmlaffiliation{ru}{Department of Mathematics \& Computer Science, Rutgers University, Newark, US}

\icmlcorrespondingauthor{Pei Wang}{peiwang@rutgers.edu}

\icmlkeywords{Machine Learning, ICML}

\vskip 0.3in
]



\printAffiliationsAndNotice{}  

\begin{abstract}
  \input{abstract}
\end{abstract}

\input{1_introduction}
\input{2_optimal_discretization}

\input{3_gradient}
\input{4_algorithm}

\input{5_parallel}

\input{6_complexity}
\input{7_conclusion_etc}

\appendix
\input{SI_1_proofs}

\input{SI_2_differentiation}
\input{SI_3_empirical}

\bibliography{references}
\bibliographystyle{icml2021}

\end{document}

%% file: custom_settings.tex
\usepackage{amsmath, amssymb, amsfonts, amsthm}
\usepackage{enumitem}
\usepackage{blkarray}
\usepackage{mathrsfs}  
\usepackage{mathtools}

\usepackage{tikz}
\usetikzlibrary{arrows}
\usetikzlibrary{calc}
\usetikzlibrary{math}

\makeatletter
\newtheorem*{rep@theorem}{\rep@title}
\newcommand{\newreptheorem}[2]{%
\newenvironment{rep#1}[1]{%
\def\rep@title{#2 \ref{##1}}%
\begin{rep@theorem}}%
{\end{rep@theorem}}}
\makeatother

\newtheorem{thm}{Theorem}[section]
\newreptheorem{thm}{Theorem}
\newtheorem*{thm*}{Theorem}

\newtheorem{prop}[thm]{Proposition}
\newreptheorem{prop}{Proposition}

\theoremstyle{definition}

\newtheorem*{definition*}{Definition}

\newtheorem*{example*}{Example}

\newtheorem*{framework*}{Unifying OT Framework}

\theoremstyle{remark}

\newtheoremstyle{smplain}
  {\topsep}   
  {\topsep}   
  {\itshape}  
  {0pt}       
  {\bfseries} 
  {.}         
  {5pt plus 1pt minus 1pt} 
  {\thmname{#1}\thmnumber{ {#2}}\thmnote{ (#3)}}          
  
\theoremstyle{smplain}
\newcounter{sdecorator}
\newcommand{\argmin}{\text{argmin }}
\newcommand{\argmax}{\text{argmax }}
\newcommand{\KL}{\text{KL}}
\newcommand{\dd}{\textrm{d}}

%% file: abstract.tex
Obtaining solutions to Optimal Transportation (OT) problems is 
typically intractable when the marginal spaces are continuous.
Recent research has focused on approximating continuous solutions with
discretization methods based on i.i.d. sampling, and has proven convergence
as the sample size increases. 
However, obtaining OT solutions with large sample sizes requires intensive computation effort, 
that can be prohibitive in practice. 
In this paper, we propose an algorithm for calculating 
discretizations with a given number of points for marginal 
distributions,
by minimizing the (entropy-regularized) Wasserstein distance,
and result in plans that are comparable to those obtained with much larger numbers of i.i.d. samples.
Moreover, a local version of such 
discretizations which is parallelizable for
large scale applications is proposed.
We prove bounds for our approximation and demonstrate performance on a wide range of problems. 

%% file: 1_introduction.tex
\section{Introduction}

%
%

Optimal transport is the problem of finding a coupling of probability distributions that minimizes cost \cite{kantorovich2006translocation}, and is a 
technique applied across various fields and literatures \cite{peyre2019computational,villani2008optimal}. 
Although many methods exist for obtaining optimal transference plans for distributions on discrete spaces, 
computing the plans is not generally possible for continuous spaces \cite{janati2020entropic}. 
Given the prevalence of continuous spaces in machine learning, this is a significant limitation for 
theoretical and practical applications. 

One strategy for approximating continuous OT plans is based on discrete approximation, e.g. via samples. 
Recent research has provided guarantees on the fidelity of discrete, sample-based approximations to continuous OT 
as $N\to \infty$ \cite{aude2016stochastic}. Specifically, by sampling large numbers of points $S_i$ from each marginal, 
one may compute discrete optimal transference plan on $S_1\times S_2$, with the cost matrix 
being derived from pointwise evaluation of the cost function on $S_1\times S_2$.

Even in the discrete case, obtaining minimal cost plans is computationally challenging. 
For example, Sinkhorn scaling, which computes an entropy-regularized approximation to OT plans, 
has a complexity that scales with $|S_1\times S_2|$ \cite{allen2017much}.
Though many comparable methods exist \cite{lin2019efficiency}, all have complexity that scales with the product of sample sizes,
and require construction of the cost matrix that also scales with $|S_1\times S_2|$. 

Building on previous sample-based approaches, we develop methods for optimizing small N approximations of OT plans. 
In Sec.~\ref{sec:optimal_discretizations}, we formulate the problem of fixed size approximation
and reduce it to discretization problems on marginals with theoretical guarantees.
In Sec.~\ref{sec:gradient}, the gradient
of entropy-regularized Wasserstein distance between a continuous distribution and its discretization 
is calculated.
In Sec.~\ref{sec:algorithm_design}, we present a stochastic gradient descent algorithm that is based on optimization of the location and weight
of the points with empirical test results.
Sec.~\ref{sec:practical_methods} introduces parellelizable algorithm via decompositions of the marginal spaces,
which can exploit the intrinsic geometry. 
In Sec.~\ref{sec:complexity}, 
we analyze time and space complexity, and provide comparison with the naive sampling.

%% file: 2_optimal_discretization.tex
\section{Optimal Discretizations}
\label{sec:optimal_discretizations}

Let $(X, d_X)$, $(Y,d_Y)$ be compact Polish spaces (complete separable metric spaces), 
$\mu\in\mathcal{P}(X)$, $\nu\in\mathcal{P}(Y)$ be probability distributions
on their Borel-algebras, 
and $c: X\times Y \to \mathbb{R}$ be a cost function. 
Denote the set of all joint probability measures on $X\times Y$ with marginals $\mu$ and $\nu$ by $\Pi(\mu, \nu)$.
For a given cost function $c$, the \textit{optimal transference (OT) plan} between $\mu$ and $\nu$ is defined as \cite{kantorovich2006translocation}:
$$\gamma(\mu,\nu) \coloneqq \argmin_{\pi \in \Pi(\mu, \nu)} \langle c, \pi \rangle$$
where $\langle c, \pi \rangle \coloneqq \int_{X\times Y} c(x, y) \dd \pi(x, y)$. 

When $X = Y$, if cost $c(x, y) = d^k_X(x, y)$, 
$W_{k}(\mu, \nu) =  \langle c, \gamma(\mu,\nu) \rangle^{1/k}$
defines the \textit{$k$-Wasserstein distance} between $\mu$ and $\nu$ for $k\geq 1$.
Here $d^k_X(x, y)$ is the $k$-th power of the metric $d_X$ on $X$. 


\textit{Entropy regularized optimal transport} (EOT) \cite{cuturi2013sinkhorn,aude2016stochastic} was introduced 
to ease the computational burden of obtaining OT plan:
\begin{equation}\label{eq:eot}
    \gamma_{\lambda}(\mu,\nu) \coloneqq \argmin_{\pi \in \Pi(\mu, \nu)} \langle c, \pi \rangle + \lambda \KL(\pi||\mu \otimes \nu)
\end{equation}

where $\lambda >0$ is a regularization parameter and 
$\KL(\pi||\mu \otimes \nu) \coloneqq \int \log(\frac{\text{d} \pi}{ \text{d}\mu \otimes\text{d}\nu}) \dd \pi$ 
is the Kullback-Leibler divergence. 
EOT smooths the classical OT into a convex problem. 
Hence, given $(\mu, \nu, c)$, there exists an unique solution to~\eqref{eq:eot}.
After $X$ and $Y$ are discretized, 
EOT  problems can be solved by Sinkhorn iteration \cite{sinkhorn1967concerning}, which can be easily parallelized. 
However, for large-scale discrete spaces, the computational cost can still be unfeasible \cite{allen2017much}.
Even worse, to even apply Sinkhorn iteration, one must have a cost matrix, which itself can be a non-trivial computational burden to obtain; 
in some cases, for example where the cost is derived from a probability model \cite{wang2020sequential}, 
it may require intractable computations \cite{tran2017variational,overstall2020bayesian}. 

We propose to efficiently estimate continuous OT with a fixed size discrete approximation.
In more details,
let $m, n\in \mathbb{Z}^*$, $\mu_m\in\mathcal{P}(X)$,
$\nu_n\in\mathcal{P}(Y)$ be discrete approximation of $\mu, \nu$ respectively with $\mu_m=\sum_{i=1}^mw_i\delta_{x_i}$
and $\nu_n=\sum_{j=1}^nu_j\delta_{y_j}$, $x_i\in X$, $y_j\in Y$ and $w_i,u_j\in \mathbb{R}^+$.
Then the EOT plan $\gamma_\lambda(\mu,\nu)\in\Pi(\mu,\nu)$ for OT problem $(\mu, \nu, c)$, 
can be approximated by the EOT plan $\gamma_\lambda(\mu_m,\nu_n)\in\Pi(\mu_m,\nu_n)$ for OT problem $(\mu_m, \nu_n, c)$.

To obtain the best estimation the OT problem of $(\mu, \nu, c)$ 
with fixed size $m, n\in \mathbb{Z}^*$, 
we introduce the objective: 
\begin{align}\label{eq:obj}
\Omega _{k,\rho} (\mu_m, \nu_n) = & W_{k}^k(\mu, \mu_m)+W_{k}^k(\nu,\nu_n) \nonumber \\
& +\rho W_{k}^k(\gamma_{\lambda}(\mu,\nu), \gamma_{\lambda}(\mu_m,\nu_n)),     
\end{align}

$W_{k}^k(\mu, \mu_m)$, $W_{k}^k(\nu,\nu_n)$ and $W_{k}^k(\gamma_{\lambda}(\mu,\nu), \gamma_{\lambda}(\mu_m,\nu_n))$
are the $k$-th power of $k$-Wasserstein distances between $\mu, \nu, \gamma_{\lambda}(\mu,\nu)$ and their approximations
$\mu_m, \nu_n, \gamma_{\lambda}(\mu_m,\nu_n)$. 
The hyperparameter $\rho>0$ balances between the estimation accuracy over marginals and that of the transference plan.

To properly compute $W_{k}^k(\gamma_{\lambda}(\mu,\nu), \gamma_{\lambda}(\mu_m,\nu_n))$, 
a metric $d_{X\times Y}$ on $X\times Y$ is needed.
Without lost of generality, we assume there exists a constant $A$ such that 
\begin{align}\label{eq:d_xy}
    \max \{d^k_{X}(x_1, x_2), &d^k_{Y}(y_1, y_2)\} \stackrel{(i)}{\leq} d^k_{X\times Y}((x_1, y_1), (x_2, y_2)) \nonumber \\
    &\stackrel{(ii)}{\leq}A (d^k_{X}(x_1, x_2) + d^k_{Y}(y_1, y_2)).
\end{align}

For instance, \eqref{eq:d_xy} holds when $d_{X\times Y}$ is the $p$-product 
metric for $1\leq p\leq \infty$.
In particular, when $p=k=1$, the equality holds at $(ii)$ with $A=1$;
when $d_{X\times Y}$ is the $\infty$-product metric, equality holds at $(i)$.

To efficiently compute \eqref{eq:obj},
three Wasserstein distances are estimated by their entropy regularized approximations \cite{luise2018differential}:
\begin{align} \label{eq:obj_approx}
\Omega _{k,\zeta,\rho} (\mu_m, \nu_n) = & W_{k,\zeta}^k(\mu, \mu_m)+W_{k,\zeta}^k(\nu,\nu_n) \nonumber \\
& +\rho W_{k,\zeta}^k(\gamma_{\lambda}(\mu,\nu), \gamma_{\lambda}(\mu_m,\nu_n)).  
\end{align}

For instance, $W_{k}^k(\mu, \mu_m) = \langle d^k_X, \gamma(\mu,\mu_m) \rangle^{1/k}$ is estimated by
$W_{k,\zeta}^k(\mu, \mu_m) = \langle d^k_X, \gamma_{\zeta}(\mu,\mu_m) \rangle^{1/k}$ with regularizer $\zeta$.
Here $\gamma_{\zeta}(\mu,\mu_m) $ is computed by optimizing:
$$\widehat{W}_{k,\zeta}^k(\mu, \mu_m)\! =\! \langle d^k_X, \gamma_{\zeta}(\mu,\mu_m) \rangle + \lambda \KL(\gamma_{\zeta}(\mu,\mu_m)||\mu \otimes \mu_m).$$

A major drawback of optimizing $\Omega _{k,\zeta,\rho} (\mu_m, \nu_n)$ is evaluating 
the Wasserstein distance between $\gamma_{\lambda}(\mu,\nu)$ and $ \gamma_{\lambda}(\mu_m,\nu_n)$. 
In fact, calculating $\gamma_{\lambda}(\mu,\nu)$ is intractable, which is the original motivation to compute $\gamma_{\lambda}(\mu_m,\nu_n)$.
To overcome this difficulty, utilizing the dual formulation of the entropy regularized OT, we show that \footnote{Proofs are included in the Supplementary Material.}:

\begin{prop}\label{prop:3w_ineq} 
When $X$ and $Y$ are two compact spaces and $c$ is $\mathcal{C}^{\infty}$, there exists a constant
$C_1 \in \mathbb{R}^{+}$ such that 
\begin{align}
    \max\{ & W_{k}^k(\mu, \mu_m), \! W_{k}^k(\nu,\nu_n)\}\!\leq\! W_{k,\zeta}^k\!(\gamma_{\lambda}(\mu,\nu),\! \gamma_{\lambda}(\mu_m,\nu_n)) \nonumber\\
     & \leq C_1[W_{k,\zeta}^k(\mu, \mu_m)+W_{k,\zeta}^k(\nu,\nu_n)].
\end{align}
\end{prop}

Proposition~\ref{prop:3w_ineq} indicates that $\Omega _{k,\zeta,\rho} (\mu_m, \nu_n)$ can be approximated by $W_{k,\zeta}^k(\mu, \mu_m)+W_{k,\zeta}^k(\nu,\nu_n)$.
Thus when the continuous marginals $\mu$ and $\nu$ are properly approximated, so is the optimal transference plan between them.
Therefore, in the next sections, we will focus on developing algorithms to obtain $\mu^*_m, \nu^*_n$ that minimize
${W}_{k,\zeta}^k(\mu, \mu_m)$ and ${W}_{k,\zeta}^k(\nu, \nu_n)$.

%% file: 3_gradient.tex
\section{Gradient of the Objective Function}
\label{sec:gradient}

Let $\nu = \sum_{i=1}^{m}w_i\delta_{y_i}$ be the 
discrete probability measure in position of ``$\mu_m$''
in the last section.
The task now is to minimize ${W}_{k,\zeta}^k(\mu, \nu)$
about the
target discrete probability measure $\nu$ on $X$,
where $\mu$ is a fixed continuous probability measure on $X$.
The entropy-approximated
Wasserstein $W_{k,\zeta}^k(\mu,\nu)$ is convex
on $w_i$'s, while its ``convexity'' on $y_i$'s are not defined for a general $X$.

In this section, we derive the gradient of ${W}_{k,\zeta}^k(\mu, \nu)$ about $\nu$
following the discrete discussions of \cite{wang2019belieftransport,luise2018differential}, for applying stochastic
gradient descent (SGD) on both the positions $\{y_i\}_{i=1}^m$ and their weights
$\{w_i\}_{i=1}^m$. The SGD on $X$ is either through 
an exponential map, or by treating $X$ as (part of) an Euclidean space.

Let $g(x,y):=d_X^k(x,y)$, and denote the joint distribution minimizing
$\widehat{W}_{k,\zeta}^k$ as $\pi$ with differential form at $(x,y_i)$
being $\dd\pi_i(x)$, which is used to
define $W_{k,\zeta}^k$ in Section~\ref{sec:optimal_discretizations}. 


By introducing Lagrange multipliers $\alpha\in L^\infty(X), \beta\in\mathbb{R}^m$,
$\widehat{W}_{k,\zeta}^k(\mu,\nu) =
\max_{\alpha,\beta}\mathcal{L}(\mu, \nu; \alpha, \beta)$
where $\mathcal{L}(\mu, \nu; \alpha, \beta) =
\int_X\alpha(x)\dd\mu(x)+\sum_{i=1}^{n}\beta
w_i-\zeta\int_X\sum_{i=1}^{n}w_iE_i(x)\dd\mu(x)$ with $E_i(x)=e^{(\alpha(x)+\beta_ig(x,y_i))/\zeta}$ (See \cite{aude2016stochastic} or Supplementary). Let $\alpha^\ast,\beta^\ast$ be
the argmax, we
have $$W_{k,\zeta}^k(\mu.\nu)=\int_X\sum_{i=1}^{n}g(x,y_i)E_i^\ast(x)w_i\dd\mu(x)$$
with $E_i^\ast(x)=e^{(\alpha^\ast(x)+\beta^\ast_i-g(x,y_i))/\zeta}$. Since
$\alpha'(x):=\alpha(x)+t$ and $\beta_i':=\beta_i-t$ produces the same $E_i(x)$
for any $t\in\mathbb{R}$, the representative with $\beta_n=0$ that is equivalent to
$\beta$ (as well as $\beta^\ast$) is denoted by $\overline{\beta}$
(similarly $\overline{\beta}^\ast$) below, in order to get uniqueness, making the differentiation possible.

From a direct differentiation of $W_{k,\zeta}^k$, we have
\begin{align}\small
  \label{eq:gradient_W_w} &\!\!\!\!\dfrac{\partial W_{k,\zeta}^k}{\partial w_i}
=\int_Xg(x,y_i)E^\ast_i(x)\dd\mu(x)+\nonumber \\
&\!\!\!\dfrac{1}{\zeta}\int_X\sum_{j=1}^{n}g(x,y_j) \left(\dfrac{\partial
\alpha^\ast\!(x)}{\partial w_i} + \dfrac{\partial \beta^\ast_j}{\partial w_i}\right) w_j
E^\ast_j(x)\dd\mu(x).
\end{align}
\begin{align}\small
  \label{eq:gradient_W_y} &\!\!\!\!\nabla_{y_i}\!W_{k,\zeta}^k
=\!\int_X\!\!\!\nabla_{y_i}g(x,y_i)\left(\!1-\dfrac{g(x,y_i)}{\zeta}\!\right)
\!E^\ast_i(x)w_i\dd\mu(x)+\nonumber \\
&\!\!\frac{1}{\zeta}\!\!\int_X\sum_{j=1}^{n}g(x,y_j)
\!\left(\nabla_{y_i}\!\alpha^\ast\!(x)\!+\! \nabla_{y_i}\!\beta^\ast_j\right) w_j
E^\ast_j(x)\dd\mu(x).
\end{align} With the transference plan $\dd\pi_i(x)=w_iE^\ast_i(x)\dd\mu(x)$ and the
derivatives of $\alpha^\ast$, $\beta^\ast$, $g(x,y_i)$ calculated, the gradient of
$W_{k,\zeta}^k$ can be assembled.

Assume that $g$ is Lipschitz and is differentiable almost everywhere (for
$k\ge1$, and $d_X$ Euclidean 
distance in $\mathbb{R}^d$, differentiability fails to hold only when $k=1$ and
$y_i=x$), and that $\nabla_y g(x,y)$ is calculated. The
derivatives of $\alpha^\ast$ and $\overline{\beta}^\ast$ can then be calculated thanks to
Implicit Function Theorem for Banach spaces (see \cite{accinelli2009generalization}).

The maximality of $\mathcal{L}$ at $\alpha^\ast$ and $\overline{\beta}^\ast$
induces that 
$\mathcal{N}:=\nabla_{\alpha,\overline{\beta}}\mathcal{L}|_{(\alpha^\ast, \bar{\beta}^\ast)}=0\in(L^\infty(X)\otimes\mathbb{R}^{m-1})^\vee$,
the Fr\'echet derivative vanishes. Differentiate (in the sense of Fr\'echet)
again on $(\alpha, \overline{\beta})$ and $y_i, w_i$, respectively, we get
\begin{equation}
  \label{eq:hessian}
  \nabla_{(\alpha,\overline{\beta})}\mathcal{N} = -\dfrac{1}{\zeta}\left[
  \begin{array}{cc}
    \dd\mu(x)\delta(x,x') & \dd\pi_j(x')\\
    \dd\pi_i(x)           & w_i\delta_{ij}\\
  \end{array}
\right]
\end{equation}
as a bilinear functional on $L^\infty(X)\times\mathbb{R}^{m-1}$ 
(note that in Eq.~\eqref{eq:hessian}, the index $i$ of
$\dd\pi_i$ cannot be $m$).
The bilinear functional $\nabla_{(\alpha,\overline{\beta})}\mathcal{N}$ is
invertible, and denote its inverse by $\mathbf{M}$ as a bilinear form on
$\left(L^\infty(X)\otimes\mathbb{R}^{m-1}\right)^\vee$.

The last ingredient for Implicit Function Theorem is $\nabla_{w_i,
  y_i}\mathcal{N}$:
\begin{equation}\small
  \label{eq:derivative_alpha_w}
  \nabla_{w_i}\mathcal{N} =
  \left(-\dfrac{1}{w_i}\int_X(\ \cdot\ )\dd\pi_i(x), \vec{0}\right)
\end{equation}
\begin{align}\small
  \label{eq:derivative_alpha_y}
  \nabla_{y_i}\mathcal{N} = &
  \left(\dfrac{1}{\zeta}\!
  \int_X\!\!(\cdot)\nabla_{y_i}g(x,y_i)\dd\pi_i(x),\right. \nonumber\\
   &\left.\dfrac{\delta_{ij}}{\zeta}\!\int_X\!
   \nabla_{y_i}g(x,y_i)\dd\pi_i(x)\right).
\end{align}
Then $\nabla_{w_i,y_i}(\alpha^\ast,\overline{\beta}^\ast)=
\mathbf{M}(\nabla_{w_i, y_i}\mathcal{N})$.

Therefore, we have gradient $\nabla_{w_i,y_i}W_{k,\zeta}^k$
calculated.

Moreover, we can differentiate Eq.~\eqref{eq:gradient_W_w}, \eqref{eq:gradient_W_y},
\eqref{eq:hessian}, \eqref{eq:derivative_alpha_w}, \eqref{eq:derivative_alpha_y} to get a Hessian matrix of 
$W_{k,\zeta}^k$ on $w_i$'s and $y_i$'s provided a
better differentiability of $g(x,y)$ (which may enable Newton's method, or a mixture of Newton's method and minibatch SGD to accelerate the convergence). 
More details about claims, calculations and proofs are provided in the Supplementary Text.

%% file: 4_algorithm.tex
\section{The Discretization Algorithm}
\label{sec:algorithm_design}


Here we provide a description of an algorithm for Efficient Discretizations of Optimal Transport (EDOT)
from a distribution $\mu$ to $\mu_m$ with integer $m$, a given cardinality of support.
In general, $\mu$ need not be explicitly accessible, and even if it is accessible, exactly computing the transference plan
is not feasible. 
Therefore, in this construction,
we assume $\mu$ is given in terms of a \textbf{random sampler}, and apply
minibatch stochastic gradient descent (SGD) through a
set of samples independently drawn from $\mu$ of size $N$ on each step to
approximate $\mu$. 


Recall that to calculate the gradient
$\nabla_{\mu_m}{W}_{k,\zeta}^k(\mu, \mu_m)=
\left(\nabla_{x_i}{W}_{k,\zeta}^k(\mu,\mu_m),
  \nabla_{w_i}{W}_{k,\zeta}^k(\mu,\mu_m)\right)_{i=1}^{m}$, 
we need \textbf{1).}~$\pi_{X,\zeta}^{\phantom{|}}$, the EOT transference plan between $\mu$ and
$\mu_m$, \textbf{2).}~the cost $g = d_{\!X}^{\phantom{|}k}$ on $X$ and, \textbf{3).}~its gradient on the second variable
$\nabla_{x'}d_{\!X}^{\phantom{|}k}(x,x')$.
From $N$ samples $\{y_i\}_{i=1}^N$, we can construct
$\mu^{\phantom{|}}_N=\frac{1}{N}\sum_{i=1}^{N}\delta_{y_i}$ and calculate the
gradients with $\mu$ replaced by $\mu^{\phantom{|}}_N$ as an estimation, 
whose effectiveness
(convergence as $N\rightarrow\infty$) is 
induced by \cite{aude2016stochastic}.
%
The algorithm is stated in Algorithm \ref{alg:simple_EDOT} (EDOT).

In simulations, we choose $k=2$ to reduce complexity in calculating
the distance function (square of Euclidean distance is quadratic)
and its derivatives. The regularizer $\zeta$ should be small enough
to reduce the error of the EOT estimation $W_{k,\zeta}$ of the
Wasserstein distance $W_k$. However, setting $\zeta$ too small
will make $e^{-g(x,y)/\zeta}$ or its byproduct in Sinkhorn algorithm indistinguishable from
$0$ in \textit{double} format. Our strategy is setting $\zeta=0.01$
for $X$ of diameter~$1$ and scales propotional
with $\mathrm{diam}(X)^k$ (see next section).

\begin{algorithm}[H]
  \caption{Simple EDOT using minibatch SGD}
  \label{alg:simple_EDOT}
  \begin{algorithmic}\small
    \STATE {\bfseries input:} $\mu$, $k$, $m$, $\zeta$, $N$ batch size,  $\epsilon$ for stopping criterion,
    $\alpha=0.2$ for momentum, $\beta=0.2$ for learning rate.
    \STATE {\bfseries output:} $x_i\in X$, $w_i>0$ such that
    $\mu_m=\sum_{i=1}^{m}w_i\delta_{x_i}$.
    \STATE {\bfseries initialize:} randomly choose
    $\sum_{i=1}^{m}w_i^{(0)}\delta^{\phantom{|}}_{x_i^{(0)}}$; set $t=0$;\\
    \rule{1cm}{0pt}set learning rate $\eta_t=0.5(1+\beta t)^{-1/2}$ (for $t>0$).
    \REPEAT
    \STATE Set $t\leftarrow t+1$;
    \STATE Sample $N$ points $\{y_j\}_{j=1}^N\subseteq X$ from $\mu$ independently;
    \STATE Construct $\mu_N^{(t)}=\frac{1}{N}\sum_{j=1}^{N}\delta_{y_j}$;
    \STATE Calculate gradients
    $\nabla_{\mathbf{x}}\widehat{W}_{k,\zeta}^k(\mu_N^{(t)},\mu_m^{(t)})$ and \\
    \rule{2.3cm}{0pt}$\nabla_{\mathbf{w}}\widehat{W}_{k,\zeta}^k(\mu_N^{(t)},\mu_m^{(t)})$;
    \STATE Update $D\mathbf{x}_t\leftarrow\alpha
    D\mathbf{x}_{t-1}+\nabla_{\mathbf{x}}\widehat{W}_{k,\zeta}^k(\mu_N^{(t)},\mu_m^{(t)})$,\\
    \rule{1cm}{0pt}$D\mathbf{w}_t\leftarrow\alpha
    D\mathbf{w}_{t-1}+\nabla_{\mathbf{w}}\widehat{W}_{k,\zeta}^k(\mu_N^{(t)},\mu_m^{(t)})$; 
    \STATE Update
    $\mathbf{x}^{(t)}\leftarrow\mathbf{x}^{(t-1)}-\eta_tD\mathbf{x}_t$,\\
    \rule{1cm}{0pt}$\mathbf{w}^{(t)}\leftarrow\mathbf{w}^{(t-1)}-\eta_tD\mathbf{w}_t$;

    \UNTIL
    {$|\nabla_{\mathbf{x}}\widehat{W}_{k,\zeta}^k|+|\nabla_{\mathbf{w}}\widehat{W}_{k,\zeta}^k|<\epsilon$;}
    \STATE Set output $x_i\leftarrow x_i^{(t)}$, $w_i\leftarrow w_i^{(t)}$.
  \end{algorithmic}
\end{algorithm}

\subsection{Examples of Discretization}
\label{subsec:example_discretization}
We demonstrate our algorithm on:
\textbf{1).} $\mu$ is the uniform distribution on $X=[0,1]$.
\textbf{2).} $\mu$ is the mixture of two truncated normal distributions on
$X=[0,1]$, the PDF is $f(x)=0.3\phi(x;0.2, 0.1)+0.7\phi(x; 0.7, 0.2)$, where
$\phi(x;\xi, \sigma)$ is the density of the truncated normal distribution
on $[0,1]$ with expectation $\xi$ and standard deviation $\sigma$.
\textbf{3).} $\mu$ is the mixture of two truncated normal distributions on
$X=[0,1]^2$, the two distributions are: 
$\phi(x;0.2,0.1)\phi(y;0.3,0.2)$ of weight $0.3$ and
$\phi(x;0.7,0.2)\phi(y;0.6,0.15)$ of weight $0.7$.

\begin{figure}[ht]
  \centering
  \includegraphics[scale=0.35]{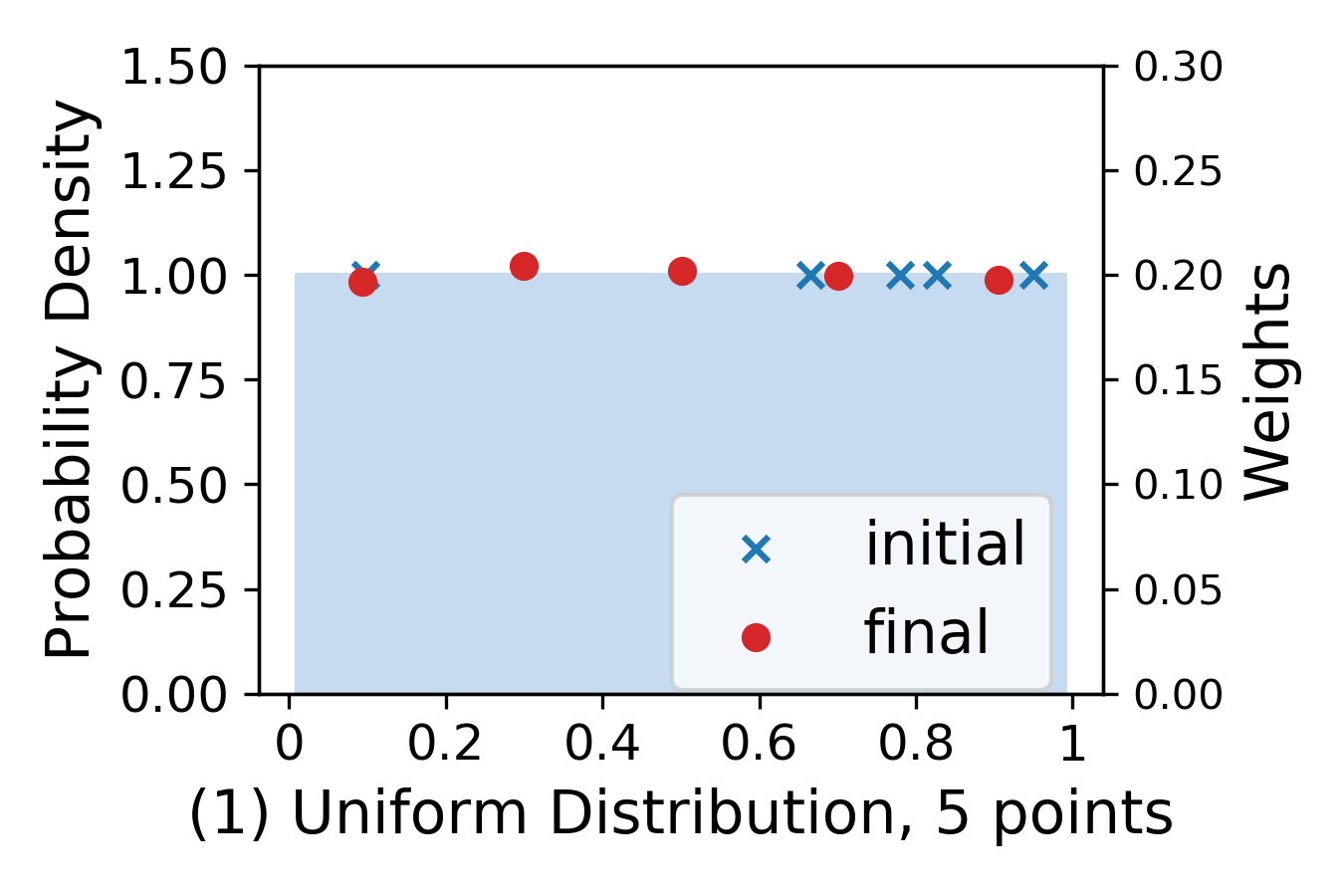}
  \includegraphics[scale=0.35]{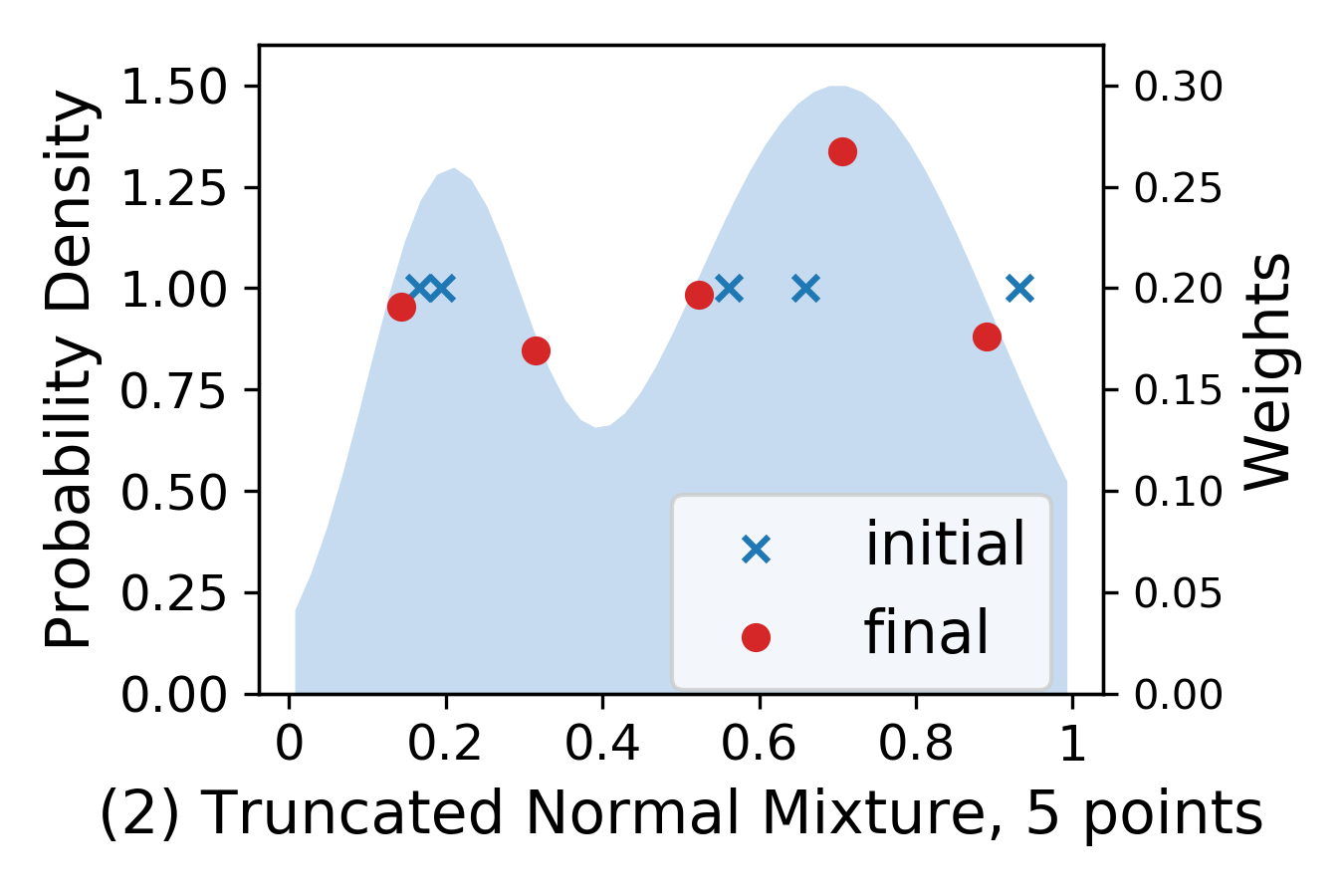}
  \includegraphics[scale=0.35]{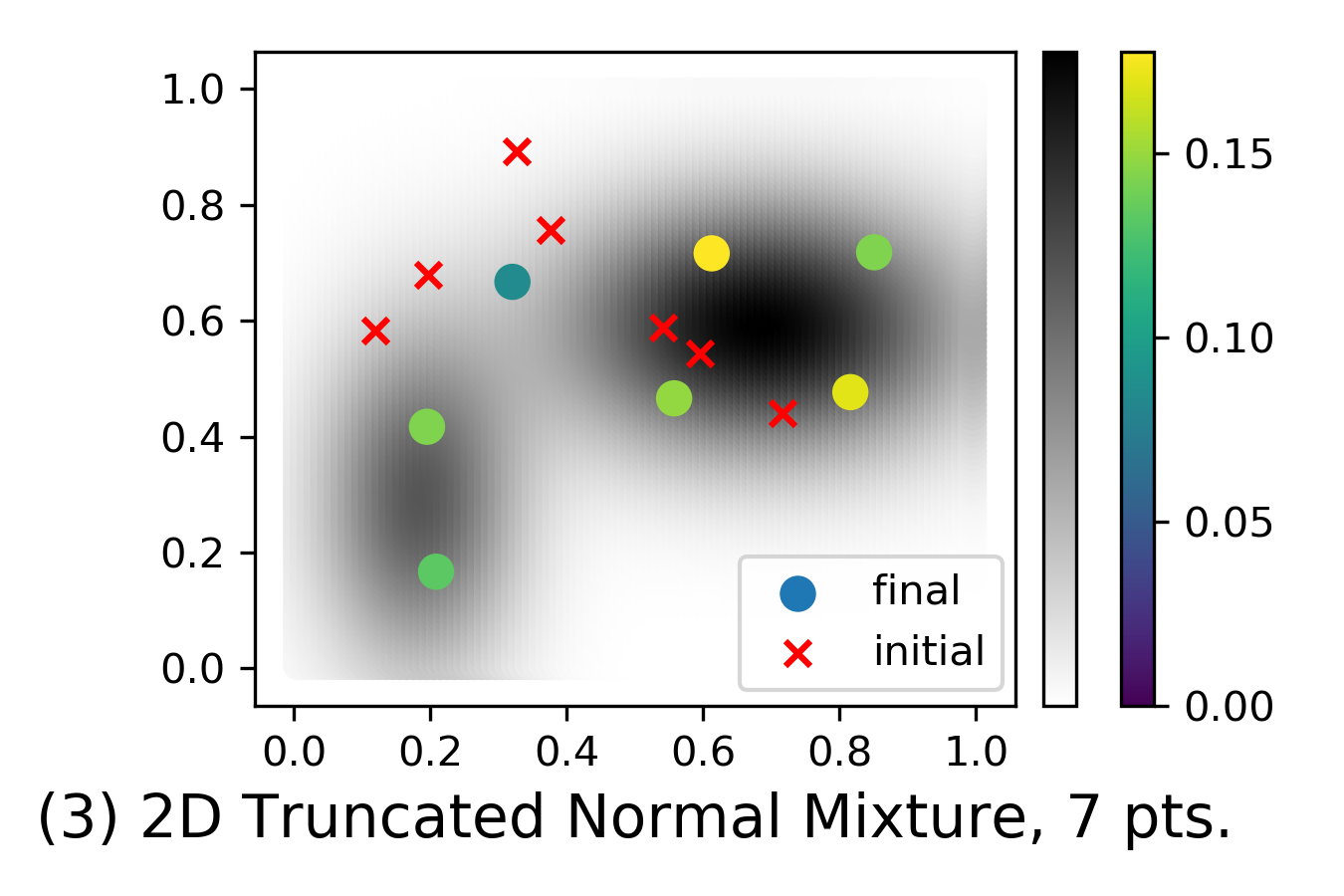}
  \includegraphics[scale=0.35]{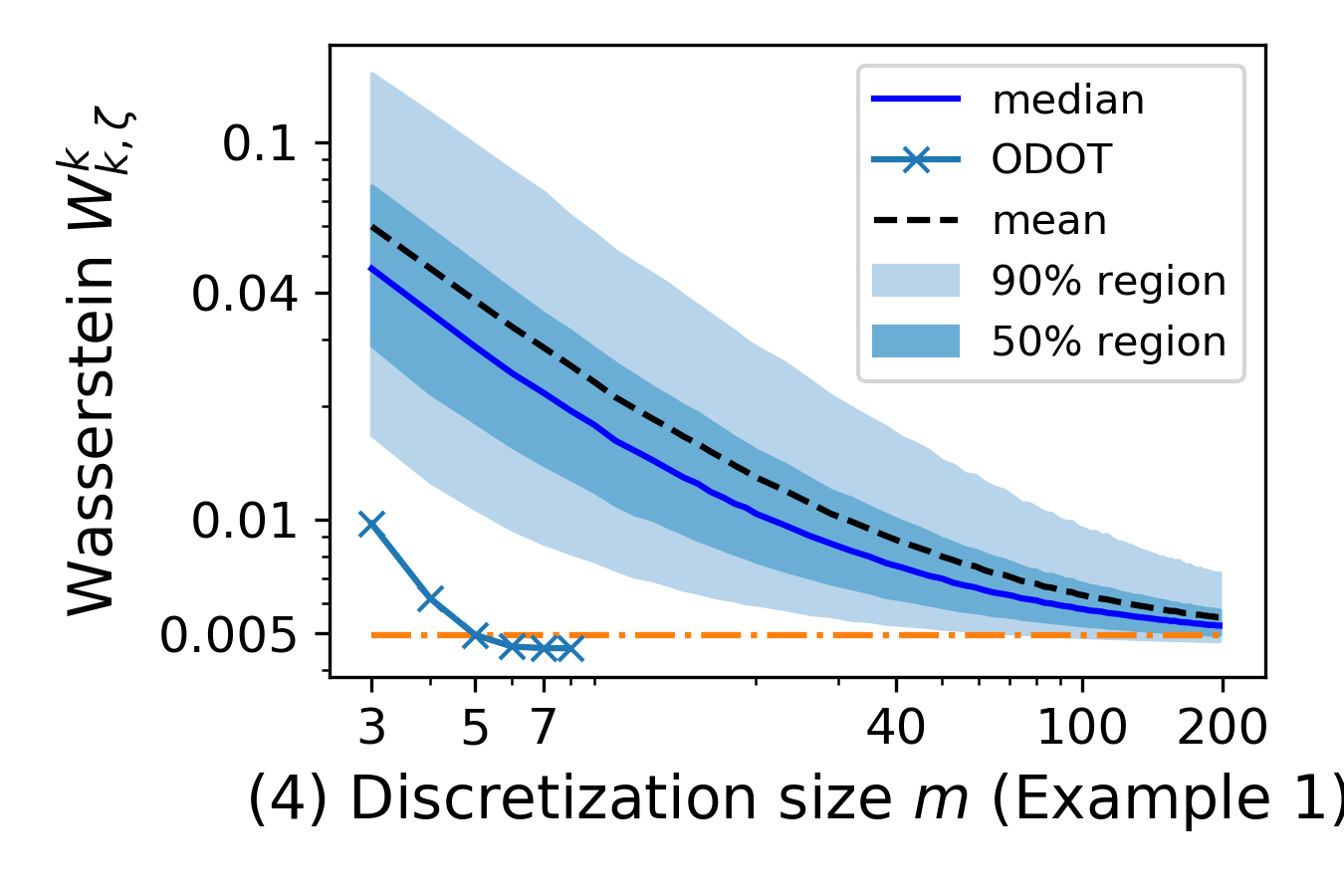}
  \includegraphics[scale=0.35]{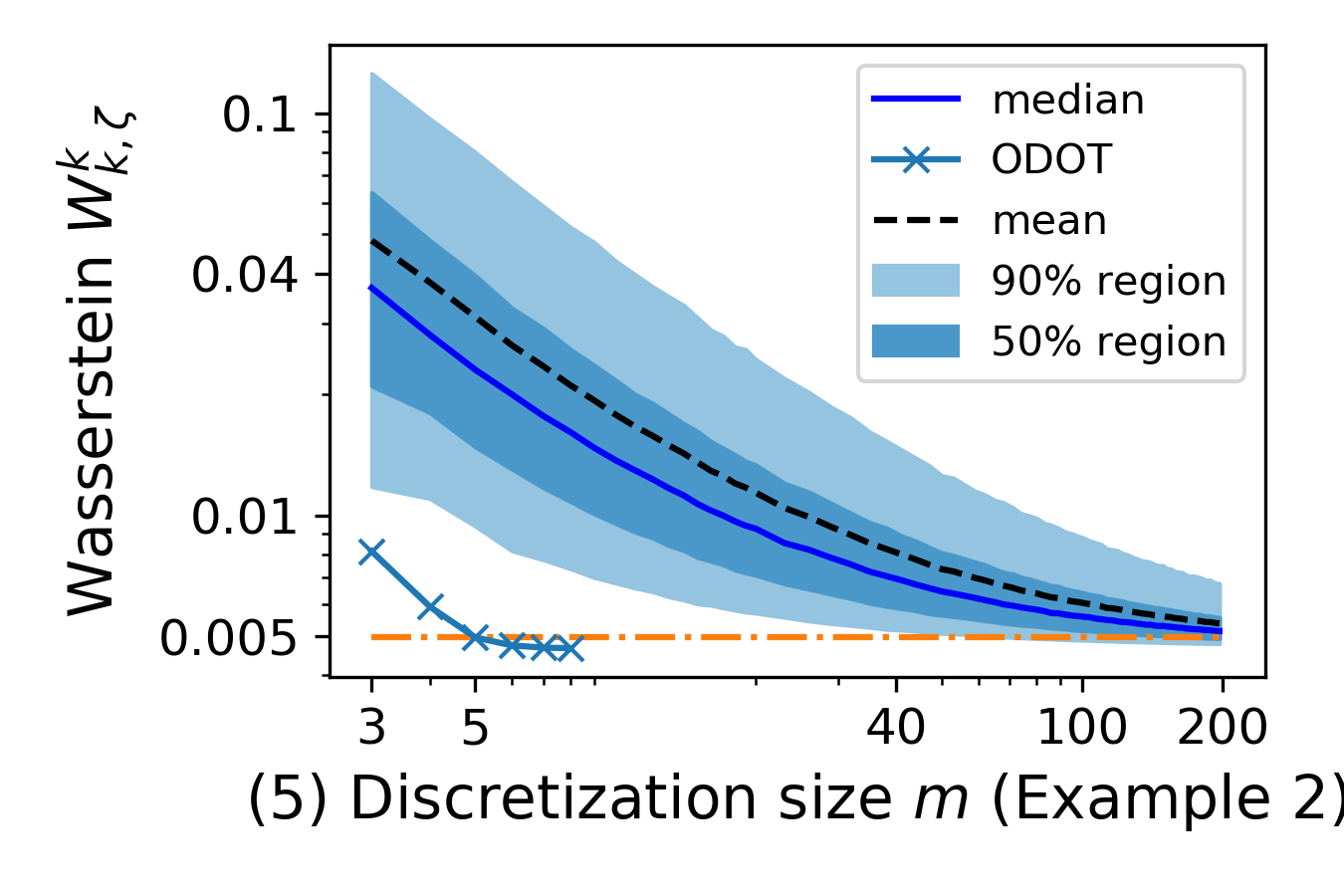}
  \includegraphics[scale=0.35]{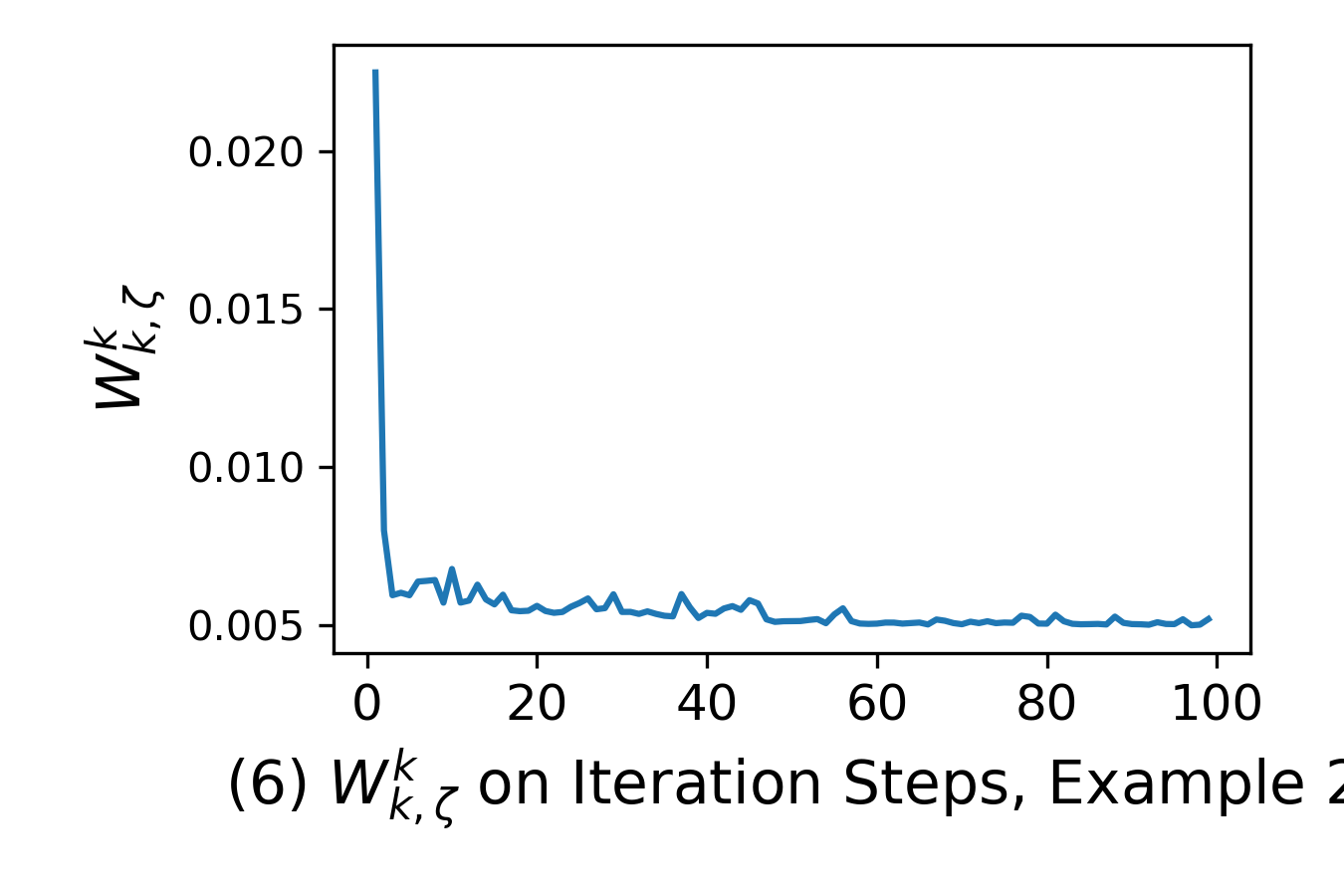}
  \caption{(1-3) plot sample discretizations of the examples (1-3). 
  In (3), $x$-axis and $y$-axis are the 2D coordinates, the probability density of $\mu$ and weights of $\mu_m$ are encoded by color. 
    (4,5) Comparisons with i.i.d. sampling on Wasserstein distance to the continuous distribution. EDOT are calculated with $m=3$ to $7$ ($3$ to $8$). The 4 boundary curves of the shaded region are $5\%$-, $25\%$-, $75\%$- and $95\%$-percentile curves, orange line represents the level of $m=5$; 
    (6) plots the entropy regularized Wasserstein distance $W_{k,\zeta}^k (\mu, \mu_m)$ versus SGD steps 
on example 2 (1D truncated normal mixture) with $\mu_m$ optimized by 5-point EDOT.\label{fig:example_disc}}
\end{figure}

Let $k=2$, $\zeta=0.01$, $N=100$ for all plots in this section. 
Fig.~\ref{fig:example_disc}~(1-3) plot 
the discretizations ($\mu_m$) for example~(1-3) with $m = 5, 5, 7$ respectively.

\textbf{Convergence.}
There is no guarantee that Algorithm~\ref{alg:simple_EDOT}
converges to the global optimal positions (though if stabilized at some 
positions, the weights are optimized).
On these examples we usually iterate 5000 to 6000 steps to 
get a result (gradients will be small, but are still affected by the fluctuation caused by the sample size $N$). 
We estimate entropy-regularized Wasserstein $W_{k,\zeta}^k$ between $\mu$ and $\mu_m$ by Richardson extrapolation (see Supplementary for details)
and use it to score the convergence rate.

Fig.~\ref{fig:example_disc} (6) illustrates the convergence rate of $W_{k,\zeta}^k (\mu, \mu_m)$ versus SGD steps 
on example 2 (1D truncated normal mixture) with $\mu_m$ obtained by 5-point EDOT. 

Generally, Wasserstein distance decays and gets close to 
final state after tens of steps, while
slight changes 
can still be
observed (in 1D examples)
even after 3000 steps. Sometimes an ill 
positioned initial distribution can make it slower
to approach a stationary state.

\textbf{Comparison with Naive Sampling.}
Fig.~\ref{fig:example_disc} (4-5) plot the entropy-regularized Wasserstein ${W}_{k,\zeta}^k(\mu,\mu_m)$ with different choices of $m$ for Example 1) and 2). 
Here $\mu_m$'s are generated by two methods: 
a). Algorithm ~\ref{alg:simple_EDOT} with $m$ chosen from 3 to 7 in Example 1 and to 8 in Example 2,
shown by $\times$ 's.
b). naive sampling (i.i.d. with equal weights) simulated using Monte Carlo of volume 20000 on each size from 3 to 200. 
Fig.~\ref{fig:example_disc} (4-5) demonstrate the effectiveness of EDOT: 
as indicated by the orange horizontal dashed line, even 5 points in these 2 examples 
excels 95\% of naive samplings of size 40 and 75\% of naive samplings of size over 100. 




\subsection{An Example on Transference Plan}
\label{subsec:example_transference}
\begin{figure}[ht]
    \centering
    \includegraphics[scale=0.39, trim={20 10 10 10}, clip]{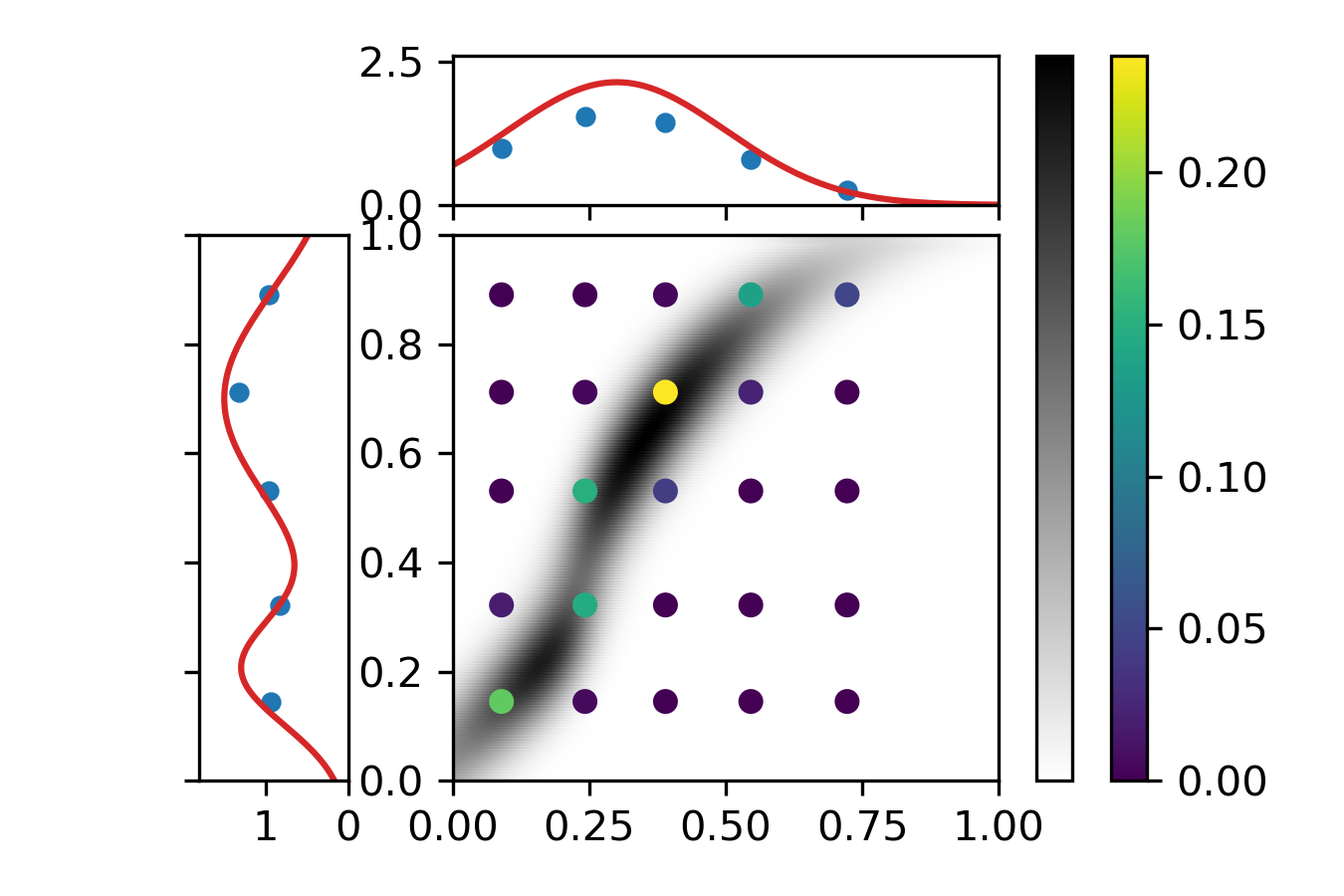}
    \includegraphics[scale=0.39, trim={20 10 10 10}, clip]{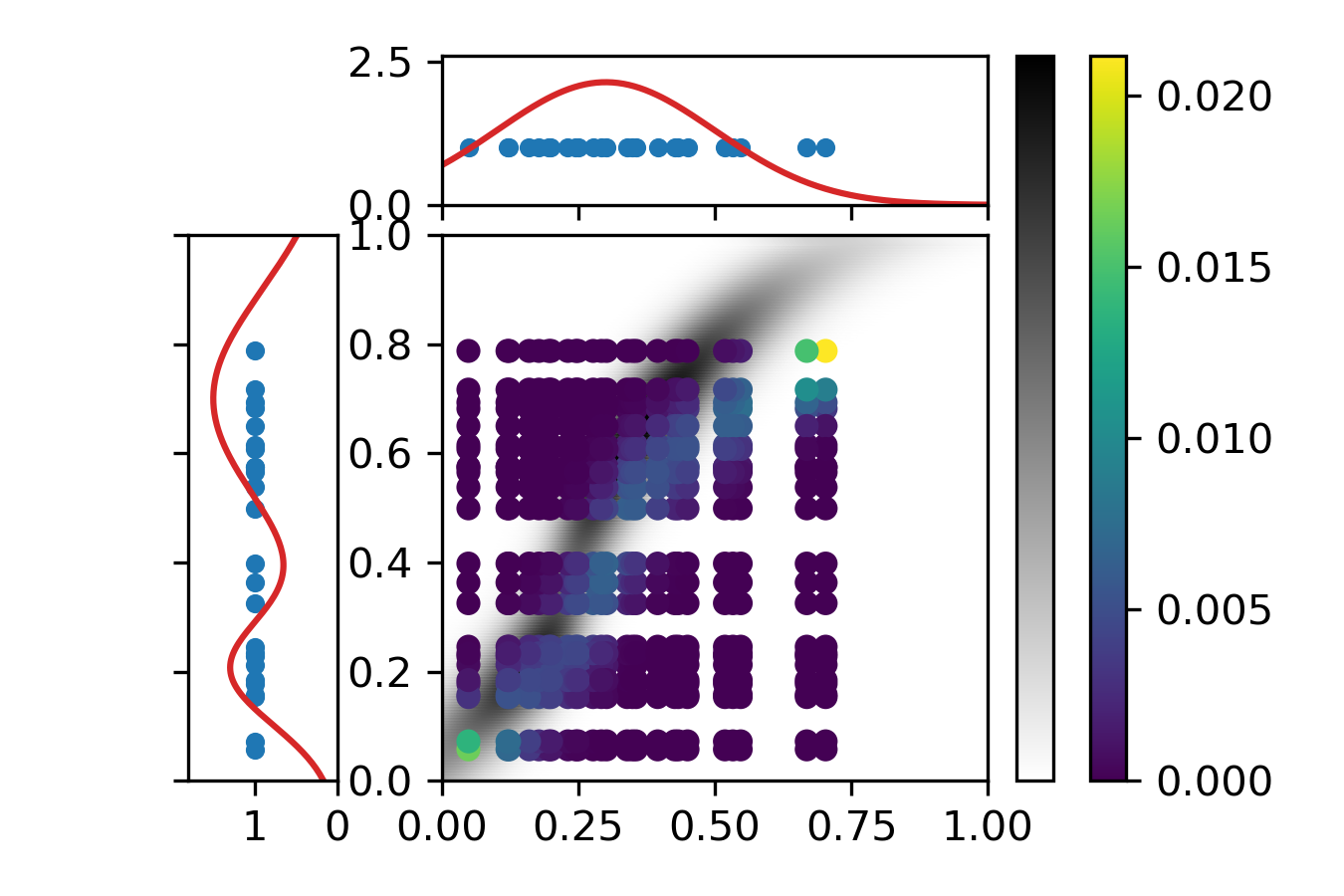}
    \caption{Approximation of the transference plan. 
    Left: EDOT on $5\times5$ grid with
    $W_{k,\zeta}^k(\mu,\mu_5)=4.792\times10^{-3}$,
    $W_{k,\zeta}^k(\nu,\nu_5)=5.034\times10^{-3}$,
    $W_{k,\zeta}^k(\gamma,\gamma_{5,5})=8.446\times10^{-3}$;
    Right: Naive on $25\times25$ grid with 
    $W_{k,\zeta}^k(\mu,\mu_7)=5.089\times10^{-3}$,
    $W_{k,\zeta}^k(\nu,\nu_7)=2.222\times10^{-2}$,
    $W_{k,\zeta}^k(\gamma,\gamma_{7,7})=2.563\times10^{-2}$. In both figures, the background grayscale 
    density represents the true EOT density.
    }
    \label{fig:joint_distribution}
\end{figure}

In Fig.~\ref{fig:joint_distribution}, we illustrate an application of EDOT on 
an optimal transport task as proposed in Section~\ref{sec:optimal_discretizations}. 
The optimal transport problem has $X=Y=[0,1]$,
and the marginal distributions $\mu$, $\nu$ are truncated normal (mixtures), where
$\mu$ has 2 components (shown on the left) and $\nu$ has only 1 component (shown on the top).
The cost function is taken as the
square of Euclidean distance on an interval.
We can see that on the $5\times 5$ EDOT example, the high
density area of the transference plan are correctly 
covered by lattice points with high weights,
while in the naive sampling, even in grid size $25\times25$,
the points on the lattice with highest weights
missed the region where true transference plan is of the most density. Comparison on Wasserstein distance
also tells us EDOT works better with a smaller grid.

%% file: 5_parallel.tex
\section{Methods of Improvement}
\label{sec:practical_methods}



\subsection{Adaptive Cell Refinement}
\label{subsec:divide_and_conquer}
The computational cost of simple EDOT increases with the dimensionality and diameter of the underlying space.
Larger $m$ discretization $\mu_m=\sum_{i=1}^mw_i\delta_{y_i}$ is needed to capture the higher dimensional distributions.
This will result an increase in parameters in SGD for calculating the gradient of ${W}_{k,\zeta}^k$:
$md$ for positions $y_i$'s  and $m-1$ for weights $w_i$'s.
Such an increment will both increase complexity in each step, and also 
require more steps for SGD to converge.
Furthermore, the calculation will have a higher complexity ($\mathcal{O}(mN)$ for each
normalization in Sinkhorn).

We propose to reduce computational complexity by ``divide and conquer''. 
The Wasserstein distance takes $k$-th power of the distance function
$d_X^{\phantom{|}k}$ as cost function. The locality of distance
$d_X^{\phantom{|}}$ makes the solution to the OT / EOT problem local,
meaning the probability mass is more likely to be transported to a close
destination than to a remote one. Thus we can ``divide and conquer'' ---
cut the space $X$ into small cells and solve the discretization problem separately.
We require an adaptive dividing procedure to balance the
accuracy and computational intensity among the cells.
Therefore, determining size of discretization
and choosing a proper regularizer for each cell are questions to answer,
after having a partition $X=X_1\sqcup\dots\sqcup X_{\mathcal{I}}$.
We first sample a large sample set $S \sim \mu$ (with $|S|\gg N$) 
and partition $S$ into $S_1\sqcup \dots\sqcup S_{\mathcal{I}}$ according to $X$.
In other words, to construct an EDOT problem for each $X_i$ and sample set
$S_i$ located in $X_i$, we must
figure out the parameters $m_i$ and $\zeta_i$.

\textbf{Choosing Size} $m$.
An appropriate choice of $m_i$ will balance contributions to the Wasserstein
among the subproblems as follows: 
Let $X_i$ be a manifold of dimension $d$, let $\mathrm{diam}(X_i)$ be its diameter,
and $p_i={\mu}(X_i)$ be the probability of $X_i$.
The entropy regularized wasserstein distance can be estimated as
$W_{k,\zeta}^k=\mathcal{O}(p_im_i^{-k/d}\mathrm{diam}(X_i)^k)$ \cite{weed2019sharp, dudley1969speed}.
The contribution to $W_{k,\zeta}^k(\mu,\mu_m)$ per point in support of $\mu_m$ is
$\mathcal{O}(p_im_i^{-(k+d)/d}\mathrm{diam}(X_i)^k)$. Therefore, to
balance each point's contribution to the Wasserstein among the divided subproblems,
we set $m_i\approx\frac{(p_i\mathrm{diam}(X_i)^k)^{d/(k+d)}}{\sum_{j=1}^{\mathcal{I}}(p_j\mathrm{diam}(X_j)^k)^{d/(k+d)}}$.
Since $m_i$ must be an integer, it is rounded properly. 
If some cell has $m_i=0$, then the probability should be added to its adjacent cells.\footnote{In general, there are 
various ways of distributing the probability to neighbors (in Algorithm~\ref{alg:refine_dfs}, it is unique).}


\textbf{Adjusting Regularizer} $\zeta$.
In the calculation of $W_{k,\zeta}^k$,
the Sinkhorn iterations on $e^{-g(x,y)/\zeta}$
is calculated. Therefore, $\zeta$ should 
scale with $g(x,y)$ (i.e., $d_X^k$)
to make the transference plan not be
affected by scaling of $d_X$.
Precisely, we may choose $\zeta_i=\mathrm{diam}(X_i)^k\zeta_0$ for
some constant $\zeta_0$.

\textbf{The Construction.}
Theoretically, the idea of division can be applied to any refinement procedure
that can be applied iteratively and eventually makes the
diameter of the cells approach $0$. In our simulation, we use an
adaptive kd-tree style cell refinement in an Euclidean space~$\mathbb{R}^d$.

Let $X$ be embedded into $\mathbb{R}^d$ within an axis-aligned rectangular region.
We choose an axis $\mathbf{x}_l$ in $\mathbb{R}^d$ and evenly split the region
along a hyperplane orthogonal to $\mathbf{x}_l$ (e.g. cut square $[0,1]^2$ along
the line $x=0.5$), thus we construct $X_1$ and $X_2$.
With a sample set $S$ given, we split it into two sample sizes $S_1$ and $S_2$
according to which subregion each sample is located in. Then the corresponding
$m_i$ and $\zeta_i$ can be calculated as discussed above. Thus two cells and
corresponding subproblems are constructed. If some of the $m_i$ is still
too large, the cell is cut along another axis to construct two other cells.
The full list of cells and subproblems can be constructed recursively, see Algorithm~\ref{alg:refine_dfs}.

After having the set of subproblems, we may apply Algorithm~\ref{alg:simple_EDOT}
for solutions in each cell (samples used there are redrawn from the sample set
$S_i$ in each subproblem), then combine the solutions 
$\mu_{m_i}^{(i)}=\sum_{j=1}^{m_i}w_j^{(i)}\delta^{\phantom{|}}_{y_j^{(i)}}$
into the final result
$\mu_m:=\sum_{i=1}^{\mathcal{I}}\sum_{j=1}^{m_i}p_iw_j^{(i)}\delta^{\phantom{|}}_{y_j^{(i)}}$. 

\begin{algorithm}[ht]
  \caption{Adaptive Refinement via Depth First Search \label{alg:refine_dfs}}
  \begin{algorithmic}\small
    \STATE {\bfseries input:} $m$, $\mu$,
    $N_0$ sample size, $m_\ast$ max number of points in a cell,
    $\mathbf{a}=(a_0, a_1, \dots, a_{d-1})$ and
    $\mathbf{b}=(b_0, b_1,\dots,b_{d-1})$ as lower/upper bounds of the region;
    \STATE {\bfseries output:} $out$: stack of subproblems $(S_i, m_i, p_i)$ with $\sum p_i=1$, $\sum m_i=m$,
    $\sum |S_i|=N_0$;
    \STATE {\bfseries initialization:} 
    $p_0\leftarrow1$;
    $T$, $out$ be empty stacks;
    \STATE Sample $N_0$ points $S_0\sim\mu$;
    \STATE $T$.push(($S_0$, $m$, $p_0$, $\mathbf{a}$, $\mathbf{b}$));
    \WHILE{$T$ is not empty}
    \STATE ($S$, $m$, $p$, $\mathbf{a}$, $\mathbf{b}$) $\leftarrow$ $T$.pop();
    \STATE $l\leftarrow \argmax\{b_i-a_i\}$, 
    $mid\leftarrow (a_{l}+ b_{l})/2$;
    \STATE $\mathbf{a}^{(1)}\leftarrow \mathbf{a}$, 
    $\mathbf{a}^{(2)}\leftarrow(a_0,\dots, mid,\dots,a_{d-1})$;
    \STATE $\mathbf{b}^{(1)}\leftarrow(b_0,\dots, mid,\dots,b_{d-1})$, $\mathbf{b}^{(2)}\leftarrow \mathbf{b}$;
    \STATE $S_1\leftarrow\{\mathbf{x}\in S: x_l <= mid\}$,
    $S_2\leftarrow\{\mathbf{x}\in S: x_l > mid\}$, $N_1\leftarrow |S_1|$, $N_2\leftarrow |S_2|$;
    \STATE $m_1\leftarrow
    Round_{0.5\uparrow}\left(\frac{N_1^{d/(k+d)}}{N_1^{d/(k+d)}+N_2^{d/(k+d)}}\right)$;
    \STATE $m_2\leftarrow
    Round_{0.5\downarrow}\left(\frac{N_2^{d/(k+d)}}{N_1^{d/(k+d)}+N_2^{d/(k+d)}}\right)$;
    \IF{$m_1 = 0$}
    \STATE $p_1\leftarrow0$, $p_2\leftarrow(N_1+N_2)/N_0$;
    \ELSIF {$m_2=0$}
    \STATE $p_1\leftarrow(N_1+N_2)/N_0$, $p_2\leftarrow0$;
    \ELSE
    \STATE $p_1\leftarrow N_1/N_0$, $p_2\leftarrow N_2/N_0$;
    \ENDIF
    \FOR {$i\gets 1, 2$}
    \IF {$m_i>m_{\ast}$}
    \STATE $T$.push($(S_i, m_i, p_i, 
    \mathbf{a}^{(i)},\mathbf{b}^{(i)})$);
    \ELSIF {$m_i>0$}
    \STATE $out$.push($(S_i, m_i, p_i$);
    \ENDIF
    \ENDFOR
    \ENDWHILE
  \end{algorithmic}
\end{algorithm}

\begin{figure}[!ht]
  \centering
  \includegraphics[scale=0.4]{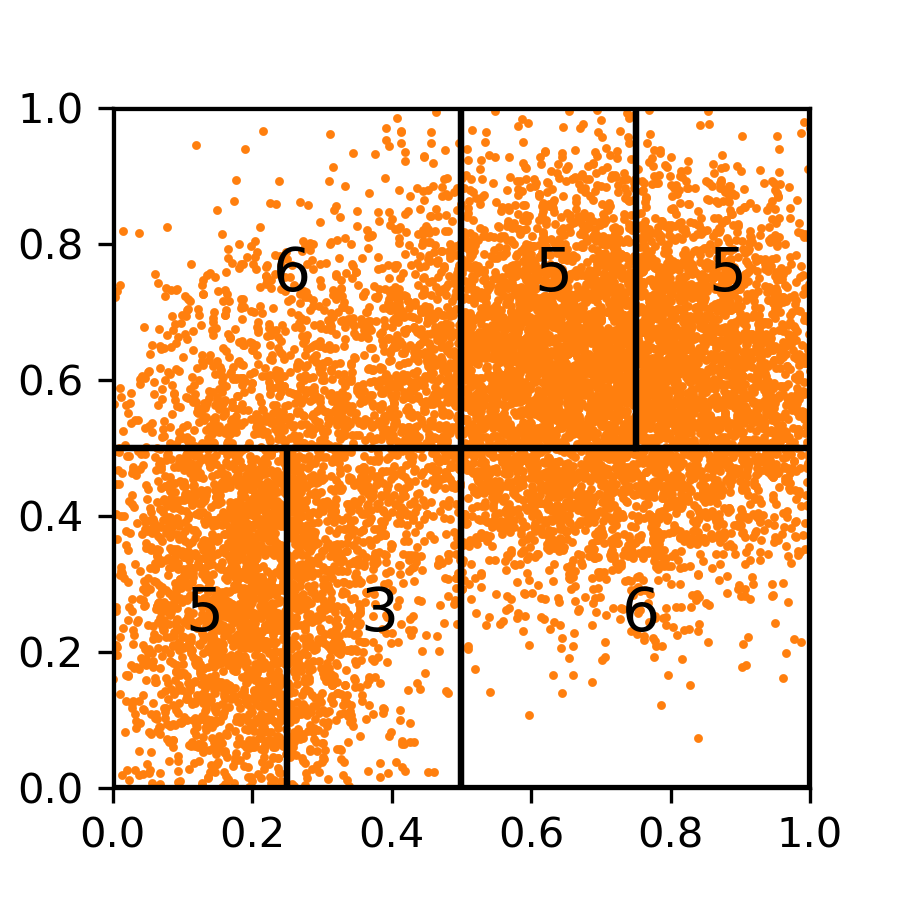}
  \includegraphics[scale=0.4]{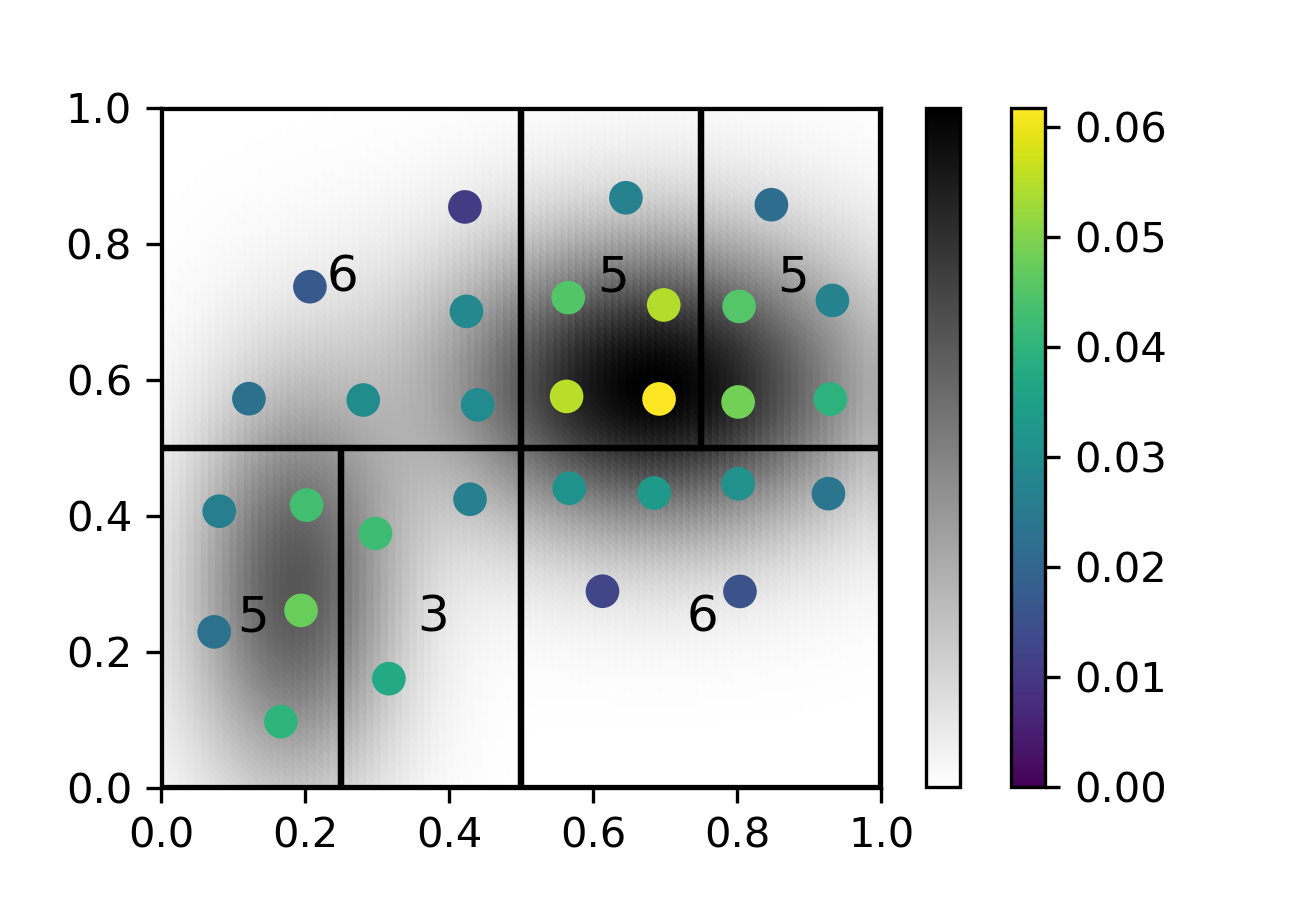}
  \caption{\label{fig:2d_refinement}An example of adaptive refinement on a
    unit square. Left: division of 10000 sample $S$ approximating a mixture of two truncated normal distributions, and the refinement for 30 discretization points.
  Right: the discretization optimized locally and combined as a probability
  measure, with $k=2$.}
\end{figure}

Fig.~\ref{fig:2d_refinement} shows the optimal discretization for the example
in Fig.~\ref{fig:example_disc}(3)) with $m=30$, obtained by applying the EDOT with adaptive cell refinement.

\subsection{On Embedded CW-Complexes}
\label{subsec:embedded}

Although the samples on space $X$ are usually represented as a vector in
$\mathbb{R}^d$, inducing an embedding $X\hookrightarrow\mathbb{R}^d$, the
space $X$ usually has its own structure as a CW-complex (or simply a manifold)
with a metric.

As the metric in the CW-complex usually bears a more intrinsic structure of $X$
than the one induced by the embedding, if the CW-complex structure and the
metric is known, even piecewise, we may apply
Algorithm~\ref{alg:simple_EDOT} or Algorithm~\ref{alg:refine_dfs} on each cell
or piece to get a precise discretization of $X$ regarding its own metric, whereas direct discretization 
as a subset in $\mathbb{R}^d$ may result in low expressing efficiency e.g. some points may
drift off from $X$.

We show two examples: truncated normal mixture distribution on a Swiss roll, and
a mixture normal distribution on a sphere mapped through stereographic
projection. 

\begin{figure}[ht]
  \centering
  \includegraphics[scale=0.3, trim={150 105 52 50}, clip]{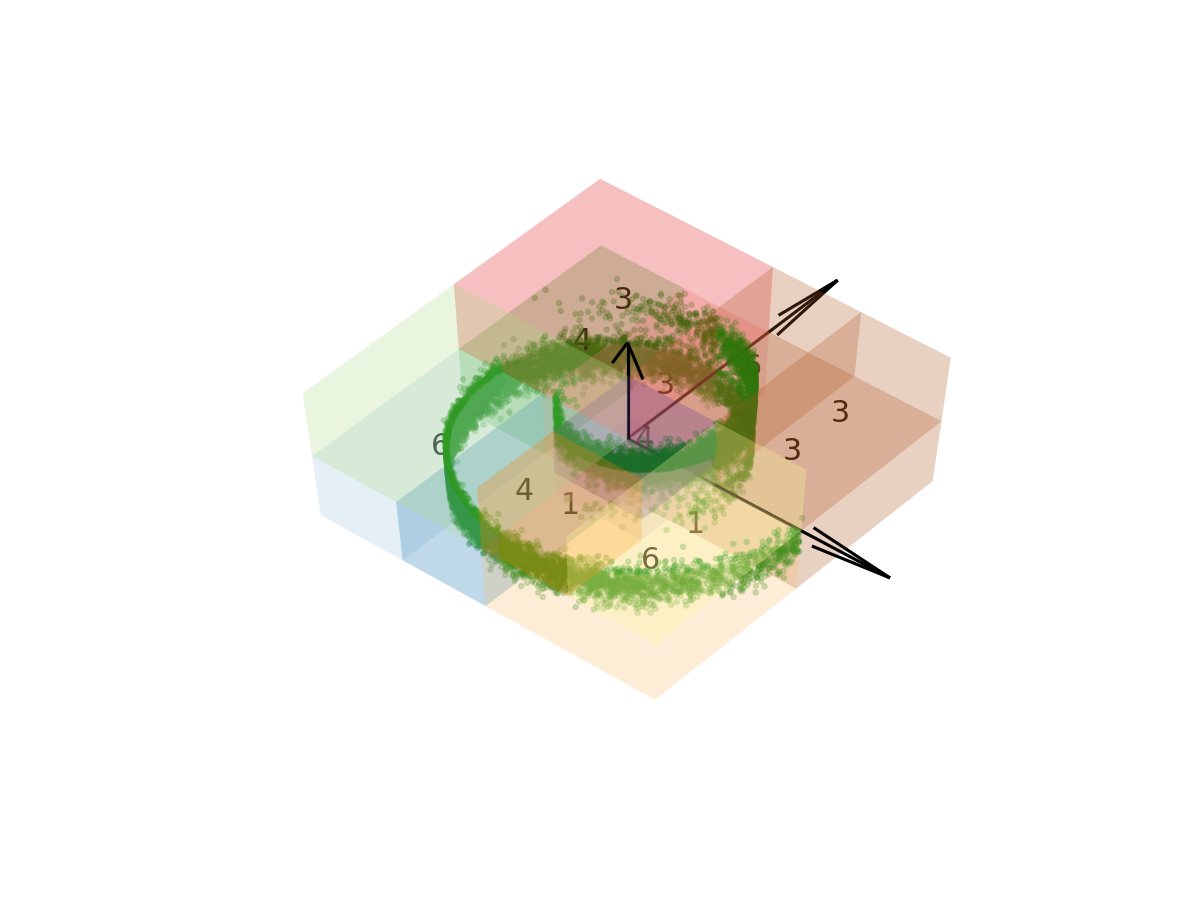} 
  \includegraphics[scale=0.24, trim={150 80 52 65}, clip]{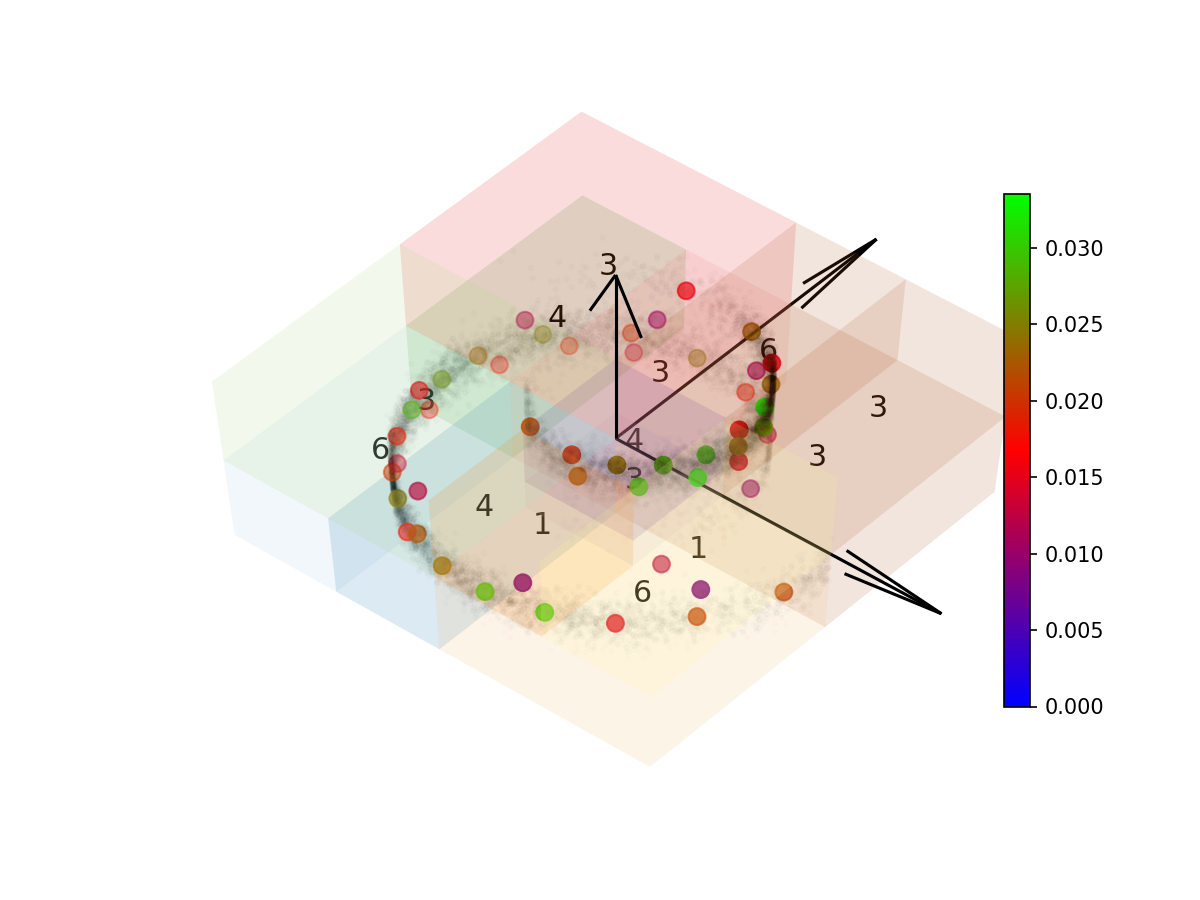}
  \caption{\label{fig:swiss_roll}On Swiss Roll. Left: 15000 samples from the truncated normal mixture distribution $\mu_{\text{swiss}}$ over $X_{\text{swiss}}$; 
  Right: 50-point 3D discretization using Algorithm~\ref{alg:refine_dfs}, the refinement cells are shown in colored boxes.}
\end{figure}

\begin{figure}[ht]
  \centering
  \includegraphics[scale=0.4]{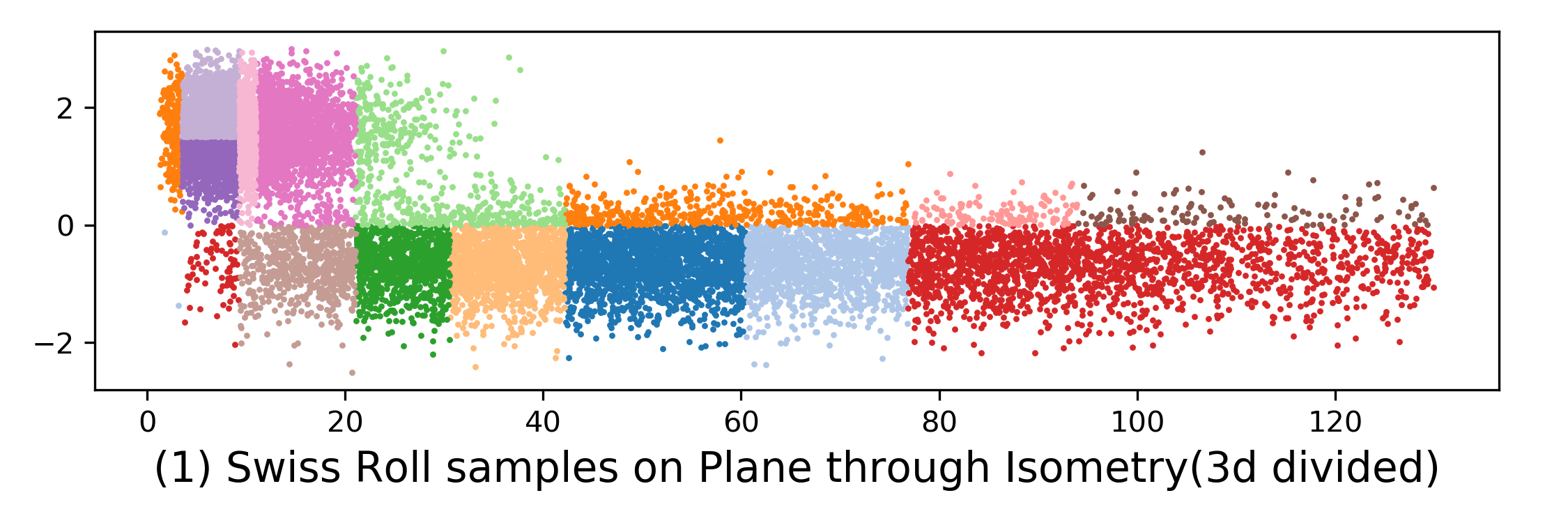}
  \includegraphics[scale=0.4]{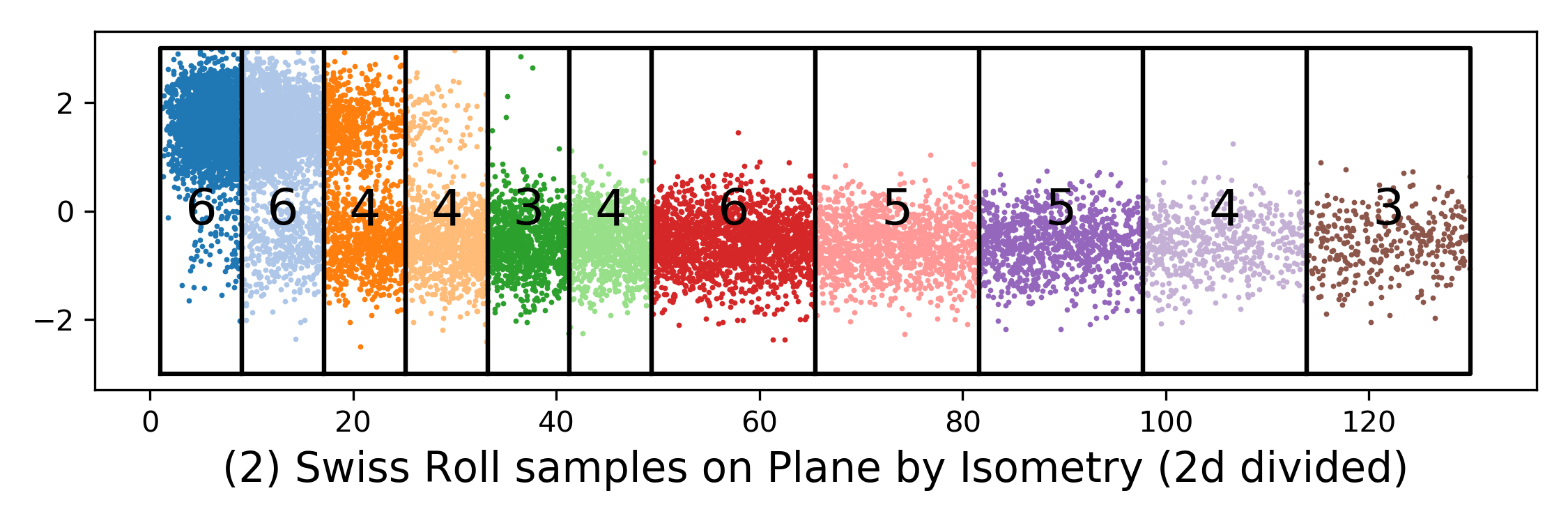}
  \includegraphics[scale=0.4]{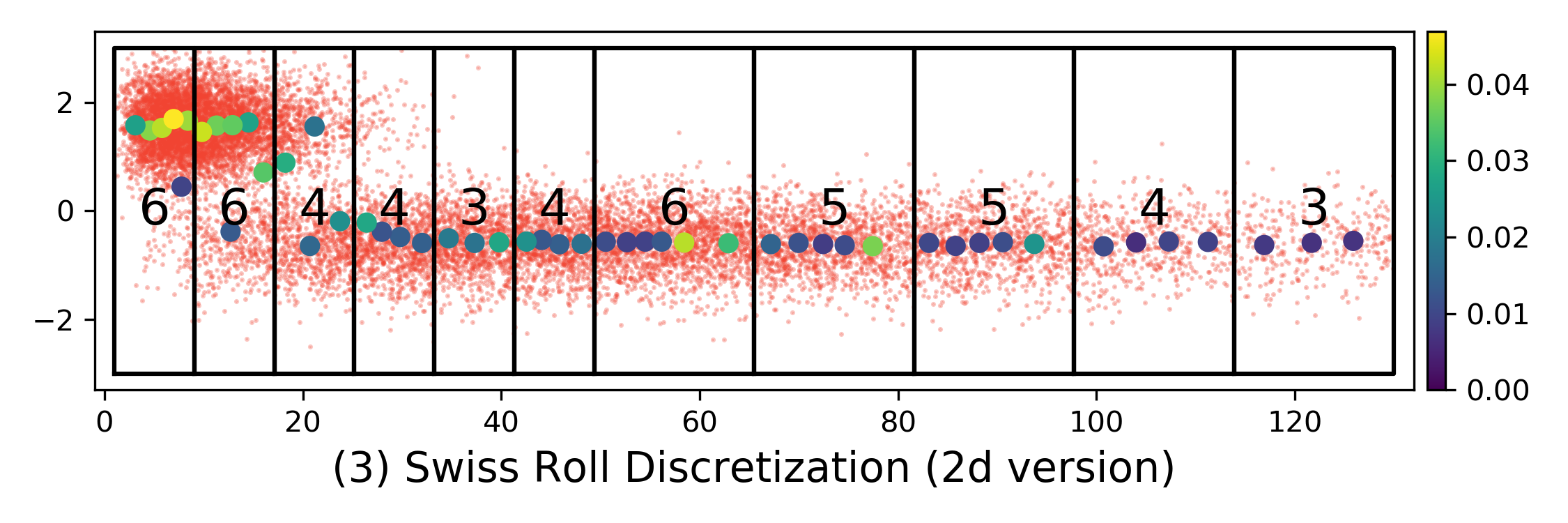}
  \caption{\label{fig:swiss_2d} Swiss Roll under Isometry. (1) Refinement Cells under 3D Euclidean metric (one color per samples from a refinement cell); 
  (2) Refinement Cells under 2D metric induced by the isometry; (3) EDOT of 50 points with respect to the 2D refinement.}
\end{figure}

\textbf{On Swiss Roll.}
%
In this case, the underlying space $X_{\text{swiss}}$ is a the Swiss Roll, a 2D rectangular strip embedded in $\mathbb{R}^3$: $(\theta,z)\mapsto(\theta,\theta,z)$ in cylindrical coordinates, $\mu_{\text{swiss}}$ is a truncated normal mixture on $(\theta,z)$-plane.
Samples $S_{\text{swiss}}\sim\mu_{\text{swiss}}$ over $X_{\text{swiss}}$ is shown on Fig.~\ref{fig:swiss_roll} (left) embedded in 3D 
and Fig.~\ref{fig:swiss_2d} (1) isometric into $\mathbb{R}^2$.

Following the Euclidean metric in $\mathbb{R}^3$,
Fig.~\ref{fig:swiss_roll} (right) plots the EDOT solution $\mu_{m}$ through adaptive cell refinement (Algorithm~\ref{alg:refine_dfs}) with $m=50$.
The resulting cell structure is shown as colored boxes. 
The corresponding partition of $S_{\text{swiss}}$ is shown on Fig.~\ref{fig:swiss_2d} (1), with samples contained in a cell marked by the same color.
According to Fig.~\ref{fig:swiss_roll} (right), points
in $\mu_m$ are mainly located on the strip with
only one point off in the most sparse cell (yellow cell located in the bottom in the figure).

On the other hand, consider the metric on $X_{\text{swiss}}$ induced by the isometry from the Swiss Roll as a manifold to 
a strip on $\mathbb{R}^2$. A more intrinsic discretization of $\mu_{\text{swiss}}$ can be obtained by applying EDOT through a refinement on the coordinate space, the (2D) strip.
The partition of $S_{\text{swiss}}$ is shown on Fig.~\ref{fig:swiss_2d} (2), and resulting discretization $\mu_{50}$ is shown in Fig.~\ref{fig:swiss_2d} (3).
Notice that all $50$ points are located on the (locally) high density region of the Swiss Roll.
We observe from Fig.~\ref{fig:swiss_2d} (1) and (2) that
the 3D partition may pull disconnected and intrinsically remote regions together
while the 2D partition maintains the intrinsic structure. 

\begin{figure}[ht]
    \centering
    \includegraphics[scale=0.4, trim={100 100 90 80}, clip]{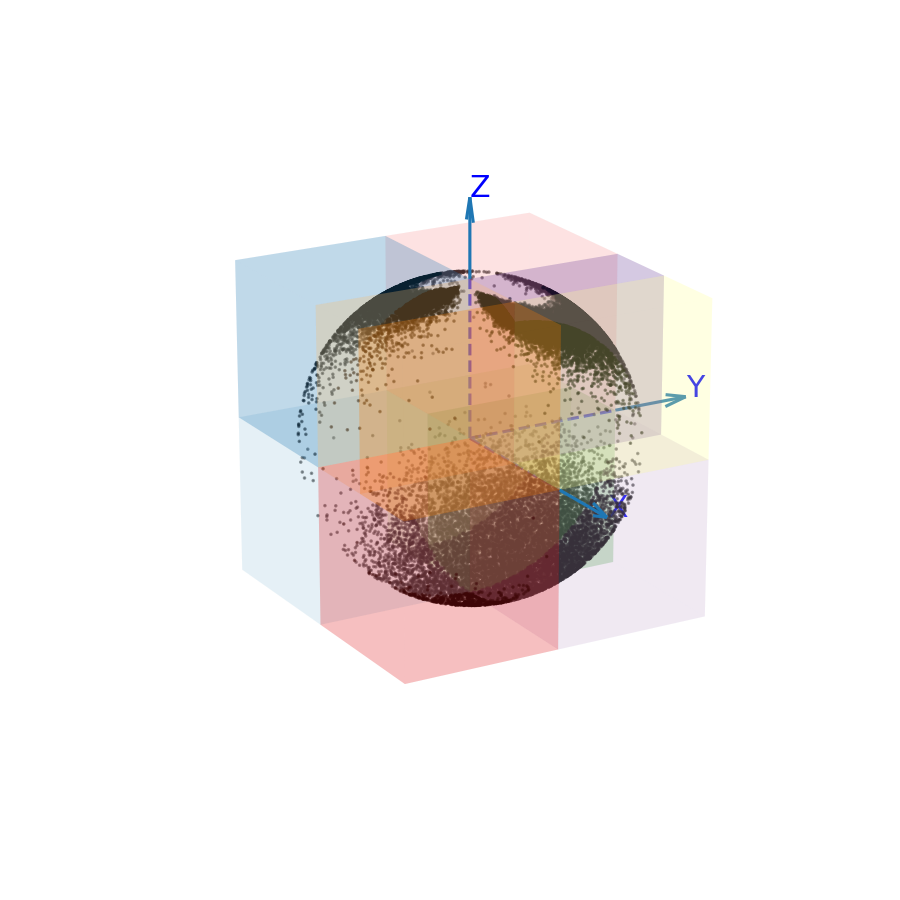}
    \includegraphics[scale=0.4, trim={100 110 30 70}, clip]{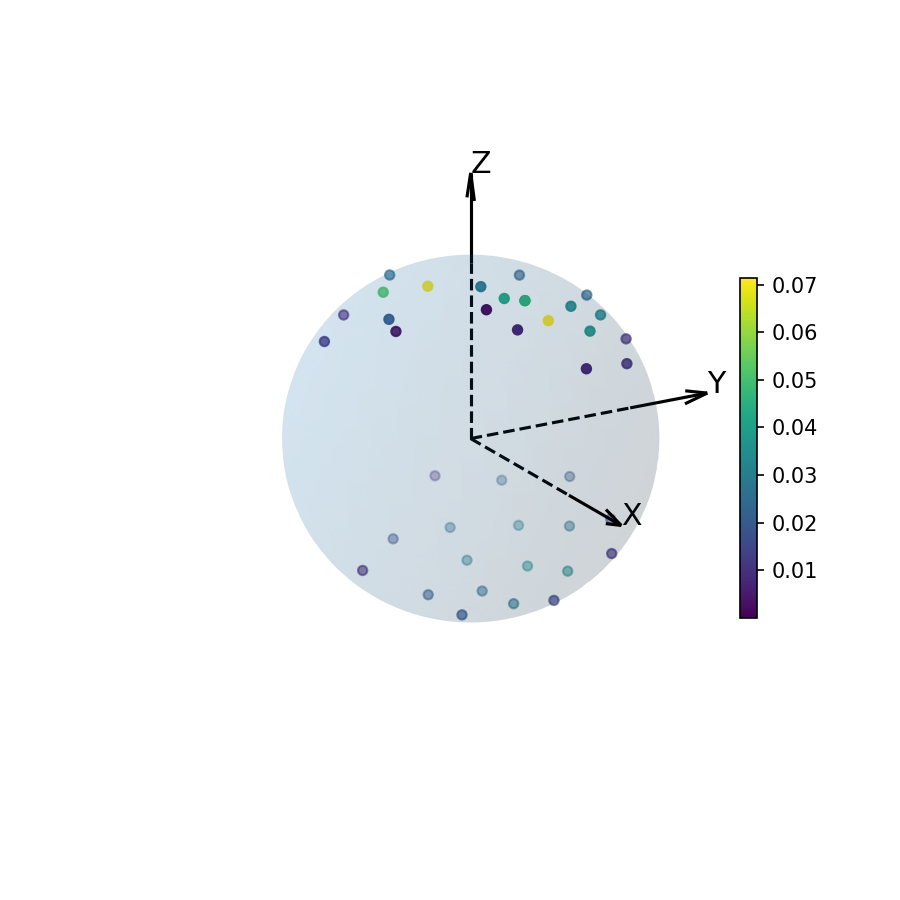}
    \caption{Left: 30000 samples on sphere; Right: EDOT of 40 points in 3D.}
    \label{fig:sphere_3d}
\end{figure}

\textbf{On Sphere.}
The underlying space $X_{\text{sphere}}$ is the unit sphere in $\mathbb R^3$.
$\mu_{\text{sphere}}$ is the pushforward of a normal mixture distribution on $\mathbb{R}^2$ by stereographic projection.
The sample set $S_{\text{sphere}}\sim\mu_{\text{sphere}}$ over $X_{\text{sphere}}$ is shown on Fig.~\ref{fig:sphere_3d} left.

Consider the (3D) Euclidean metric on $X_{\text{sphere}}$ induced by the the embedding,  
Fig.~\ref{fig:sphere_3d} (right) plots the EDOT solution with refinement for $\mu_{m}$ with $m=40$.
The resulting cell structure is shown as colored boxes.

To consider the intrinsic metric, a CW-complex is constructed.
The structure is built with a point on the
equator as a $0$-cell structure, the rest of the equator as a $1$-cell,
and the upper hemisphere and lower hemisphere as two dimension-$2$ (open) cells.
We take the upper and lower hemispheres and map them onto unit disc through
stereographic projection with respect to south pole and north pole,
respectively. Then we take the metric from spherical geometry, and rewrite
the distance function and its gradient using the natural coordinate on the unit
disc (see Supplementary for details.) 
Fig.~\ref{fig:sphere_2d} shows the refinement of EDOT on the samples (in red) and corresponding discretizations in colored points. 
More figures can be found in the Supplementary.

\begin{figure}[ht]
    \centering
    \includegraphics[scale=0.35, trim={10 10 10 10}, clip]{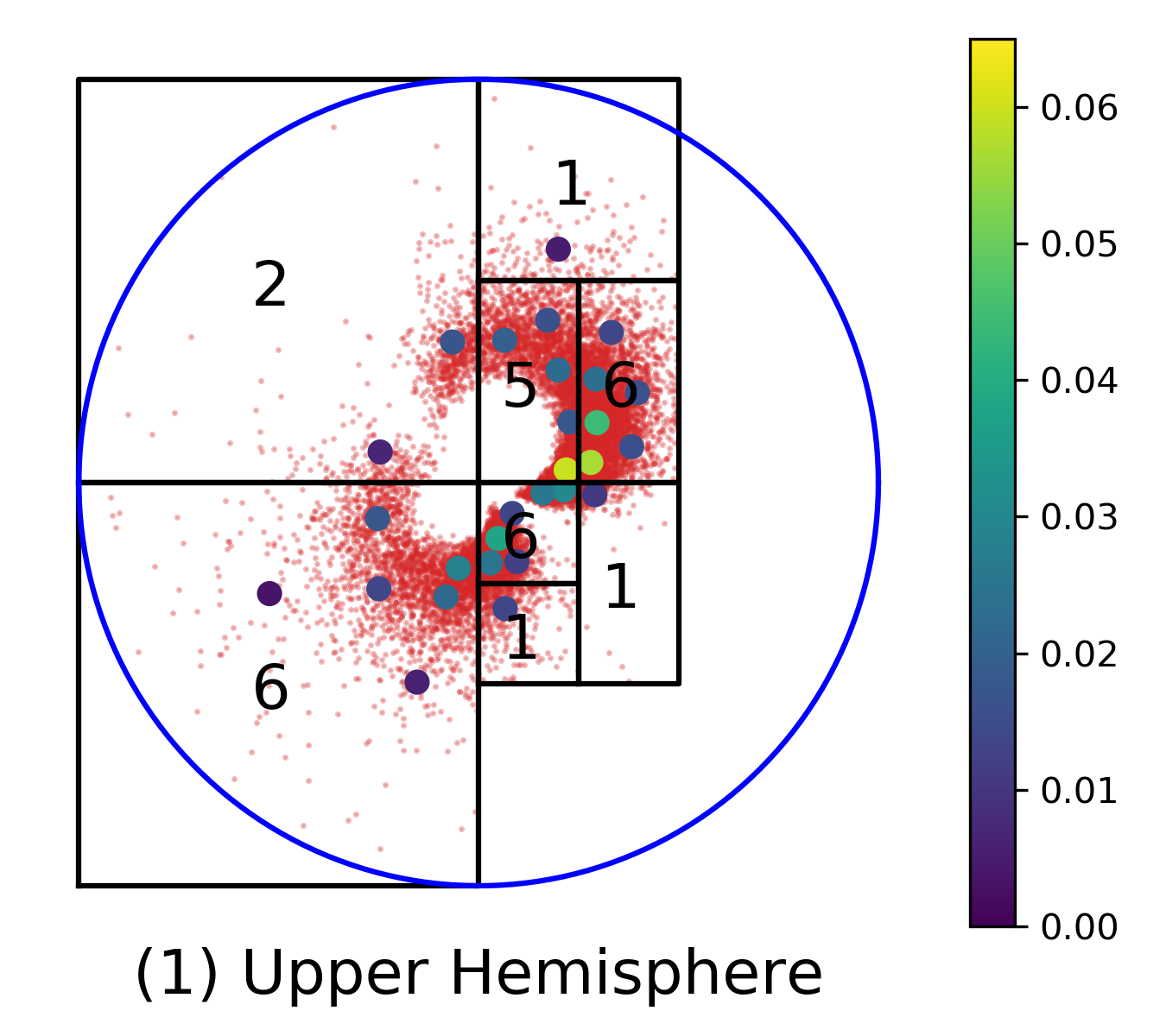}
    \includegraphics[scale=0.35, trim={10 10 10 10}, clip]{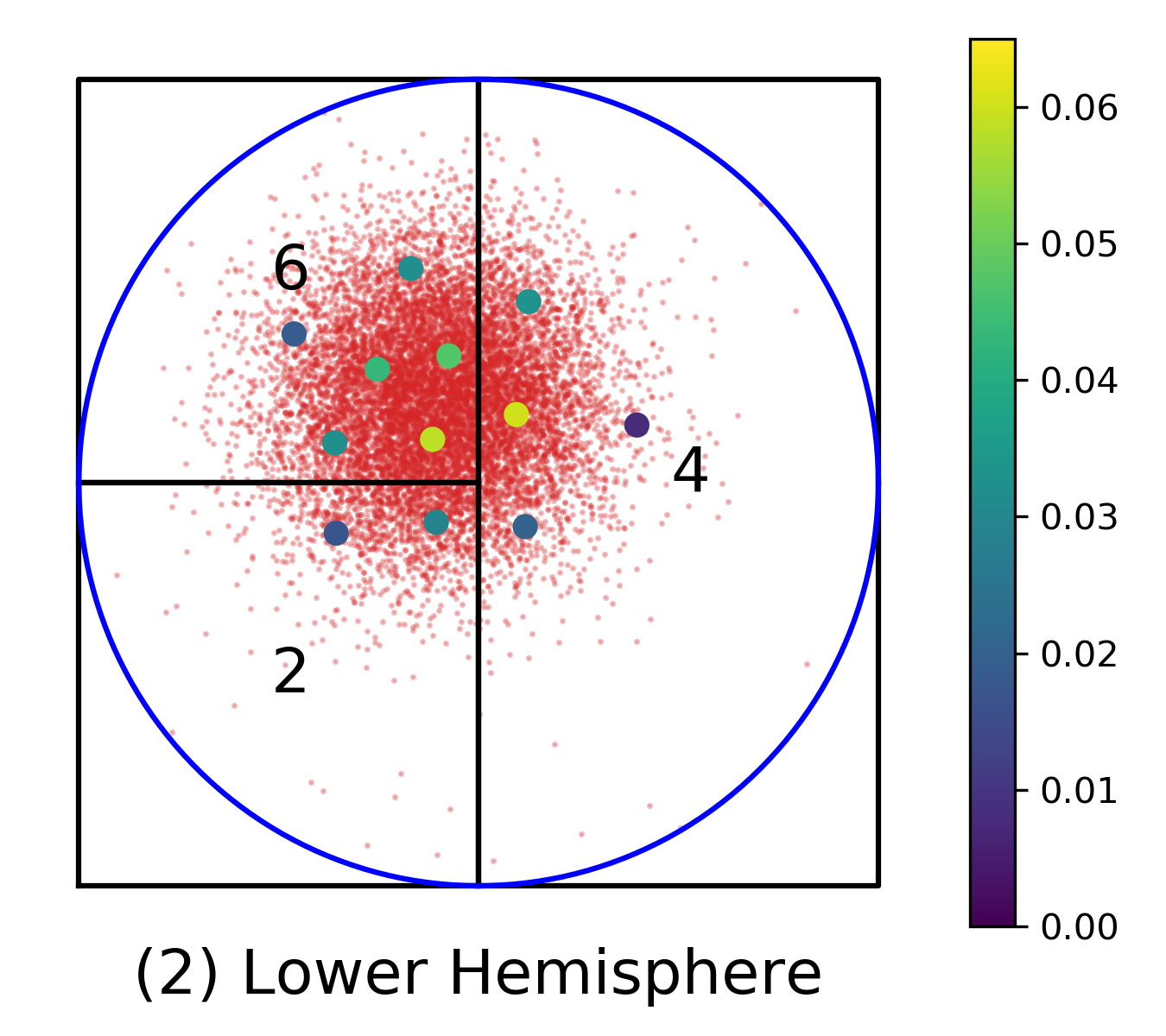}
    \caption{CW-EDOT of 40 points in 2D. Left: upper hemisphere; Right: lower hemisphere, by stereographic projections about poles.} 
    \label{fig:sphere_2d}
\end{figure}

\subsection{An Example on Transference Plan with ACR}
We now illustrate the performance of Adaptive EDOT on 
a 2D optimal transport task. Let $X=Y=[0,1]^2$, $c = d_X$ be the Euclidean distance, 
$g = d^2_X$, and the marginal $\mu$, $\nu$ be truncated normal (mixtures), where
$\mu$ has only 2 components 
and $\nu$ has 5 components. Fig.~\ref{fig:high_dim_edot} plots the McCann Interpolation of the OT plan of
between $\mu$ and $\nu$ (shown in red dots) and its discrete approximations (weights are color coded) with $m=n=30$.
With $m=n=10$, Adaptive EDOT results:
$W_{k,\zeta}^k(\mu,\mu_{10})=1.33\times10^{-2}$,
$W_{k,\zeta}^k(\nu,\nu_{10})=1.30\times10^{-2}$,
$W_{k,\zeta}^k(\gamma,\gamma_{10,10})=2.71\times10^{-2}$.
With $m = n=30$, Adaptive EDOT results, 
$W_{k,\zeta}^k(\mu,\mu_{30})=9.62\times10^{-3}$,
$W_{k,\zeta}^k(\nu,\nu_{30})=9.18\times10^{-3}$,
$W_{k,\zeta}^k(\gamma,\gamma_{30,30})=1.516\times10^{-2}$.
Where as naive sampling results:
$W_{k,\zeta}^k(\mu,\mu_{30}')=1.75\times10^{-2}$,
$W_{k,\zeta}^k(\nu,\nu_{30}')=1.58\times10^{-2}$,
$W_{k,\zeta}^k(\gamma,\gamma_{30,30}')=3.95\times10^{-2}$.
Adaptive EDOT approximated the quality of 900 naive samples with only 100 points on a 4 dimensional transference plan.

%% file: 6_complexity.tex
\section{Analysis of the Algorithms}
\label{sec:complexity}
\subsection{On Algorithm~\ref{alg:simple_EDOT}}
First, for each iteration in the minibatch SGD, 
let $N$ be the sample
(minibatch) size of $\mu_N$ for approximating $\mu$. 
Let $m$ be the size of target discretization $\mu_m$ (the output).
Further let $d$ be the dimension of $X$, $\epsilon$ be the error bound in 
Sinkhorn calculation for the entropy-regularized optimal transference plan between 
$\mu_N$ and $\mu_m$ .
\begin{figure}[ht]
    \centering
    \includegraphics[scale=0.36,
    trim={5 0 5 0}, clip]{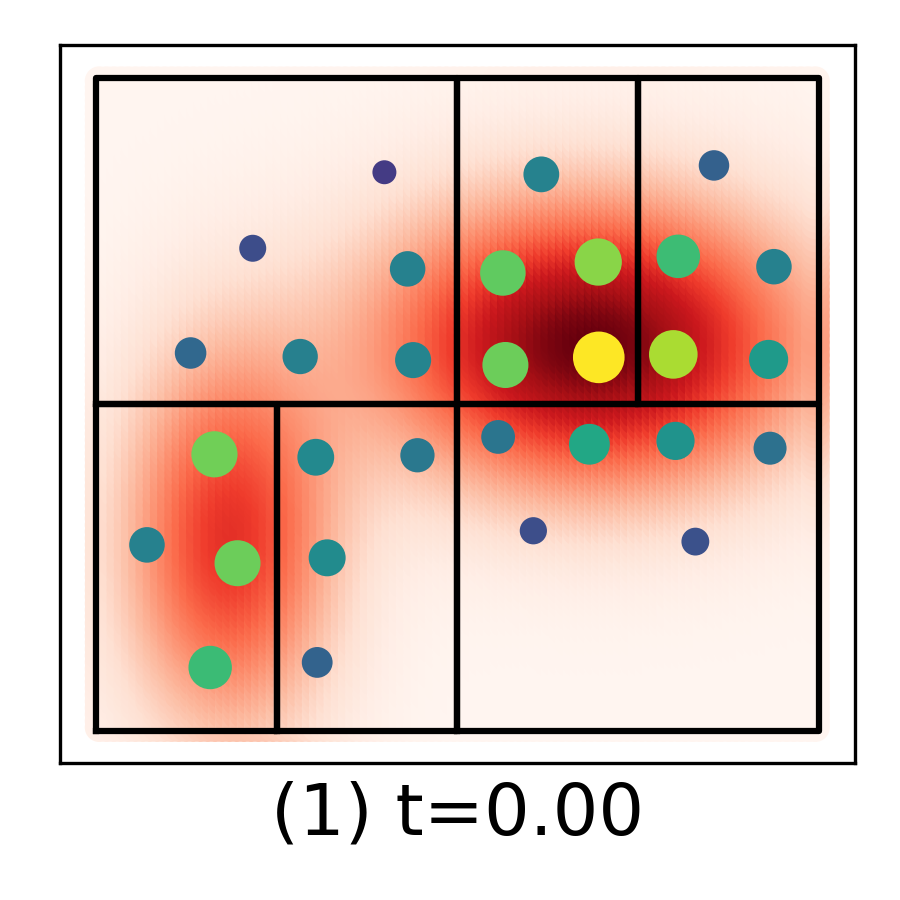}
    \includegraphics[scale=0.36,
    trim={5 0 5 0}, clip]{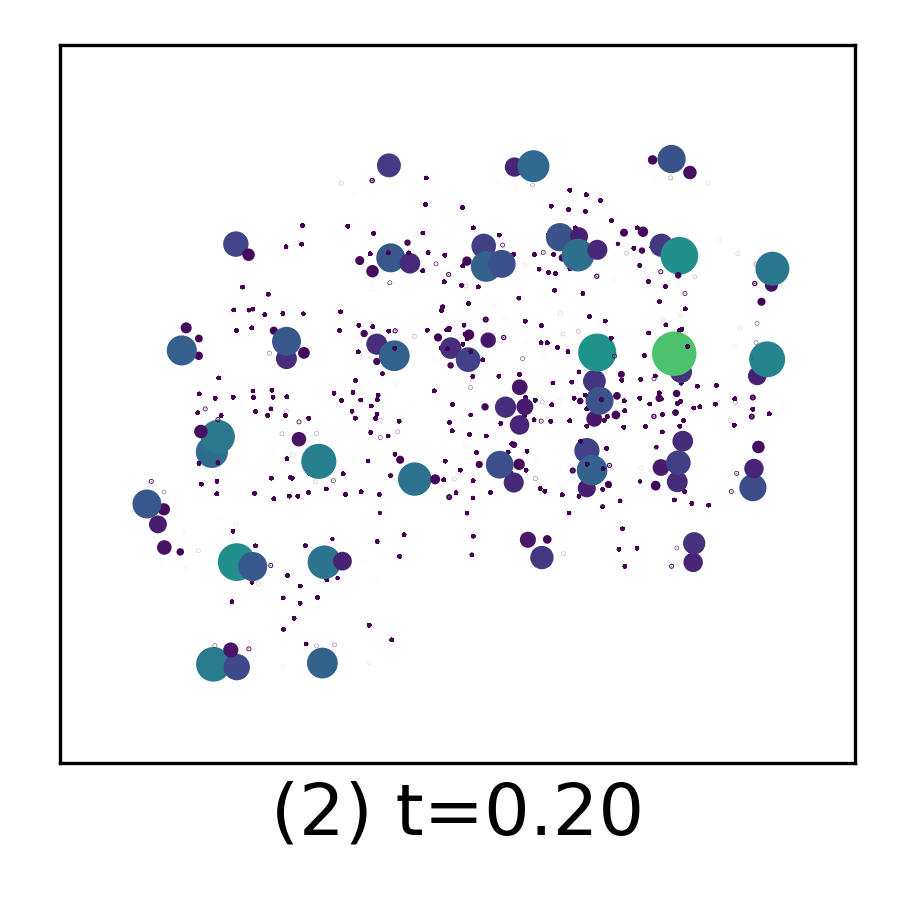}
    \includegraphics[scale=0.36,
    trim={5 0 5 0}, clip]{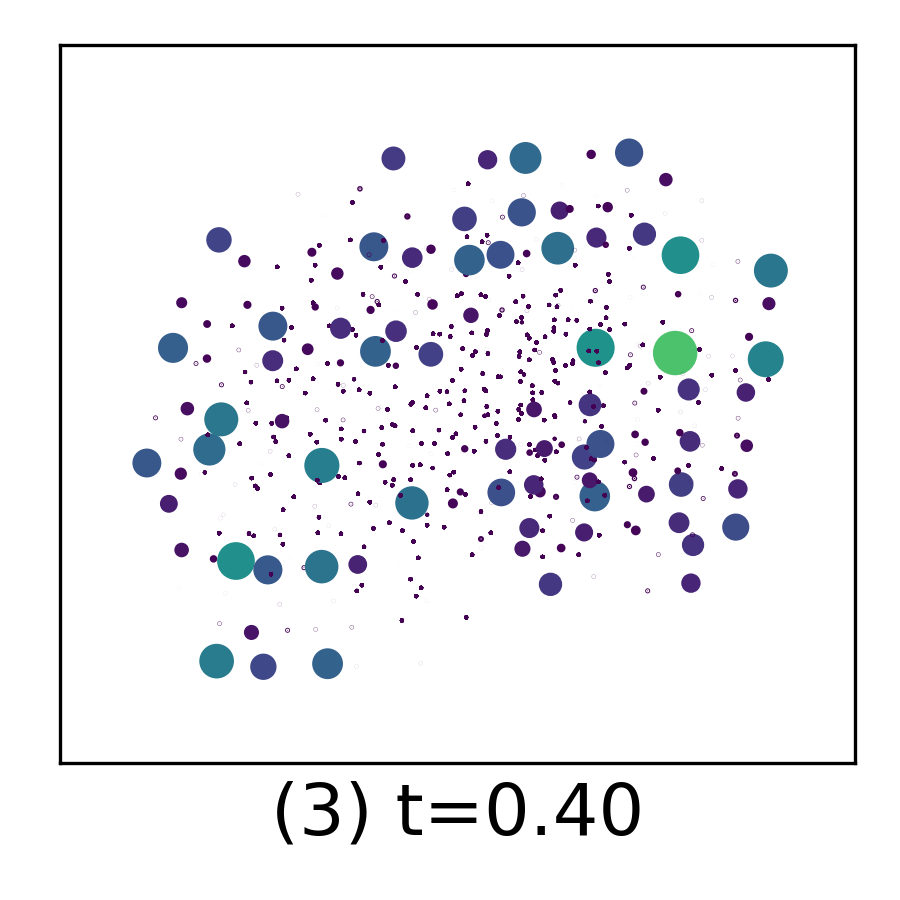}
    \includegraphics[scale=0.36,
    trim={5 0 5 0}, clip]{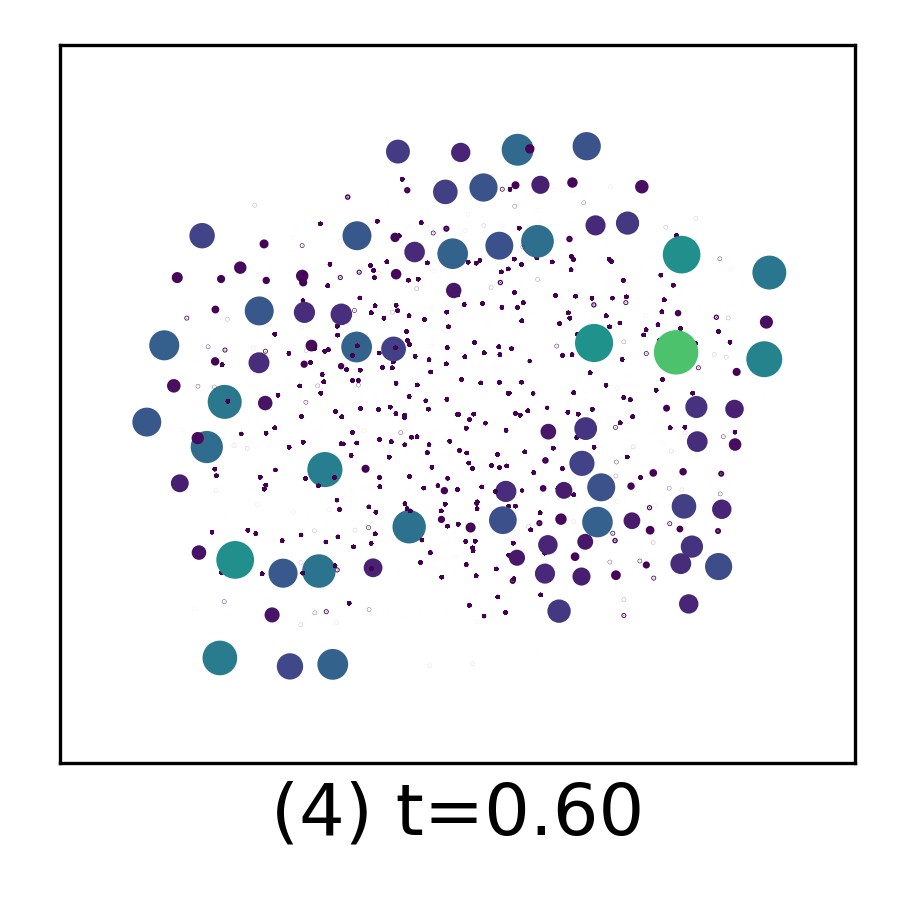}
    \includegraphics[scale=0.36,
    trim={5 0 5 0}, clip]{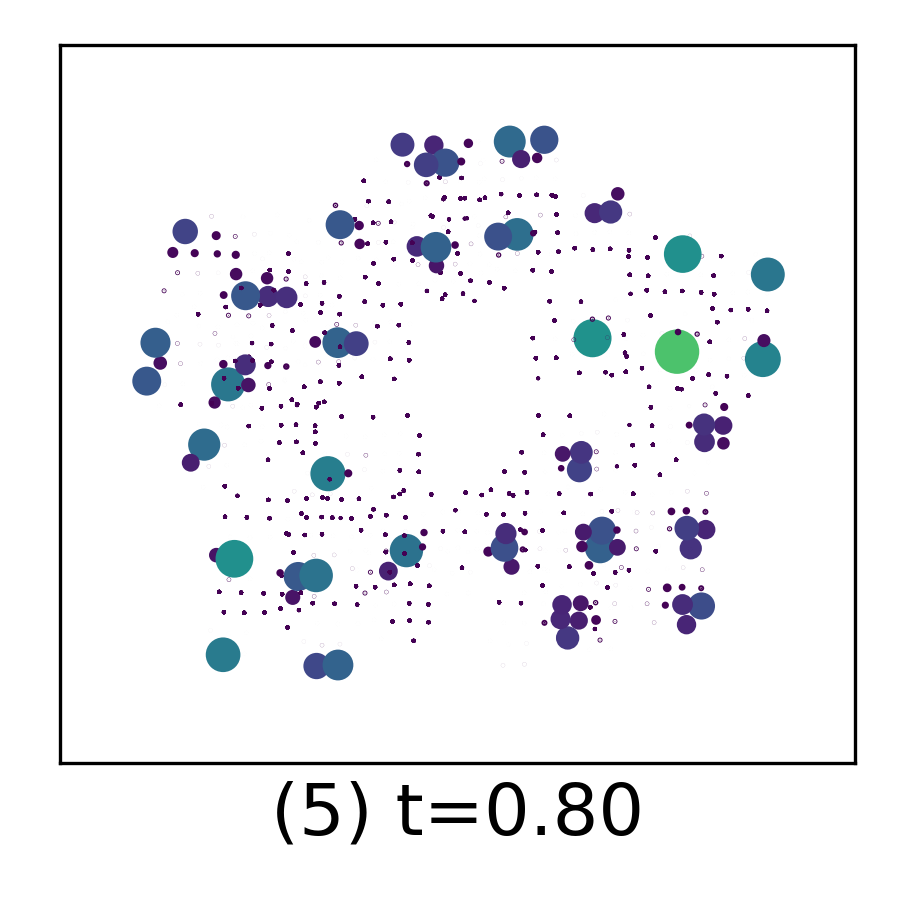}
    \includegraphics[scale=0.36,
    trim={5 0 5 0}, clip]{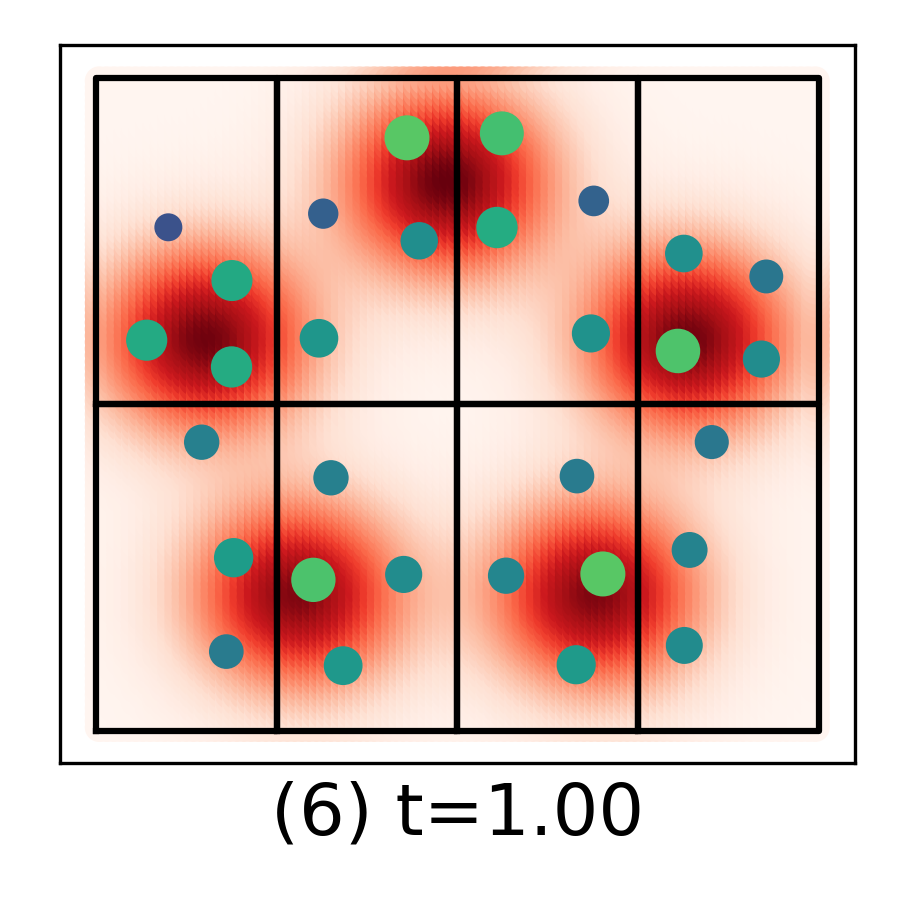}
    \caption{McCann Interpolation (See Supplementary for more slices and animation) of $30 \times 30$  transference plan.}
    \label{fig:high_dim_edot}
\end{figure}
The Sinkhorn algorithm for the positive matrix $e^{-g(x,y)/\zeta}$ (of size $N\times m$)
converges
linearly, which takes $\mathcal{O}(\log (1/\epsilon))$ steps to fall into a
region of radius $\epsilon$, contributing
$\mathcal{O}(Nm\log(1/\epsilon))$ in time complexity.
The inverse matrix $\mathbf{M}$ of
$\nabla_{(\alpha,\bar{\beta})}\mathcal{N}=-\frac{1}{\zeta}\left[
\begin{array}{cc}
    A & B \\
    B^T & D
\end{array}\right]$
(Eq.~\eqref{eq:hessian}) is taken block-wise
(see Supplementary for details): 
$$\mathbf{M}=-\zeta\left[
\begin{array}{cc}
    A^{-1}+A^{-1}BE^{-1}B^TA^{-1} &  -A^{-1}BE^{-1} \\
    -E^{-1}B^TA^{-1} & E^{-1}
\end{array}
\right]$$
where $E=D-B^TA^{-1}B$. Block
$E$ is constructed in $\mathcal{O}(Nm^2)$ and
inverted in $\mathcal{O}(m^3)$; block $A^{-1}BE^{-1}$ takes $\mathcal{O}(Nm^2)$ as
$A$ is diagonal; and the block $A^{-1}+A^{-1}BE^{-1}B^TA^{-1}$ takes
$\mathcal{O}(N^2m)$ to construct. When $m\ll N$, the time complexity in
constructing $\mathbf{M}$ is $\mathcal{O}(N^2m)$. 
From $\mathbf{M}$ to gradient of dual variables: the tensor contractions have complexity
$\mathcal{O}((N+m)^2md)$.
Finally, to get the gradient, the 
complexity is dominated by the second term 
of Eq.~\eqref{eq:gradient_W_y},
which is a contraction between a matrix $Nm$ (i.e., $g\dd\pi(x)$) with tensors
of sizes $Nmd$ and $m^2d$ (two gradients on dual variables
$\alpha$, $\beta$) along $N$ and $m$
respectively. Thus the final step contributes $\mathcal{O}((N+m)md)$.

The time complexity of increment steps in SGD turns out to be
$\mathcal{O}(md)$.
Therefore, for $L$ steps of minibatch SGD, the time complexity is
$\mathcal{O}((N+m)^2mdL + NmL \log(1/\epsilon))$.

For space complexity, Sinkhorn algorithm (which can be done in position $\mathcal{O}(Nm)$) is the
only iterative computation in a 
single SGD step, and between two SGD steps, only the resulting distribution is
passed to next step. Therefore, the space complexity is $\mathcal{O}((N+m)^2)$
coming from the $\mathbf{M}$, others are at most of size $\mathcal{O}(m(N+m))$.

From the above discussion, we can further see that when
Algorithm~\ref{alg:refine_dfs} is applies, the total complexities (both in time
and space) are reduced as the magnitudes of both $N$ and $m$ are much smaller in each refined cell.

Finally, the convergence rate affects the time complexity.
Since minibatch SGD is applied, the convergence rate
depends on the choice of the initial distribution, and is also
affected by the size $m$, $N$, and the diameter of the cell.

\subsection{On Algorithm~\ref{alg:refine_dfs}}

The procedure of dividing sample set $S$ into subsets in Algorithm~\ref{alg:refine_dfs} 
is similar to Quicksort, thus the space and time complexities are similar.
The similarity comes from the binary divide-and-conquer structure as well as
that each split action is based on comparing each sample with a target (the
``mid'' in Algorithm~\ref{alg:refine_dfs}).

The time complexity of Algorithm~\ref{alg:refine_dfs} is $\mathcal{O}(N_0\log N_0)$
in best and average case, $\mathcal{O}(N_0^2)$ in worst case (which would happen
when each cell contains only one sample). The time complexity depends only on
the distribution, not on the order the samples are stored.

The space complexity is $\mathcal{O}(N_0d+m)$, or simply $\mathcal{O}(N_0d)$ as $m\ll
N_0$. Furthermore, if we use swap to make the split of sample size ``in
position'' (as in Quicksort), we may only need space for one copy of set $S$
(of size $N_0d$), those in DFS stack and output stack (proportional to $m$) and
some fixed size space for calculations.

The applications of Algorithm~\ref{alg:simple_EDOT} on the cells are independent
from each other, making the computation parallelizable.
And the post-processing
of such parallelization takes only $m$ multiplications and a concatenation to construct $\mu_m$.

\subsection{Comparison with Naive Sampling}

After having a size-$m$ discretization on $X$ and a size-$n$ discretization on
$Y$, the EOT solution (Sinkhorn algorithm) has time compleixty
$\mathcal{O}(mn\log(1/\epsilon))$. In EDOT, two discretization problems
must be solved before applying Sinkhorn, while the naive sampling requires nothing but
sampling. 

According to the previous analysis, solving a single continuous EOT problem
using size-$m$ EDOT method directly may
result in higher time complexity than naive sampling with even a larger sample
size (than $m$), when the cost function $c:X\times Y\rightarrow\mathbb{R}$ is
completely known. However, in real applications, the cost function may be from
real world experiments (or from extra computations) done for each pair $(x,y)$
in the discretization, thus 
the size of discretized distribution is critical for cost-control, but the
distance function $d_X$ and $d_Y$ usually come along with the spaces $X$ and
$Y$ and are easy to calculate. Another scenario of application of EDOT is that the
marginal distributions $\mu_X$ and $\nu_Y$ are fixed for different cost
functions, then discretizations can be reused, thus the cost of discretization
is one-time and improvement it brings accumulates in each repeat.

%% file: 7_conclusion_etc.tex
\section{Conclusion}
\label{sec:conclusion}
We developed methods for efficiently approximating OT couplings with fixed size $m\times n$ approximations. 
We provided bounds on the relationship between the discrete approximation and the original continuous
problem. 
We implemented two algorithms and demonstrated their efficacy as compared to naive sampling and analyzed computational complexity. 
Our approach provides a new approach to efficiently computing OT plans. 

%% file: SI_1_proofs.tex
\section{Proof of Proposition~\ref{prop:3w_ineq}}

\begin{repprop}{prop:3w_ineq}
When $X$ and $Y$ are two compact spaces and $c$ is $\mathcal{C}^{\infty}$, there exists a constant
$C_1 \in \mathbb{R}^{+}$ such that 
\begin{align}
    \max\{ & W_{k}^k(\mu, \mu_m), \! W_{k}^k(\nu,\nu_n)\}\!\stackrel{(i)}{\leq}\! W_{k,\zeta}^k\!(\gamma_{\lambda}(\mu,\nu),\! \gamma_{\lambda}(\mu_m,\nu_n)) \nonumber\\
     & \stackrel{(ii)}{\leq} C_1[W_{k,\zeta}^k(\mu, \mu_m)+W_{k,\zeta}^k(\nu,\nu_n)].
\end{align}
\end{repprop}

\begin{proof}
We will adopt notations: $Z =  X \times Y$ and $z_i = (x_i, y_i) \in Z$.

For inequality~(i), without loss assume that 
\[\max\{ W_{k}^k(\mu, \mu_m), W_{k}^k(\nu,\nu_n)\} = W_{k}^k(\mu, \mu_m)\]

Denote the optimal $\pi_{Z} \in \Pi(\gamma_{\lambda}(\mu,\nu), \gamma_{\lambda}(\mu_m,\nu_n))$ that achieves $ W_{k}^k (\gamma_{\lambda}(\mu,\nu), \gamma_{\lambda}(\mu_m,\nu_n))$ by 
$\pi^*_{Z}$ similarly for $\pi^*_{X}$. Then we have:
\begin{align*}
   W_{k}^k & (\gamma_{\lambda}(\mu,\nu), \gamma_{\lambda}(\mu_m,\nu_n)) \nonumber\\
   & =\int_{Z^2} d^k_Z(z_1, z_2) \dd \pi^*_{Z} \geq A \int_{Z^2} d^k_{X}(x_1, x_2) \dd \pi^*_Z \nonumber\\
   & =A \int_{X^2}\! d^k_{X}(x_1, x_2) \int_{Y^2}\! \dd \pi^*_Z \stackrel{(a)}{=}A\! \int_{X^2}\! d^k_{X}(x_1, x_2) \dd \pi'_X \nonumber\\
   &  \stackrel{(b)}{\geq} A \int_{X^2} d^k_{X}(x_1, x_2) \dd \pi^*_X = W_{k}^k(\mu, \mu_m)
\end{align*}
Here $\pi'_X \in \Pi(\mu, \mu_m)$, eq~(a) holds since $\pi^*_Z \in \Pi(\gamma_{\lambda}(\mu,\nu), \gamma_{\lambda}(\mu_m,\nu_n))$ 
and ineq~(b) holds since $\pi^*_X$ is the optimal choice.

For inequality~(ii),
we will use the following to simplify the notations, 
$$\dd\gamma \otimes \dd \gamma_{mn} \coloneqq \dd\gamma_{\lambda} (\mu,\nu) \otimes \dd \gamma_{\lambda}(\mu_m,\nu_n)$$
$d^{k}_{Z} \coloneqq d^{k}_{Z}(z_1, z_2),  d^{k}_{X} \coloneqq d^{k}_{X}(x_1, x_2), d^{k}_{Y} \coloneqq d^{k}_{Y}(y_1, y_2)$

\begin{align*}
    &W_{k,\zeta}^k(\gamma_{\lambda}(\mu,\nu), \gamma_{\lambda}(\mu_m,\nu_n)) =  \int_{Z^2} d^{k}_{Z}(z_1, z_2) \dd \pi^*_{z,\zeta}\\
    & \stackrel{(a)}{=} \int_{Z^2} d^{k}_{Z}\cdot \exp(\frac{\alpha(z_1)+\beta(z_2)-d^{k}_{Z}(z_1, z_2)}{\zeta})\dd\gamma \otimes \dd \gamma_{mn}\\
    & \stackrel{(b)}{\leq} A \int_{Z^2} (d^{k}_{X}+d^{k}_{Y})\cdot \exp(\frac{\alpha(z_1)+\beta(z_2)-d^{k}_{Z}}{\zeta})\dd\gamma \otimes \dd \gamma_{mn}\\
    & \stackrel{(c)}{\leq} C_1 \int_{Z^2} (d^{k}_{X}+d^{k}_{Y})\cdot \exp(\frac{-d^{k}_{Z}}{\zeta})\dd\gamma \otimes \dd \gamma_{mn}\\
    & \stackrel{(d)}{\leq} C_1 [\int_{Z^2} d^{k}_{X}\cdot \exp(\frac{-d^{k}_{X}}{\zeta}) + d^{k}_{Y}\cdot \exp(\frac{-d^{k}_{Y}}{\zeta})\dd\gamma \otimes \dd \gamma_{mn}]\\
    & \stackrel{(e)}{=} C_1 [\int_{X^2} d^{k}_{X}\cdot \exp(\frac{-d^{k}_{X}}{\zeta}) \dd\mu \otimes \dd \mu_{m} \\
    &\hspace{0.4in} + \int_{Y^2} d^{k}_{Y}\cdot \exp(\frac{-d^{k}_{Y}}{\zeta})\dd\nu \otimes \dd \nu_{n}]\\
    & \stackrel{(f)}{\leq} C_1 \int_{X^2} d^{k}_{X}\cdot \exp(\frac{s(x_1)+t(x_2)-d^{k}_{X}}{\zeta}) \dd\mu \otimes \dd \mu_{m}\\
    & \hspace{0.4in}+ C_1\int_{Y^2} d^{k}_{Y}\cdot \exp(\frac{s'(y_1)+t'(y_2)-d^{k}_{Y}}{\zeta})\dd\nu \otimes \dd \nu_{n}\\
    & = C_1[W_{k,\zeta}^k(\mu, \mu_m)+W_{k,\zeta}^k(\nu,\nu_n)] \\
\end{align*}

Justifications for the derivations:

(a) Based on the dual formulation, it is shown in \cite{aude2016stochastic}[Proposition 2.1] that 
for $\zeta>0$, there exist $\alpha(z_1), \beta(z_2) \in \mathcal C (Z)$ such that 
$\dd \pi^*_{z,\zeta} =\exp(\frac{\alpha(z_1)+\beta(z_2)-d^{k}_{Z}(z_1, z_2)}{\zeta})\dd\gamma \otimes \dd \gamma_{mn} $;

(b) Inequality (ii) of eq~\eqref{eq:d_xy};

(c) According to \cite{genevay2019sample}[Theorem~2], when $X, Y$ are compact and $c$ is smooth, $\alpha,\beta$ are uniformly bounded; 
moreover both $d^{k}_{X}$ and $d^{k}_{Y}$ are uniformly bounded by the diameter of $X$ and $Y$ respectively, hence constant $B$ exists;

(d) Inequality (ii) of eq~\eqref{eq:d_xy};

(e) $\gamma_{\lambda} (\mu,\nu) \in \Pi(\mu, \nu)$ and $\gamma_{\lambda}(\mu_m,\nu_n)\in \Pi(\mu_m, \nu_n)$;

(f) Similarly as in (a), for $\zeta>0$, there exist $s(x_1), t(x_2) \in \mathcal C (X)$ and  $s'(y_1), t'(y_2) \in \mathcal C (Y)$ such that 
$\exp(\frac{-d^{k}_{X}}{\zeta}) \dd\mu \otimes \dd \mu_{m} = \dd \pi^*_{X}$ and $\exp(\frac{-d^{k}_{Y}}{\zeta})\dd\nu \otimes \dd \nu_{n} = \dd \pi^*_{Y}$.
Moreover, $\int_{X^2} d^{k}_{X}\cdot \exp(\frac{s(x_1)+t(x_2)}{\zeta})\dd\mu \otimes \dd \mu_{m} \geq 0$ and 
$\int_{Y^2} d^{k}_{Y}\cdot \exp(\frac{s'(y_1)+t'(y_2)}{\zeta})\dd\nu \otimes \dd \nu_{n} \geq 0$.
\end{proof}

%% file: SI_2_differentiation.tex
\section{Gradient of $W_{k,\zeta}^k$}
\label{sec:supp:gradient}

\subsection{The Gradient}
\label{subsec:supp:gradient}

Following the definitions and notations in Section~2 and Section~3 of the paper,
we calculate the gradient of $W_{k,\zeta}^k(\mu,\mu_m)$ about parameters
of $\mu_m:=\sum_{i=1}^{m}w_i\delta_{y_i}$ in detail.

$W_{k,\zeta}^k(\mu,\mu_m)=\int_{X^2}g(x,y)\dd{\gamma_{\zeta}(x,y)}$ where
\begin{equation}
  \label{eq:supp:primal_problem}
  \small \gamma_{\zeta}=\argmin_{\gamma\in\Pi(\mu,\mu_m)}\int_{X^2}\!\!\!g(x,y)\dd{\gamma(x,y)}+\zeta\mathrm{KL}(\gamma||\mu\otimes\mu_m).
\end{equation}
  
Let $\alpha\in L^\infty(X)$ and $\beta\in\mathbb{R}^{m}$,
denote ${\beta}=\sum_{i=1}^{m}\beta_i\delta_{y_i}$, let
\begin{align}
  &\mathcal{F}(\gamma;\mu,\mu_m,\alpha,\beta)\nonumber\\
  :=&\int_{X^2}\!\!\!g(x,y)\dd{\gamma(x,y)}+
      \zeta\mathrm{KL}(\gamma||\mu\otimes\mu_m)\nonumber\\
  &+\int_X\alpha(x)\left(\int_X\dd\gamma(x,y)-\dd\mu_m(y)\right) \nonumber\\
  &+ \int_X{\beta}(y)\left(\int_{X}\dd\gamma(x,y)-\dd\mu(x)\right)\nonumber\\
  =&\int_{X^2}\!\!\!g(x,y)\dd{\gamma(x,y)}+
     \zeta\mathrm{KL}(\gamma||\mu\otimes\mu_m)\nonumber\\
  &+\int_X\alpha(x)\left(\sum_{i=1}^{m}\dd\gamma(x,y_i)-\dd\mu_m(y_i)\right)\nonumber\\
  &+ \sum_{i=1}^{m}\beta_i\left(\int_{X}\dd\gamma(x,y_i)-\dd\mu(x)\right)
\end{align}

Since $\gamma$ on second component $X$ is discrete and supported on $\{y_i\}_{i=1}^m$,
we may denote $\dd\gamma(x,y_i)$ by $\dd\pi_i(x)$, thus
\begin{align}
  &\mathcal{F}(\gamma;\mu,\mu_m,\alpha,\beta)\nonumber\\
  =&\int_{X}\sum_{i=1}^{m}g(x,y_i)\dd{\pi_i(x)}\nonumber\\
    &+ \zeta\sum_{i=1}^{m}\int_{X}\left(\ln\dfrac{\dd\pi_i(x)}{\dd\mu(x)}-\ln(\mu_m(y_i))\right)\dd\pi_i(x)\nonumber\\
  &+\int_X\alpha(x)\left(\sum_{i=1}^{m}\dd\pi_i(x)-\dd\mu_m(y_i)\right)\nonumber\\
  &+ \sum_{i=1}^{m}\beta_i\left(\int_{X}\dd\pi_i(x)-\dd\mu(x)\right)
\end{align}

Then the Fenchel duality of Problem~\eqref{eq:supp:primal_problem} is
\begin{align}
  \label{eq:supp:fenchel_dual}
  &\mathcal{L}(\mu,\mu_m;\alpha,\beta)\nonumber\\
  =&\int_{X}\alpha(x)\dd\mu(x)+\sum_{i=1}^{m}\beta_iw_i\nonumber\\
  &-\zeta\int_X\sum_{i=1}^{m}e^{(\alpha(x)+\beta_i-g(x,y_i))/\zeta}w_i\dd\mu(x).
\end{align}
Let $\alpha^\ast$, $\beta^\ast$ be the argmax of the Fenchel dual \eqref{eq:supp:fenchel_dual},
the primal is solved by
$\dd\gamma^\ast(x,y_i)=e^{(\alpha^\ast(x)+\beta^\ast_i-g(x,y_i))/\zeta}w_i\dd\mu(x).$
To make the solution unique, we restrict the freedom of solution (where we see
that $\mathcal{L}(\mu,\mu_m;\alpha,\beta)=\mathcal{L}(\mu,\mu_m;\alpha+t,\beta-t)$
for any $t\in\mathbb{R}$). We use the condition $\beta_m=0$ to narrow the
choices down to only one, and denote the dual variable $\beta$ having the property
by $\overline{\beta}$ and $\overline{\beta}^\ast$.

We first calculate $\nabla_{w_i, y_i}\mathcal{L}(\mu,\mu_m;\alpha^\ast,\overline{\beta}^\ast)$
with $\alpha^\ast$ and $\overline{\beta}^\ast$ as functions of $\mu_m$.
(From the paper)

\begin{align}\small
  \label{eq:supp:gradient_W_w} &\!\!\!\!\dfrac{\partial \mathcal{L}}{\partial w_i}
=\int_Xg(x,y_i)E^\ast_i(x)\dd\mu(x)+\nonumber \\
&\!\!\!\dfrac{1}{\zeta}\int_X\sum_{j=1}^{n}g(x,y_j) \left(\dfrac{\partial
\alpha^\ast\!(x)}{\partial w_i} + \dfrac{\partial \beta^\ast_j}{\partial w_i}\right) w_j
E^\ast_j(x)\dd\mu(x).
\end{align}
\begin{align}\small
  \label{eq:supp:gradient_W_y} &\!\!\!\!\nabla_{y_i}\!\mathcal{L}
=\!\int_X\!\!\!\nabla_{y_i}g(x,y_i)\left(\!1-\dfrac{g(x,y_i)}{\zeta}\!\right)
\!E^\ast_i(x)w_i\dd\mu(x)+\nonumber \\
&\!\!\frac{1}{\zeta}\!\!\int_X\sum_{j=1}^{n}g(x,y_j)
\!\left(\nabla_{y_i}\!\alpha^\ast\!(x)\!+\! \nabla_{y_i}\!\beta^\ast_j\right) w_j
E^\ast_j(x)\dd\mu(x).
\end{align}

Next, we calcuate the derivatives of $\alpha^\ast$ and $\overline{\beta}^\ast$,
by finding their defining equation and then using Implicit Function Theorem.

The optimal solution to the dual variables $\alpha^\ast$, $\beta^\ast$
is obtained by solving the stationary state equation $\nabla_{\alpha,\beta}\mathcal{L}=0$.
The derivatives are taken in the sense of Fr\'echet derivative.
The Fenchel dual function on $\alpha$, $\overline{\beta}$,
has its domain and codomain
$\mathcal{L}(\mu,\mu_m;\cdot,\cdot):L^\infty(X)\times\mathbb{R}^{m-1}\rightarrow\mathbb{R}$,
the derivatives are
\begin{align}
  \label{eq:supp:d_fenchel_a}
  &\nabla_{\alpha}\mathcal{L}(\mu,\mu_m;\alpha,\overline{\beta})\nonumber\\
  =&\int_X\left(1-\sum_{i=1}^mw_iE_i(x)\right)(\ \cdot\ )\dd\mu(x),
\end{align}
\begin{align}
  \label{eq:supp:d_fenchel_b}
  \dfrac{\partial}{\partial\overline{\beta}_i}\mathcal{L}(\mu,\mu_m;\alpha,\overline{\beta})
  =w_i\left(1-\int_XE_i(x)\dd\mu(x)\right)
\end{align}
where $E_i(x)=e^{(\alpha(x)+\beta_i-g(x,y_i))/\zeta}$ is defined as in the paper,
and $\nabla_{\alpha}\mathcal{L}(\mu,\mu_m;\alpha,\overline{\beta})\in (L^\infty(X))^\vee$
(as a linear functional),
$\dfrac{\partial}{\partial\overline{\beta}_i}\mathcal{L}(\mu,\mu_m;\alpha,\overline{\beta})\in\mathbb{R}$.
Next, we need to show $\mathcal{L}$ is differentiable 
in the sense of Fr\'echet derivative, i.e.,
\begin{align}
  \lim_{||h||\rightarrow0}\dfrac{1}{||h||}&\left( 
  \mathcal{L}(\mu,\mu_m;\alpha+h,\overline{\beta})-
  \mathcal{L}(\mu,\mu_m;\alpha,\overline{\beta})\right.\nonumber\\
  &\left.-\nabla_{\alpha}\mathcal{L}(\mu,\mu_m;\alpha,\overline{\beta})(h)
  \right)=0.
\end{align}
By definition of $\mathcal{L}$ (we write $\mathcal{L}(\alpha)$ for
$\mathcal{L}(\mu,\mu_m;\alpha,\overline{\beta})$),
\begin{align}
  \label{eq:supp:frechet_remainder}
  &\mathcal{L}(\alpha+h)-
  \mathcal{L}(\alpha)
  -\nabla_{\alpha}\mathcal{L}(\alpha)(h)\nonumber\\
  =&\int_X\!h(x)\dd\mu(x)\!-\!\zeta\!\int_X\sum_{i=1}^{m}\left(e^{h(x)/\zeta}\!-\!1\right)w_iE_i(x)\dd\mu(x)\nonumber\\
  &-\int_X\left(1-\sum_{i=1}^mw_iE_i(x)\right)h(x)\dd\mu(x)\nonumber\\
  =&\zeta\int_X\sum_{i=1}^{m}\left(
     1+\dfrac{h(x)}{\zeta}-e^{h(x)/\zeta}
     \right)E_i(x)w_i\dd\mu(x)\nonumber\\
  =&\zeta\int_X\left(\sum_{k=2}^{\infty}\dfrac{1}{k!}\dfrac{h(x)^k}{\zeta^k}\right)\sum_{i=1}^{m}
     E_i(x)w_i\dd\mu(x),
\end{align}
the last equality is from the Taylor expansion of exponential function.
Consider that $||h||_{\infty}=\mathrm{ess\ sup\ }|_{x\in X}h(x)|$ the
essential supremum of $|h(x)|$ for $x\in X$ given measure $\mu$.

Denote $\mathcal{N}:=\nabla_{\alpha,\overline{\beta}}\mathcal{L}$,
\begin{align}
  \label{eq:supp:frechet_L_final}
   &\dfrac{1}{||h||}\left(\mathcal{L}(\alpha+h)-
     \mathcal{L}(\alpha)
     -\nabla_{\alpha}\mathcal{L}(\alpha)(h)\right)\nonumber\\
  \le&\dfrac{\zeta}{||h||}\int_X\left(\sum_{k=2}^{\infty}\dfrac{1}{k!}\dfrac{|h(x)|^k}{\zeta^k}\right)\sum_{i=1}^{m}
       E_i(x)w_i\dd\mu(x)\nonumber\\
  \le&\dfrac{\zeta}{||h||}\int_X\left(\sum_{k=2}^{\infty}\dfrac{1}{k!}\dfrac{||h||^k}{\zeta^k}\right)\sum_{i=1}^{m}
       E_i(x)w_i\dd\mu(x)\nonumber\\
  =&\zeta\left(\sum_{k=2}^{\infty}\dfrac{1}{k!}\dfrac{||h||^{k-1}}{\zeta^k}\right)\int_X\sum_{i=1}^{m}
     E_i(x)w_i\dd\mu(x)\nonumber\\
  =&\zeta\left(\sum_{k=2}^{\infty}\dfrac{1}{k!}\dfrac{||h||^{k-1}}{\zeta^k}\right)  
\end{align}
Therefore,
\begin{align}
  \label{eq:supp:lim_frechet}
     \lim_{||h||\rightarrow0}\dfrac{1}{||h||}\left(\mathcal{L}(\alpha+h)-
     \mathcal{L}(\alpha)
     \!-\!\nabla_{\alpha}\mathcal{L}(\alpha)(h)\right)=0,
\end{align}
which shows that the expression of $\nabla_\alpha\mathcal{L}(\alpha)$ in
Eq.~\eqref{eq:supp:d_fenchel_a} gives the correct Fr\'echet derivative.
Note here that $\alpha\in L^\infty(X)$ is critical in Eq.~\eqref{eq:supp:frechet_L_final}.

Let
$\mathcal{N}:=\nabla_{\alpha,\overline{\beta}}\mathcal{L}$
values in $(L^\infty(X))^\vee\times\mathbb{R}^{m-1}$,
then $\mathcal{N}=0$ defines $\alpha^\ast$ and $\overline{\beta}^\ast$,
which makes it possible to differentiate them about $\mu_m$ using Implicit
Function Theorem for Banach spaces. From now on, $\mathcal{N}$ take values
at $\alpha=\alpha^\ast$, $\overline{\beta} = \overline{\beta}^\ast$, i.e.,
the marginal conditions on $\dd\pi_i(x)=w_iE_i(x)\dd\mu(x)$ hold.

Thus we need $\nabla_{\alpha,\overline{\beta}}\mathcal{N}$ and $\nabla_{w_i,y_i}\mathcal{N}$
calculated, and prove that $\nabla_{\alpha,\overline{\beta}}\mathcal{N}$ is invertible
(and give the inverse).

It is necessary to make sure which form  $\nabla_{\alpha,\overline{\beta}}\mathcal{N}$
is in according to Fr\'echet derivative.
Start from the map $\mathcal{N}(\mu,\mu_m;\cdot,\cdot):
(L^\infty(X))\times\mathbb{R}^{m-1}\rightarrow(L^\infty(X))^\vee\times(\mathbb{R}^{m-1})^\vee$
where $\mathbb{R}^{m-1}$ is isomorphic to its dual Banach space $(\mathbb{R}^{m-1})^\vee$.
Then $\nabla_{\alpha,\overline{\beta}}\mathcal{N}\in\mathrm{Hom}^b_{\mathbb{R}}
(L^\infty(X)\times\mathbb{R}^{m-1},
(L^\infty(X))^\vee\times(\mathbb{R}^{m-1})^\vee)$, where $\mathrm{Hom}^b$ represents the set of bounded linear operators.
Moreover, recall that $(\cdot)\otimes A$ is the left adjoint functor
of $\mathrm{Hom}^b_{\mathbb{R}}(A,\cdot)$, then
for $\mathbb{R}$-vector spaces,
$\mathrm{Hom}^b_{\mathbb{R}}
(L^\infty(X)\times\mathbb{R}^{m-1},
(L^\infty(X))^\vee\times(\mathbb{R}^{m-1})^\vee)\cong
\mathrm{Hom}^b_{\mathbb{R}}\left(
  (L^\infty(X)\times\mathbb{R}^{m-1})^{\otimes2},\mathbb{R}\right)$.
Thus, we can write $\nabla_{\alpha,\overline{\beta}}\mathcal{N}$
in terms of a bilinear form on vector space
$L^\infty(X)\times\mathbb{R}^{m-1}$.

From the expression of $\mathcal{N}$, we may differentiate
(similarly as the calculations \eqref{eq:supp:frechet_remainder} to \eqref{eq:supp:lim_frechet}):
\begin{align}
  \label{eq:supp:sim_part_a}
  \nabla_{\alpha}\mathcal{N}
  =&\left(
     -\dfrac{1}{\zeta}\int_X(\cdot)(-)\sum_{i=1}^{m}w_iE_i(x)\dd\mu(x)\right.,\nonumber\\
  &\left.
    -\dfrac{1}{\zeta}\int_X (\cdot)w_iE_i(x)\dd\mu(x)
    \right)
\end{align}
\begin{align}
    \label{eq:supp:sim_part_b}
  \nabla_{\overline{\beta}}\mathcal{N}
  =&\left(
     -\dfrac{1}{\zeta}\int_X (\cdot)w_iE_i(x)\dd\mu(x)\right.,\nonumber\\
    &\left.
    -\dfrac{\delta_{ij}}{\zeta}\int_Xw_iE_i(x)\dd\mu(x)
    \right)
\end{align}

Consider the boundary conditions:
$\sum_{i=1}^{m}w_iE_i(x)\dd\mu(x)=\sum_{i=1}^{m}\dd\pi_i(x)=\mu(x)$
and
$\int_Xw_iE_i(x)\dd\mu(x)=\int_X\pi_i(x)=w_i$,
the $\nabla_{\alpha,\overline{\beta}}\mathcal{N}$ as the Hessian
form of $\mathcal{L}$, can be written as
\begin{align}
  \nabla_{\alpha,\overline{\beta}}\mathcal{N}
  =-\dfrac{1}{\zeta}\left[
     \begin{array}{cc}
       \left\langle-,\cdot\right\rangle & \left\langle\pi_j, \cdot\right\rangle\\
       \left\langle-,\pi_i\right\rangle & w_i\delta_{ij}
     \end{array}
\right]
\end{align}
with $\left\langle\phi_1,\phi_2\right\rangle=\int_X\phi_1(x)\phi_2(x)\mu(x)$,
or further
\begin{align}
  \label{eq:supp:hessian}
  \nabla_{\alpha,\overline{\beta}}\mathcal{N}
  =-\dfrac{1}{\zeta}\left[
     \begin{array}{cc}
       \dd\mu(x)\dd\mu(x')\delta(x,x') & \dd\pi_j(x)\\
       \dd\pi_i(x') & w_i\delta_{ij}
     \end{array}
\right]
\end{align}
over the basis $\{\delta(x), \mathbf{e}_i\}_{x\in X, i<m}$.

By the inverse of $\nabla_{\alpha,\overline{\beta}}\mathcal{N}$
we mean the element in $\mathrm{Hom}^b_{\mathbb{R}}\left(
  (L^\infty(X))^\vee\times(\mathbb{R}^{m-1})^\vee,
  L^\infty(X)\times\mathbb{R}^{m-1}
\right)$ which compose with $\nabla_{\alpha,\overline{\beta}}\mathcal{N}$
(on left and on right) are identities.
By the natural identity between double dual $V^{\vee\vee}\cong V$ and
the Tensor-Hom adjunction,

\begin{align}
  &\mathrm{Hom}^b_{\mathbb{R}}\left(
  (L^\infty(X))^\vee\!\times\!(\mathbb{R}^{m-1})^\vee,
  L^\infty(X)\times\mathbb{R}^{m-1}
  \right)\nonumber\\
  \cong&
         \mathrm{Hom}^b_{\mathbb{R}}\left(
         (L^\infty(X))^\vee\!\times\!(\mathbb{R}^{m-1})^\vee,
         (L^\infty(X)\!\times\!\mathbb{R}^{m-1})^{\vee\vee}
         \right)\nonumber\\
  \cong&\mathrm{Hom}^b_{\mathbb{R}}\left(
         ((L^\infty(X))^\vee\times(\mathbb{R}^{m-1})^\vee)^{\otimes2},
         \mathbb{R}
         \right), 
  \end{align}
we can write the inverse of $\nabla_{\alpha,\overline{\beta}}\mathcal{N}$
as a bilinear form again.

Denote $\nabla_{\alpha,\overline{\beta}}\mathcal{N}$
in the block form $\left[
\begin{array}{cc}
  A & B\\
  B^T & D
\end{array}\right]$.
According to the block-inverse formula
\begin{equation}\small
  \label{eq:supp:block_inverse_formula}
\left[
\begin{array}{cc}
  A & B\\
  B^T & D
\end{array}\right]^{-1}
\!\!=\! \left[\!
  \begin{array}{cc}
    A^{-1}\!\!+\!A^{-1}BF^{-1}B^TA^{-1} & -A^{-1}BF^{-1}\\
    -F^{-1}B^TA^{-1}             & F^{-1}
  \end{array}
\!\right]
\end{equation}
where $F=D-B^TA^{-1}B$ whose invertibility
determines the invertibility of $\left[
\begin{array}{cc}
  A & B\\
  B^T & D
\end{array}\right]$.

Consider that $A^{-1}\in \mathrm{Hom}^b_{\mathbb{R}}((L^\infty(X)^\vee)^{\otimes2},\mathbb{R})$,
explicitly, $A^{-1}(x,y) = \delta(x,y)$. Therefore,
from Eq.~\eqref{eq:supp:hessian},
\begin{align}
 &F_{ij}\nonumber\\
 =&\delta_{ij}w_i-\int_{X^2}\delta(x,x')\pi_i(x)\pi_j(x')\nonumber\\
 =&\delta_{ij}w_i-\int_Xw_iw_jE_i(x)E_j(x)\mu(x).
\end{align}
The matrix $F$ is symmetric, of rank $m-1$ and strictly diagonally dominant, therefore
it is invertible.
To see the strictly diagonal dominance, consider
$\sum_{j=1}^{m}\int_Xw_iw_jE_i(x)E_j(x)\mu(x)=\int_Xw_iE_i(x)\sum_{j=1}^{m}w_jE_j(x)\mu(x)
=\int_Xw_iE_i(x)\mu(x)=w_i$ by applying the marginal conditions, 
and that the matrix $F$ is of size $(m-1)\times(m-1)$ (there is no $i=m$ or $j=m$ for $F_{ij}$),
then the matrix $F$ is strictly diagonally dominant.

With all ingradients known in formula \eqref{eq:supp:block_inverse_formula},
we can calculate the inverse of $\nabla_{\alpha,\overline{\beta}}\mathcal{N}$.

Following the implicit function theorem, we need $\nabla_{w_i,y_i}\mathcal{N}$,
each partial derivative is an element in $L^{\infty}(X)^\vee\times\mathbb{R}^{m-1}$.

\begin{align}
  &\dfrac{\partial\mathcal{N}}{\partial w_i}\nonumber\\
  =&\left(-\int_XE_i(x)(\cdot)\dd\mu(x),
     \delta_{ij}\left(1-\int_XE_i(x)\dd\mu(x)\right)\right)\nonumber\\
  =&\left(-\int_X(\cdot)E_i(x)\dd\mu(x),0\right).
\end{align}
Note that if we apply the constraint $\sum_{i=1}^{m}w_i=1$ to $w_i$'s,
we may set $w_m=1-\sum_{i=1}^{m-1}w_i$ and recalculate the above derivatives
as $\nabla_{w_i'}\mathcal{N}=\nabla_{w_i}\mathcal{N}-\nabla_{w_m}\mathcal{N}$ when
$i\ne m$, and $\nabla_{w_m'}\mathcal{N}=\sum_{i=1}^{m-1}\nabla_{w_i}\mathcal{N}$.
\begin{align}
  \nabla_{y_i}\mathcal{N}
  =&\left(\dfrac{1}{\zeta}\int_X(\cdot)\nabla_{y_i}g(x,y_i)w_iE_i(x)\dd\mu(x)\right.\nonumber\\
   &\left.\dfrac{\delta_{ij}}{\zeta}\int_{X}\nabla_{y_i}g(x,y_i)w_iE_i(x)\dd\mu(x)
     \right)
\end{align}

Finally, by implicit function theorem,
$\nabla_{w_i,y_j}(\alpha^\ast,\overline{\beta}^\ast)=-\left(\nabla_{\alpha,\overline{\beta}}\mathcal{N}|_{\alpha^\ast,\overline{\beta}^\ast}\right)^{-1}\left(\nabla_{w_i,y_j}\mathcal{N}|_{\alpha^\ast,\overline{\beta}^\ast}\right)$ which fits in \eqref{eq:supp:gradient_W_w} and \eqref{eq:supp:gradient_W_y}.

\subsection{Second Derivatives}
\label{subsec:second_derivatives}

In this part, we calculate the second derivatives
of $W_{k,\zeta}^k(\mu,\mu_m)$ with respect to
the ingredients of $\mu_m$, i.e., $w_i$'s
and $y_i$'s, for the potential of applying
Newton's method in EDOT (while we have not implemented yet).

Using the previous results, we can further calculate the second derivatives of $W_{k,\zeta}^k$
about $w_i$'s and $y_i$'s. Differentiating \eqref{eq:supp:gradient_W_w} and
\eqref{eq:supp:gradient_W_y} results in
\begin{align}
  \label{eq:supp:diff2nd_W_ww}
  &\dfrac{\partial^2\mathcal{L}}{\partial w_i\partial w_j}\nonumber\\
  =&\dfrac{1}{\zeta}\int_Xg(x,y_j)\sum_{k=i,j}
     \left(\dfrac{\partial\alpha(x)}{\partial w_k}+\dfrac{\partial\overline{\beta}_k}{\partial w_k}\right)E_k(x)\dd\mu(x)\nonumber\\
  &+\dfrac{1}{\zeta}\int_X\sum_{k=1}^{n}\left(\dfrac{\partial^2\alpha(x)}{\partial w_i\partial w_j}+\dfrac{\partial^2\overline{\beta}_k}{\partial w_i\partial w_j}\right)w_kE_k(x)\dd\mu(x)\nonumber\\
  &+\dfrac{1}{\zeta^2}\int_X\sum_{k=1}^{m}\prod_{l=i,j}\left(\dfrac{\partial\alpha(x)}{\partial w_l}+
    \dfrac{\partial \overline{\beta}_k}{\partial w_l}\right)w_kE_k(x)\dd\mu(x)
\end{align}

\begin{align}
  \label{eq:supp:diff2nd_W_wy}
  &\nabla_{y_j}\dfrac{\partial\mathcal{L}}{\partial w_i}\nonumber\\
  =&\delta_{ij}\left[\int_X\left(1-\dfrac{g(x,y_i)}{\zeta}\right)\nabla_{y_i}
     E_i(x)\dd\mu(x)\right]\nonumber\\
  &+\dfrac{1}{\zeta}\int_X\nabla_{y_j}g(x,y_j)
    \left(\dfrac{\partial\alpha(x)}{\partial w_i}+\dfrac{\partial\overline{\beta}_j}{\partial w_i}\right)
    w_jE_j(x)\dd\mu(x)\nonumber\\
  &+\dfrac{1}{\zeta}\int_Xg(x,y_j)
    \nabla_{y_j}\left(\dfrac{\partial\alpha(x)}{\partial w_i}+\dfrac{\partial\overline{\beta}_j}{\partial w_i}\right)
    w_jE_j(x)\dd\mu(x)\nonumber\\
  &+\dfrac{1}{\zeta^2}\int_X\sum_{k=1}^{m}g(x,y_k)\left(\dfrac{\partial\alpha(x)}{\partial w_i}+\dfrac{\partial\overline{\beta}_k}{\partial w_i}\right)\nonumber\\
   &\ \ \left(\nabla_{y_j}\alpha(x)+\nabla_{y_j}\overline{\beta}_k\right)w_kE_k(x)\dd\mu(x)
\end{align}

\begin{align}
  \label{eq:supp:diff2nd_W_yy}
  &\nabla_{y_j}\nabla_{y_i}\mathcal{L}\nonumber\\
  =&\delta_{ij}\left[
     \int_{X}\nabla_{y_i}^2g(x,y_i)(1-\dfrac{g(x,y_i)}{\zeta})w_iE_i(x)\dd\mu(x)\right.\nonumber\\
  &\left.-\dfrac{1}{\zeta}\int_X\left(\nabla_{y_i}g(x,y_i)\right)^2(2-\dfrac{g(x,y_i)}{\zeta})w_iE_i(x)\dd\mu(x)
    \right]\nonumber\\
  &+\dfrac{1}{\zeta}\int_X\left(1-\dfrac{g(x,y_i)}{\zeta}\right)\nabla_{y_i}g(x,y_i)\cdot\nonumber\\
  &\quad \left(\nabla_{y_j}\alpha(x)+\nabla_{y_j}\overline{\beta_i}\right)w_iE_i(x)\dd\mu(x)\nonumber\\
  &+\dfrac{1}{\zeta}\int_X\left(1-\dfrac{g(x,y_j)}{\zeta}\right)\nabla_{y_j}g(x,y_j)\cdot\nonumber\\
  &\quad \left(\nabla_{y_i}\alpha(x)+\nabla_{y_i}\overline{\beta_j}\right)w_jE_j(x)\dd\mu(x)\nonumber\\
  &+\dfrac{1}{\zeta}\int_X\sum_{k=1}^{m}g(x,y_k)\left(\nabla_{y_i}\nabla_{y_j}\alpha(x)+\nabla_{y_i}\nabla_{y_j}\overline{\beta}_k\right)\cdot \nonumber\\
  &\quad w_kE_k(x)\dd\mu(x)\nonumber\\
  &+\dfrac{1}{\zeta^2}\int_X\sum_{k=1}^{m}g(x,y_k)\prod_{l=i,j}\left(\nabla_l\alpha+\nabla_l\overline{\beta}_k\right)\cdot\nonumber\\
  &\quad w_kE_k(x)\dd\mu(x)
\end{align}

Once we have the second derivatives of $g(x,y)$ on $y_i$'s, we need the second derivatives
of $\alpha^\ast$ and $\overline{\beta}^\ast$ to build the above second derivatives.
From the formula
$\nabla_{w_i,y_j}(\alpha^\ast,\overline{\beta}^\ast)=-\left(\nabla_{\alpha,\overline{\beta}}\mathcal{N}|_{\alpha^\ast,\overline{\beta}^\ast}\right)^{-1}\left(\nabla_{w_i,y_j}\mathcal{N}|_{\alpha^\ast,\overline{\beta}^\ast}\right)$, we can differentiate

\begin{align}
  &\nabla_{w_k,y_l}\nabla_{w_i,y_j}(\alpha^\ast,\overline{\beta}^\ast)\nonumber\\
  =&-\nabla_{w_k,y_l}\left(\nabla_{\alpha,\overline{\beta}}\mathcal{N}|_{\alpha^\ast,\overline{\beta}^\ast}\right)^{-1}\left(\nabla_{w_i,y_j}\mathcal{N}|_{\alpha^\ast,\overline{\beta}^\ast}\right)\nonumber\\
   &-\left(\nabla_{\alpha,\overline{\beta}}\mathcal{N}|_{\alpha^\ast,\overline{\beta}^\ast}\right)^{-1}\left(\nabla_{w_k,y_l}\nabla_{w_i,y_j}\mathcal{N}|_{\alpha^\ast,\overline{\beta}^\ast}\right).
\end{align}

Here from the formula that $\nabla A^{-1}=-A^{-1}\nabla AA^{-1}$ (this is the product rule for
$AA^{-1}=I$), we have

\begin{align}
  &\nabla_{w_k,y_l}\left(\nabla_{\alpha,\overline{\beta}}\mathcal{N}|_{\alpha^\ast,\overline{\beta}^\ast}\right)^{-1}\nonumber\\
  =&-\left(\nabla_{\alpha,\overline{\beta}}\mathcal{N}|_{\alpha^\ast,\overline{\beta}^\ast}\right)^{-1}
     \nabla_{w_k,y_l}\left(\nabla_{\alpha,\overline{\beta}}\mathcal{N}|_{\alpha^\ast,\overline{\beta}^\ast}\right)\nonumber\\
  &\quad \left(\nabla_{\alpha,\overline{\beta}}\mathcal{N}|_{\alpha^\ast,\overline{\beta}^\ast}\right)^{-1}
\end{align}
and
\begin{align}
  &\nabla_{w_k}\left(\nabla_{\alpha,\overline{\beta}}\mathcal{N}|_{\alpha^\ast,\overline{\beta}^\ast}\right)\nonumber\\
  =&-\dfrac{1}{\zeta}\left[
  \begin{array}{cc}
    0 & \delta_{jk}E_k(x)\dd\mu(x)\\
    \delta_{ik}E_k(x')\dd\mu(x') & \delta_{ij}\delta_{jk}
  \end{array}
\right]
\end{align}
\begin{align}
  &\nabla_{y_k}\left(\nabla_{\alpha,\overline{\beta}}\mathcal{N}|_{\alpha^\ast,\overline{\beta}^\ast}\right)\nonumber\\
  =&\dfrac{1}{\zeta^2}\left[
  \begin{array}{cc}
    0 & \delta_{jk}\nabla_{y_k}g(x,y_k)\dd\pi_k(x)\\
    \delta_{ik}\nabla_{y_k}g(x',y_k)\dd\pi_k(x') & 0
  \end{array}
\right].
\end{align}

The last piece we need is 
$\left(\nabla_{w_k,y_l}\nabla_{w_i,y_j}\mathcal{N}|_{\alpha^\ast,\overline{\beta}^\ast}\right)$:

\begin{align}
  \dfrac{\partial^2\mathcal{N}}{\partial w_j\partial w_i}=(0,0),
\end{align}

\begin{align}
  &\nabla_{y_j}\dfrac{\partial\mathcal{N}}{\partial w_i}\nonumber\\
  =&\left(\dfrac{1}{\zeta}\int_X(\cdot)\nabla_{y_j}g(x,y_j)E_i(x)\dd\mu(x),0\right),
\end{align}

\begin{align}
  &\nabla_{y_j}\nabla_{y_i}\mathcal{N}\nonumber\\
  =&\dfrac{\delta_{ij}}{\zeta}\left(
     \int_X(\cdot)\nabla_{y_i}^2g(x,y_i)w_iE_i(x)\dd\mu(x)
     \right.\nonumber\\
  &\quad+\int_X(\cdot)\left(\nabla_{y_i}g(x,y_i)\right)^2w_iE_i(x)\dd\mu(x),\nonumber\\
  & \delta_{ik}\int_X\nabla_{y_i}^2g(x,y_i)w_iE_i(x)\dd\mu(x)\nonumber\\
  &\quad+\left.\int_X\left(\nabla_{y_i}g(x,y_i)\right)^2w_iE_i(x)\dd\mu(x)
     \right),
\end{align}
where in the last one, $k$ represent the $k$-th component in $\mathcal{N}$'s second part
(about $\overline{\beta}$).

%% file: SI_3_empirical.tex
\section{Empirical Parts}
\label{sec:supp_empirical}

\subsection{Estimate $W_{k,\zeta}^k$: Richardson Extrapolation and Others}
\label{subsec:supp_richardson}

In the analysis, we may need $W_{k,\zeta}^k(\mu, \mu_m)$ to compare how 
discretization methods behave. However, when the $\mu$ is not discrete,
generally we are not able to obtain the analytical solution to the Wasserstein
distance.

In certain cases including all examples this paper contains, the Wasserstein can
be estimated by finite samples (with a large size).
According to \cite{mensch2020online}, for $\mu\in\mathcal{P}(X)$ in our setup (a
probability measure on a compact Polish space with Borel algebra) and
$g=d_X^k\in\mathcal{C}(X^2)$ being a continuous function, the the Online
Sinkhorn methods can be used to estimate $W_{k,\zeta}^k$. Online Sinkhorn needs
a large number of samples for $\mu$ (in batch) to be accurate.

In our paper, as $X$ are compact subsets in $\mathbb{R}^n$ and $\mu$ has a
continuous probability density function, 
we may use Richardson Extrapolation method to estimate the Wasserstein distance
between $\mu$ and $\mu_m$, 
which may require fewer samples and fewer computations (Sinkhorn twice with
different sizes).

Our examples are on intervals or rectangles, in which two grids $\Lambda_1$ of
$N$ points and $\Lambda_2$ of $rN$ points ($N$ and $rN$ are both integers) can
be constructed naturally for each. With $\mu$ determined by a smooth probability
density function $\rho$, let $\mu_{(N)}$ be normalization of
$\sum_{i=1}^{N}\rho(\Lambda_i)\delta_{\Lambda_i}$ (this may not be a probability
distribution, so we use its normalization). From continuity of $\rho$ and the
boundedness of dual variables $\alpha, \beta$, we can conclude that
$$\lim_{N\rightarrow\infty}W_{k,\zeta}^k(\mu_{(N)},\mu_m)=W_{k,\zeta}^k(\mu,\mu_m).$$
Let $W_{k,\zeta}^k(\mu_{(N)},\mu_m)$ be a function of $1/N$, to apply Richardson
extrapolation, we need the exponent of lowest term of $1/N$ in the expansion
$W_{k,\zeta}^k(\mu_{(N)},\mu_m)=W^\ast +
\mathcal{O}(1/N^h)+\mathcal{O}((1/N)^{h+1})$, where
$W^\ast = W_{k,\zeta}^k(\mu,\mu_m)$. 

\begin{figure}[ht]
    \centering
    \includegraphics[scale=0.55]{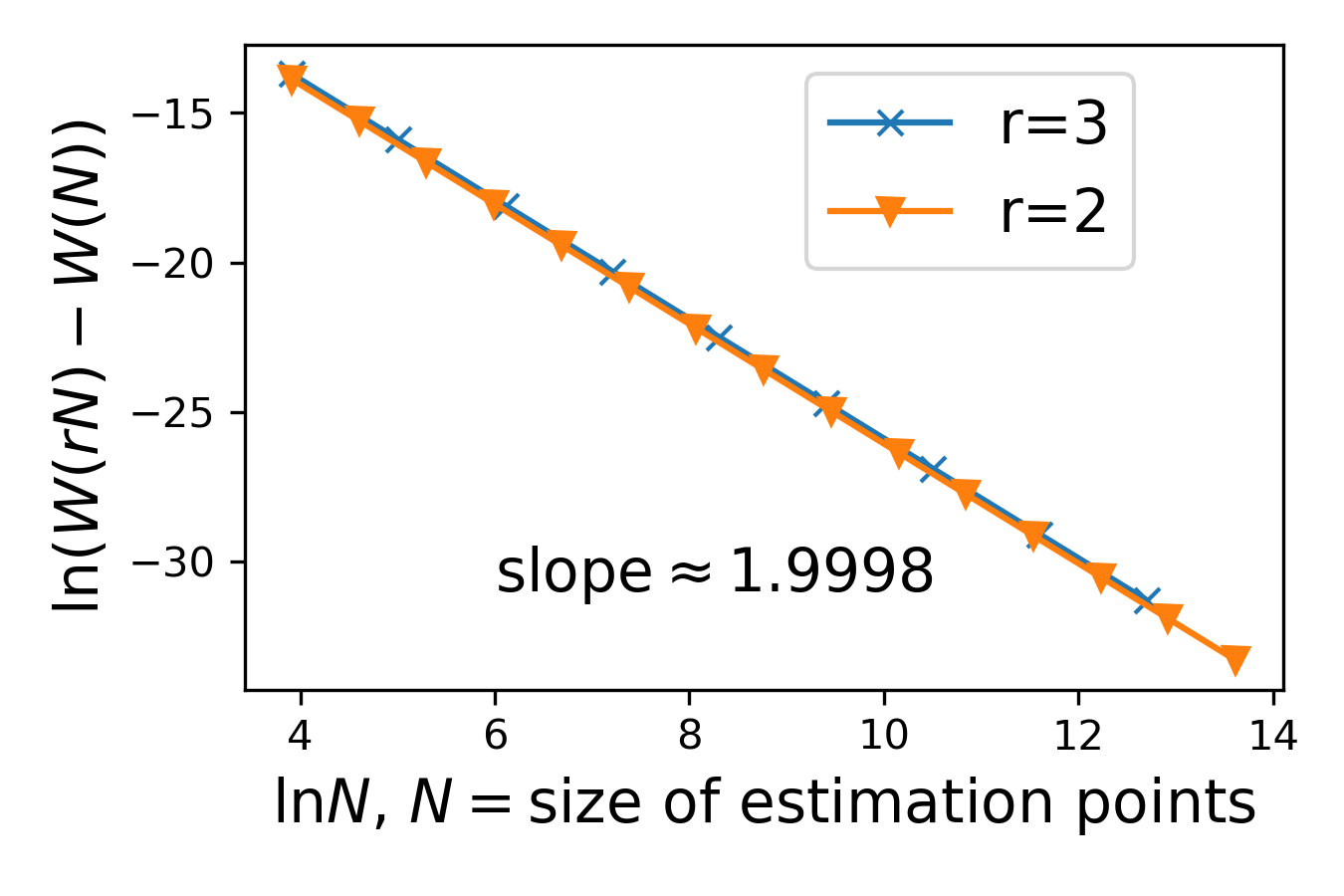}
    \caption{Richardson Extrapolation: the power of the expansion
    about $N^{-1}$. We take the EDOT $\mu_5$ of example 2 (1-dim 
    truncated normal mixture) as the target $\mu_m$, use 
    evenly positioned $\mu_N$ for different $N$'s to estimate.
    The $y$-axis is
    $\ln(W_{2,0.01}^2(\mu_{(rN)},\mu_5)-W_{2,0.01}^2(\mu_{(N)},\mu_5))$,
    where $r=2$ and $r=3$ are calculated.
    With $\ln N$ as $x$-axis, linearity can be observed. The slopes are both 
    about $-1.9998$, representing the exponent of the leading
    non-constant term of $W_{2,0.01}^2(\mu_{(N)},\mu_5)$ on $N$, 
    while the theoretical result is $r=-k/d=-2$. The differences are 
    from higher order terms on $N$.}
    \label{fig:supp:richardson}
\end{figure}

Consider that $$\left|W_{k,\zeta}^k(\mu,\mu_m)-
W_{k,\zeta}^k(\mu_{(N)},\mu_m)\right|\le W_{k,\zeta}^k(\mu,\mu_{(N)})$$
Since $W_{k,\zeta}(\mu_{(N)},\mu)\propto N^{-1/d}$, we may conclude that
$h=k/d$, where $d$ is the dimension of $X$. Figure~\ref{fig:supp:richardson} shows an
empirical example in $d=1$, $k=2$ situation.

\subsection{Example: The Sphere}
\label{subsec:supp_sphere}

The CW-complex structrue of the unit sphere $S^2$ is constructed as:
let $(1,0,0)$, the point on the equator be the only dimension-0
structure, and let let the equator be the dimension-1 structure
(line segment $[0,1]$ attached to the dimension-0 structure by identifying
both end points to the only point $(1,0,0)$).
The dimension-2 structure is the union of two unit discs, identified to the
south / north hemisphere of $S^2$ by stereographic projection
\begin{equation}\small
  \label{eq:stereo}
\pi_{\pm N}:(X,Y)\rightarrow\dfrac{1}{1+X^2+Y^2}\left({2X}, {2Y},\pm({X^2+Y^2-1})\right)
\end{equation}
with respect to the north / south pole.

\paragraph{Spherical Geometry.}
The spherical geometry is the Riemannian manifold structure induced by the
embedding onto the unit sphere in $\mathbb{R}^3$.

The geodesic between two points is the shorter arc along the great circle
determined by the two points. In their $\mathbb{R}^3$ coordinates,
$d_{S^2}(\mathbf{x}, \mathbf{y}) = \arccos(\left\langle\mathbf{x}, \mathbf{y}\right\rangle)$.
Composed with stereographic projections, the distance in terms of CW-complex
coordinates can be calculated (and be differentiated).

The gradient about $\mathbf{y}$ (or its CW-coordinate) can be calculated via
above formulas. In practice, the only
problem is that when $\mathbf{x}=\pm\mathbf{y}$ function $\arccos$
at $\pm1$ is singular. 
From the symmetry of sphere on the rotation along axis $\mathbf{x}$, the derivatives
of distance along all directions are the same. Therefore, we may choose the
radial direction on the CW-coordinate (unit disc). And the differentiations are primary
to calculate.

\subsection{A Note on the Implementation of SGD with Momentum}
There is a slight difference between our implementation 
of the SGD and the Algorithm provided in the paper.
In the implementation, we give two different learning
rates to the positions ($y_i$'s) and the weights ($w_i$'s),
as moving along positions is usually observed much slower than
moving along weights. Empirically, we make the learning
rates on position be exactly 3 times of the learning rates
on weights, at each SGD iteration. With this change,
the convergence is faster, but we do not have a theory
or empirical evidences to show a fixed ratio 3 is the 
best choice.

Implementing and testing the Newton's method (hybrid with SGD)
and other improved SGD methods could be good problems
to work on.

\subsection{Some Figures from Empirical Results}
In this part, we post some figures revealing the results 
from the simulations, mainly the 3d figures in different directions.
\begin{figure*}[ht]
    \centering
    \includegraphics[scale=0.4]{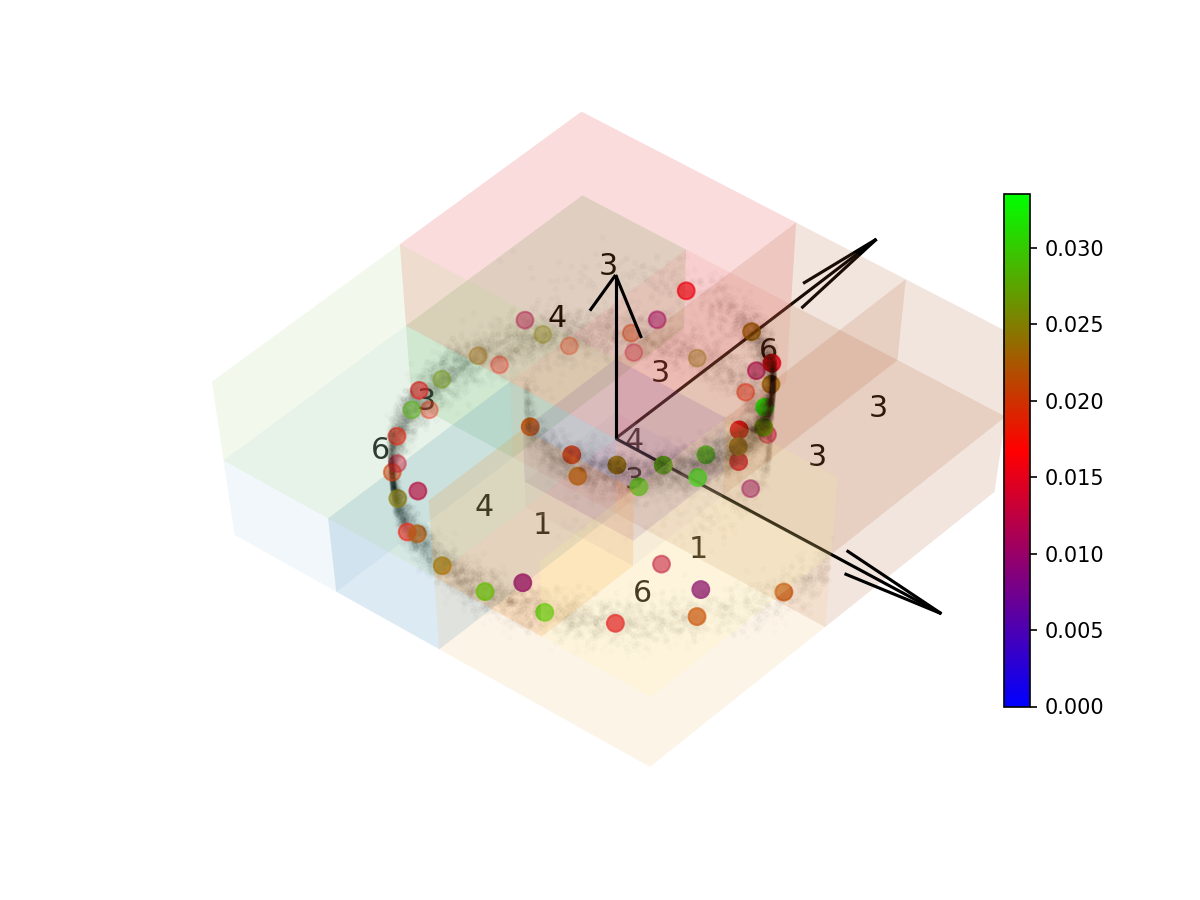}
    \includegraphics[scale=0.4]{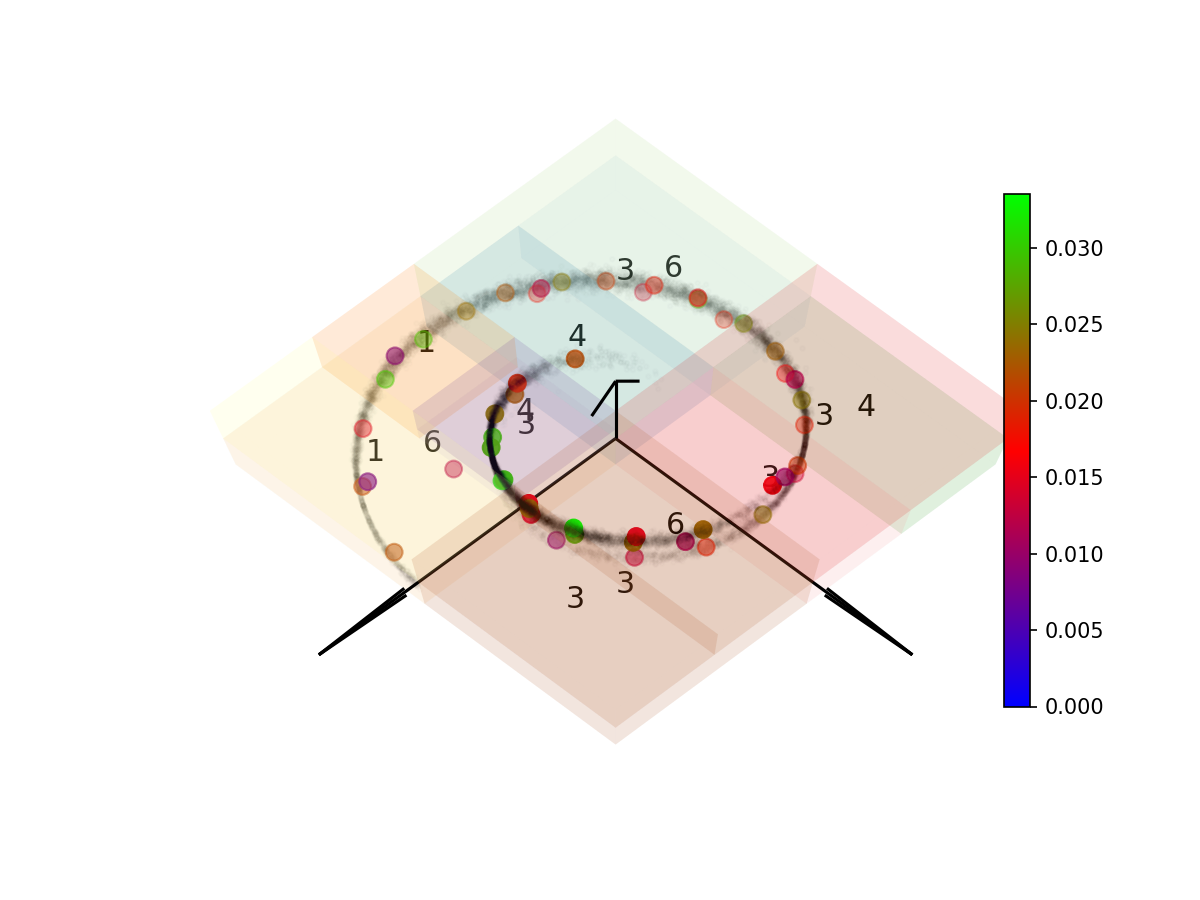}
    \caption{The Swiss roll example, 3D discretization (same as the paper) in different view directions.}
    \label{fig:supp:swiss_3d}
\end{figure*}

\begin{figure*}[ht]
    \centering
    \includegraphics[scale=0.4]{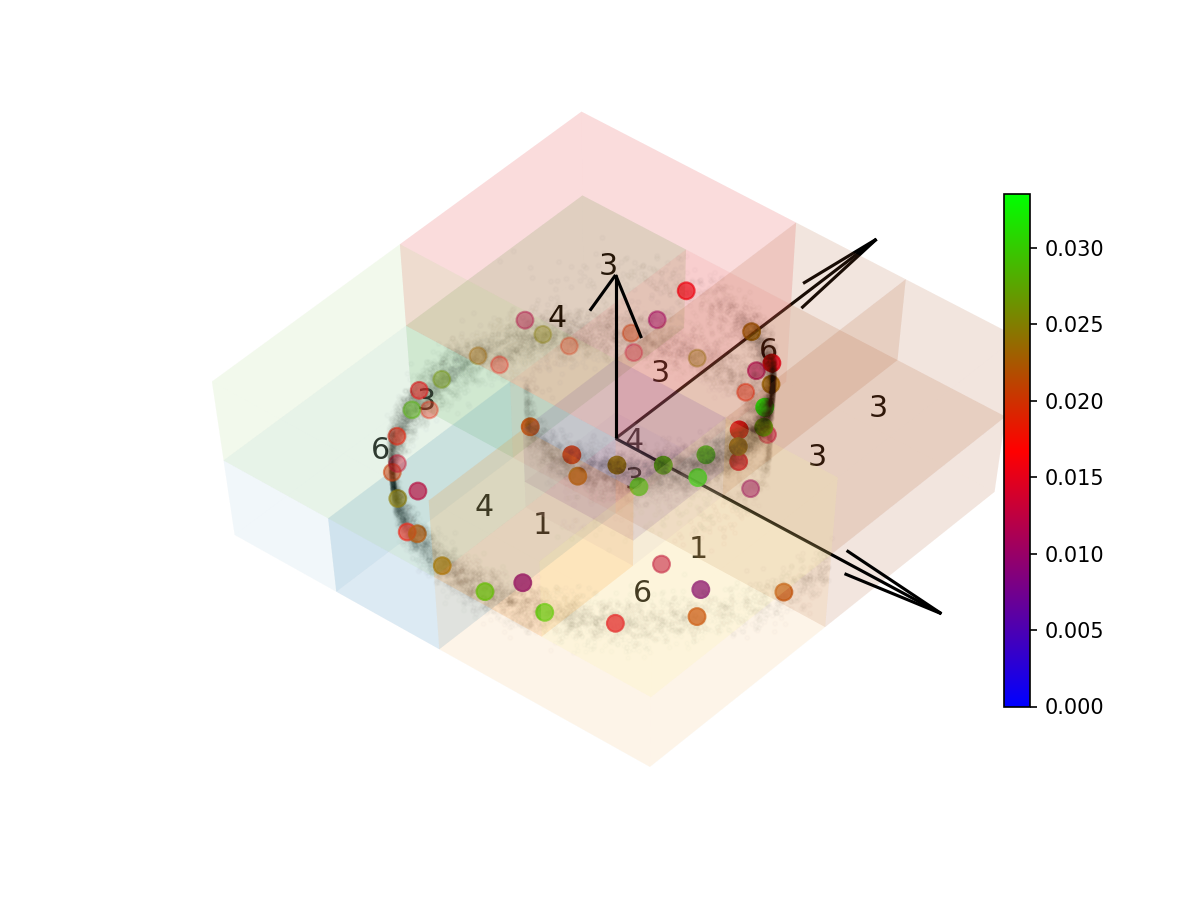}
    \includegraphics[scale=0.4, trim={100 0 0 0},
    clip]{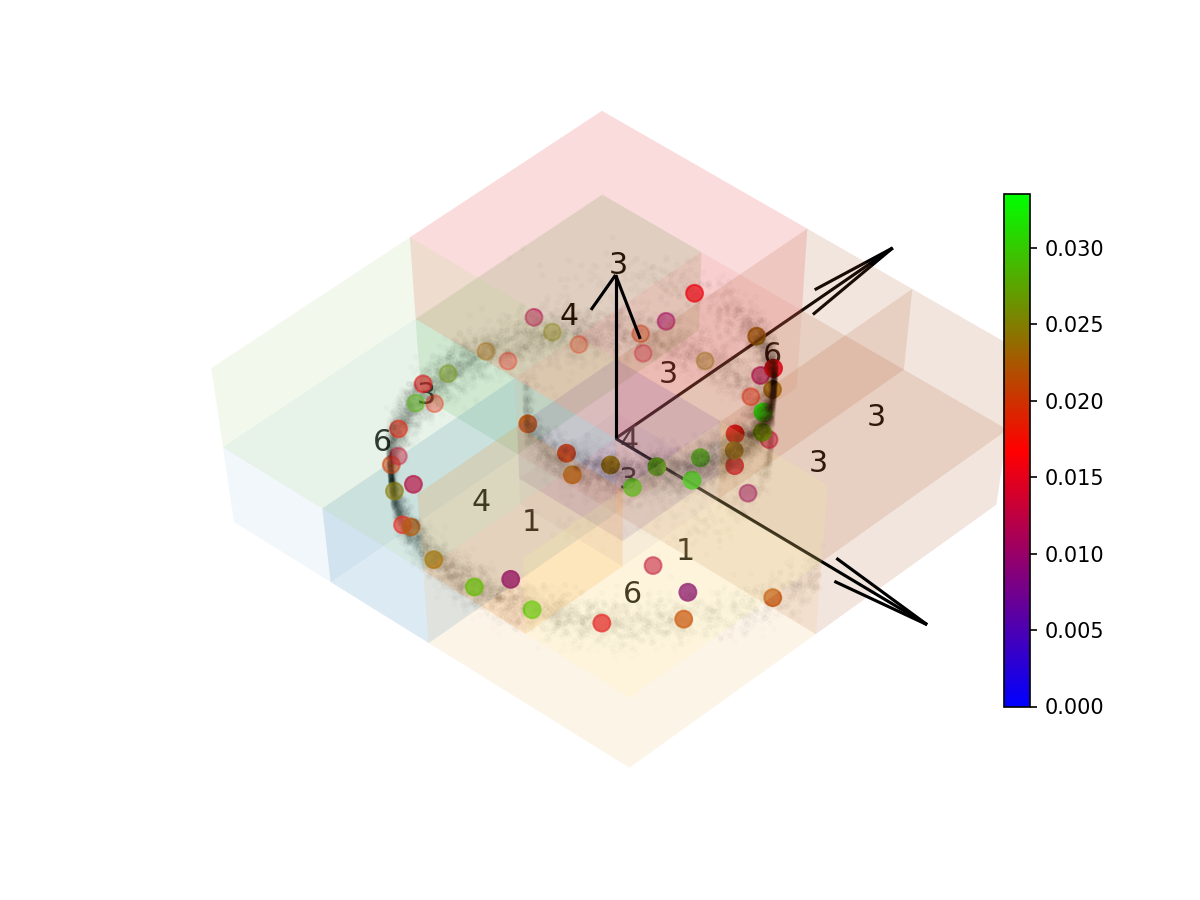}
    \caption{The Swiss roll example, 3D discretization (same as the paper) in terms of an autostereogram.}
    \label{fig:supp:swiss_3d_autostereogram}
\end{figure*}

\begin{figure*}[ht]
    \centering
    \includegraphics[scale=0.45]{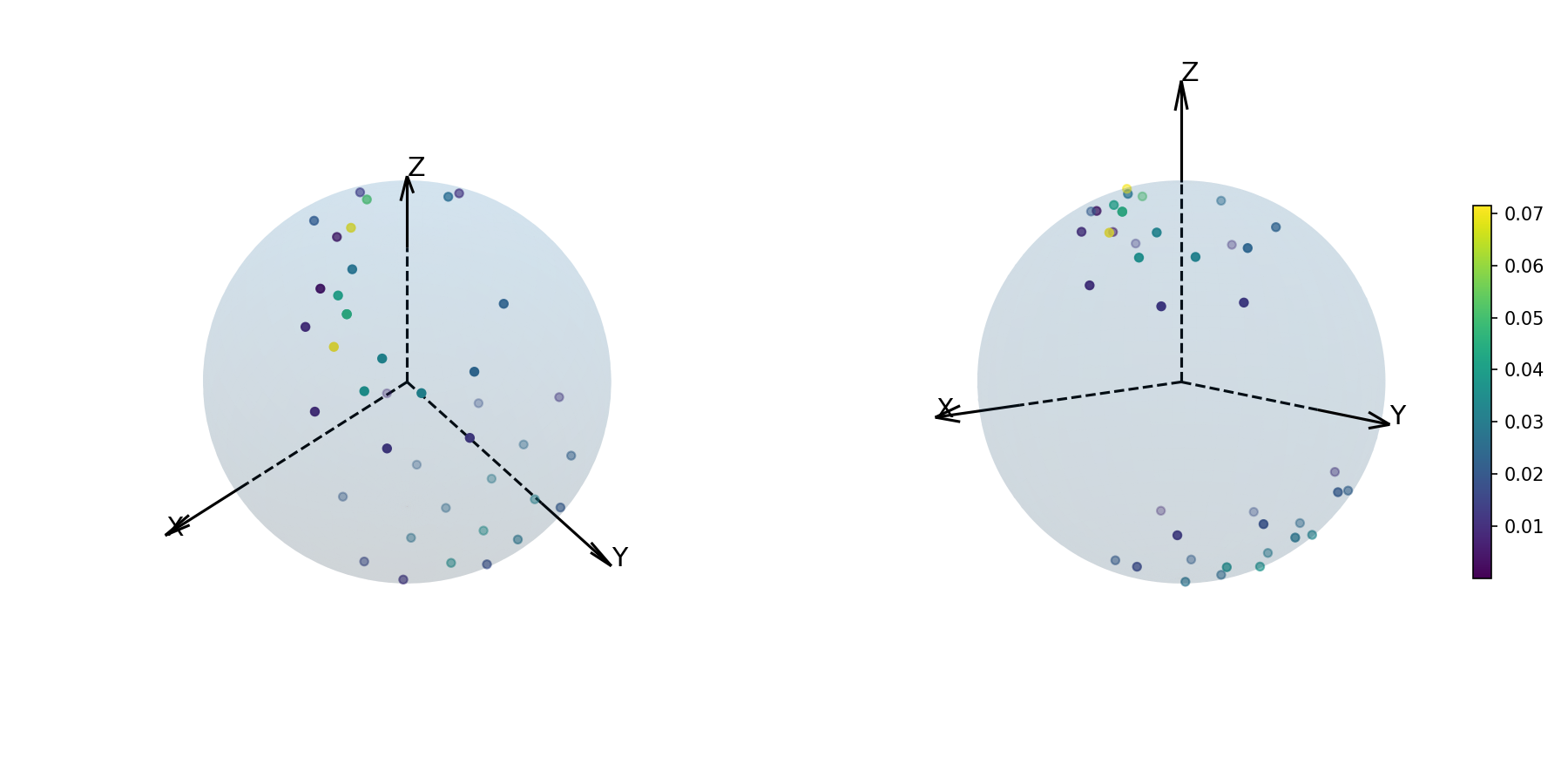}
    \caption{The Sphere example, 3D discretization (same as the paper) on two view directions.}
    \label{fig:supp:spherical_3d}
\end{figure*}

\begin{figure*}[ht]
    \centering
    \includegraphics[scale=0.45]{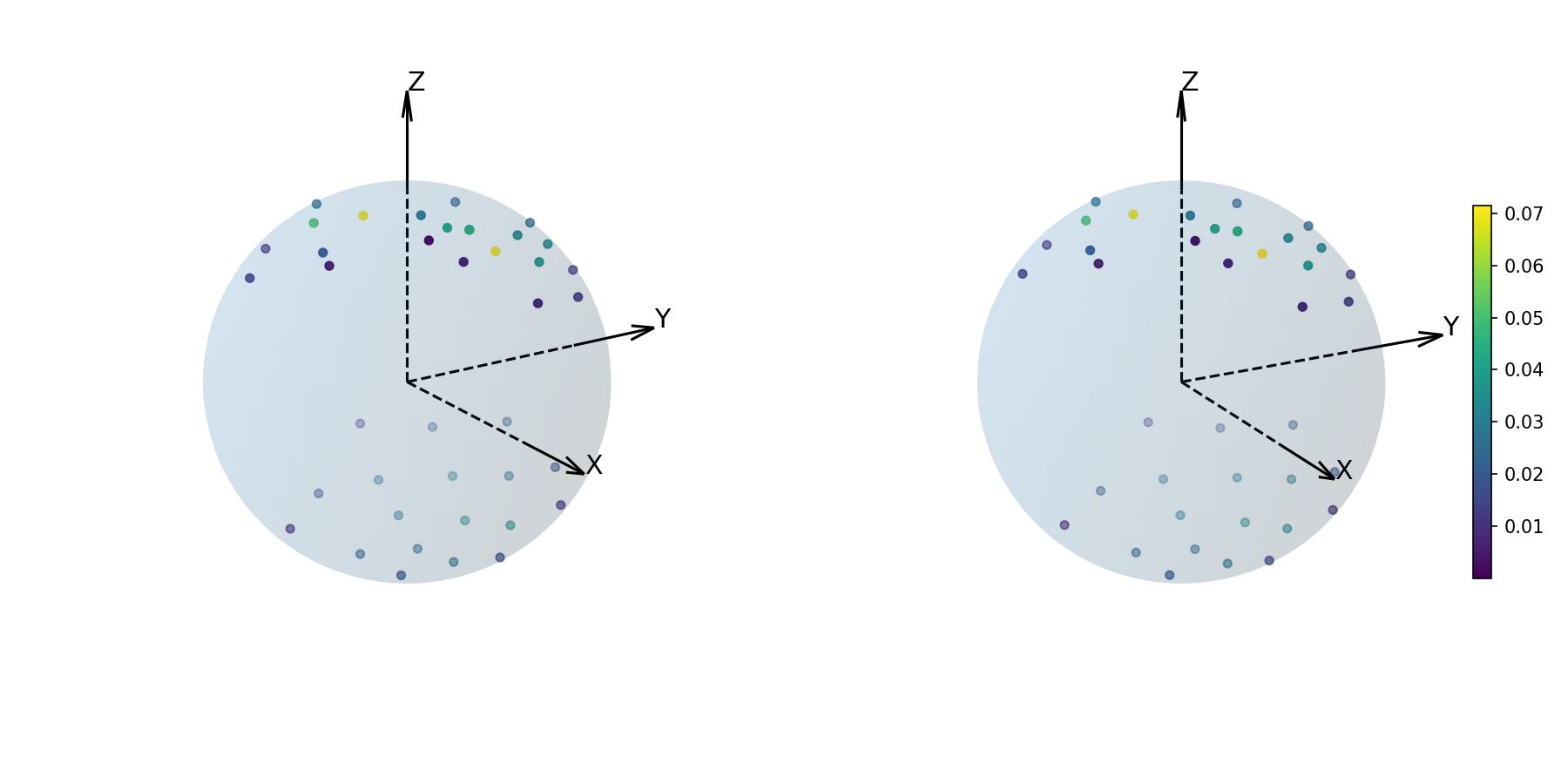}
    \caption{The Sphere example, 3D discretization (same as the paper) in terms of an autostereogram.}
    \label{fig:supp:spherical_3d_autostereogram}
\end{figure*}

\begin{figure*}
    \centering
    \foreach \i in {1,2,...,11}{
    \includegraphics[scale=0.5]{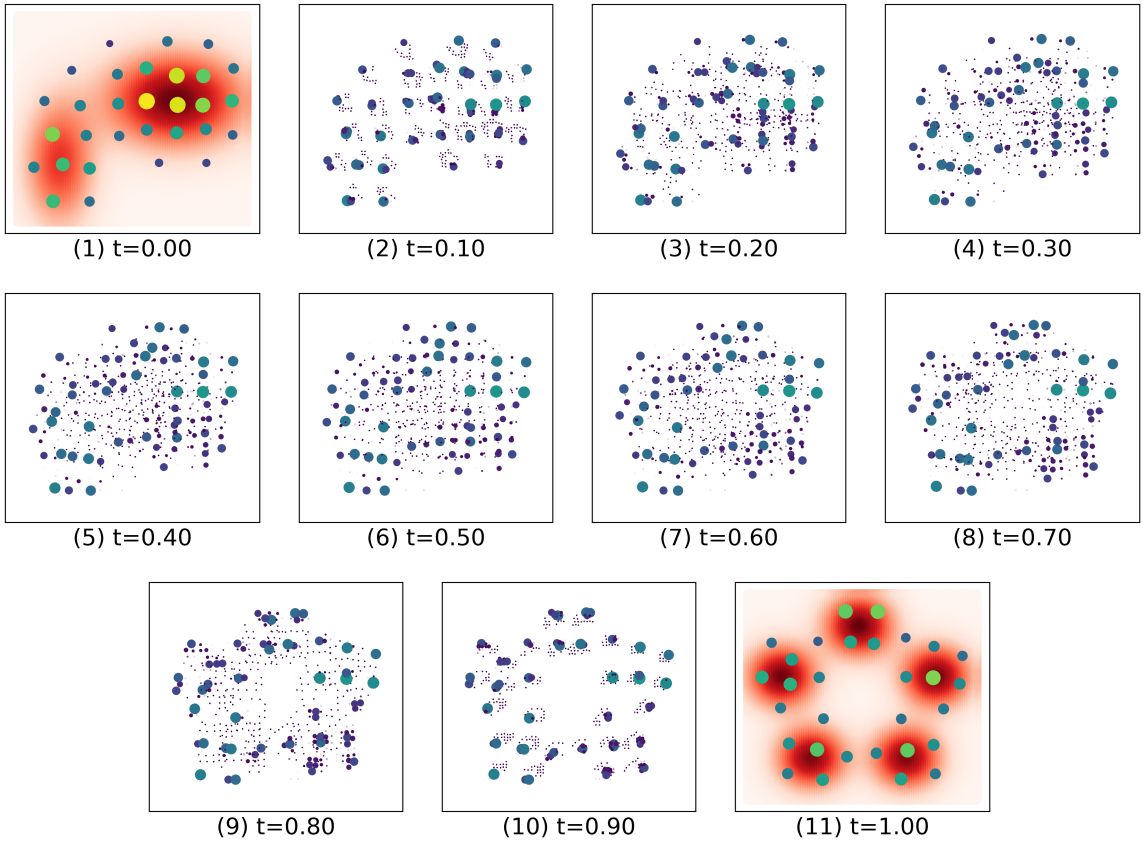}}
    \caption{The McCann Interpolation figures in finer time resolution, for
    visualizing the transference plan from (1) to (11).
    It is a refined figure of the original paper. We see can see that the 
    larger bubbles (representing a large probability mass) moves in a 
    short distance and smaller pieces moves longer.}
    \label{fig:supp:mccann}
\end{figure*}